\documentclass[11pt]{amsart}
\usepackage{mathrsfs}
\usepackage{geometry}
\usepackage{graphicx}
\usepackage{xcolor} 
\usepackage{amssymb}
\usepackage{epstopdf}
\usepackage{amsmath,amscd}
\usepackage{amsthm}
\usepackage{enumerate}
\usepackage{url,verbatim}
\usepackage{subfig}
\usepackage{enumitem}
\usepackage{tikz}

\newcommand{\compliancerevision}[1]{{\color{black}#1}}
\newenvironment{introrevision}{\begingroup\color{black}}{\endgroup}

\calclayout

\RequirePackage[colorlinks,citecolor=blue,urlcolor=blue]{hyperref}
\hypersetup{
  pdftitle={Asymptotics of Bounded Lecture-Hall Tableaux},
  pdfauthor={Zhongyang Li, David Keating, Istv\'an Prause},
  pdfkeywords={lecture-hall tableaux, random tilings, dimer models,
    limit shapes, complex Burgers equation, Gaussian free field,
    Schur generating functions}
}
\theoremstyle{plain}
\newcommand\LHLL[2]{
  \def\X{#1} \def\Y{#2}
  \foreach \i in {0,...,\X}
  {
\pgfmathsetmacro{\m}{\Y-1/(\i+1)};
    \draw[gray,very thin] (\i,0) -- (\i,\m);
    \node at (\i,-.3) {\i};
  }
  \foreach \j in {0,...,\Y}
  {
    \node at (-.3,\j) {\j};
  }
\pgfmathsetmacro{\m}{\Y-1};
  \foreach \x in {1,...,\X}
{ \foreach \j in {0,...,\m}
  {
\draw[gray,very thin] (\x-1,\j) -- (\x,\j);
}}
 \foreach \i in {2,...,\X}
 \foreach \j in {0,...,\m} 
{\foreach \y in {2,...,\i}
{
\pgfmathsetmacro{\w}{\j+(\y-1)/\i};
\pgfmathsetmacro{\z}{\j+(\y-1)/(\i+1)};
\draw[gray,very thin] (\i-1,\w) -- (\i,\z);
}}
\node at (\X+1,1) {$\cdots$};
\node at (\X+1,\Y-1) {$\cdots$};
}

\newcommand\BLHLL[2]{
  \def\X{#1} \def\Y{#2}
  \foreach \i in {0,...,\X}{
    \pgfmathsetmacro{\m}{\Y-1/(\i+1)}
    \draw[gray,thin] (\i,0)--(\i,\m);
  }
  \pgfmathsetmacro{\m}{\Y-1}
  \foreach \x in {1,...,\X}{
    \foreach \j in {0,...,\m}{
      \draw[gray,thin] (\x-1,\j)--(\x,\j);
    }
  }
  \foreach \i in {2,...,\X}{
    \foreach \j in {0,...,\m}{
      \foreach \y in {2,...,\i}{
        \pgfmathsetmacro{\w}{\j+(\y-1)/\i}
        \pgfmathsetmacro{\z}{\j+(\y-1)/(\i+1)}
        \draw[gray,thin] (\i-1,\w)--(\i,\z);
      }
    }
  }
}

\newcommand\DLL[2]{
  \def\X{#1} \def\Y{#2}
  \pgfmathsetmacro{\m}{\Y-1}
  \foreach \x in {1,...,\X}{
    \foreach \j in {0,...,\m}{
      \node at (\x,\j) {$\circ$};
      \node at (\x-.1,\j+.1) {$\bullet$};
      \draw[gray,thin] (\x-1,\j)--(\x-.1,\j+.1);
    }
  }
  \foreach \j in {0,...,\m}{
    \node at (0,\j) {$\circ$};
    \node at (-.1,\j+.1) {$\bullet$};
    \draw[gray,thin] (0,\j)--(-.1,\j+.1);
  }
  \foreach \i in {1,...,\X}{
    \foreach \j in {0,...,\m}{
      \foreach \y in {1,...,\i}{
        \pgfmathsetmacro{\w}{\j+\y/(\i+1)}
        \node at (\i,\w) {$\circ$};
        \node at (\i-.1,\w+.1) {$\bullet$};
        \draw[gray,thin] (\i,\w)--(\i-.1,\w+.1);
      }
    }
  }
  \pgfmathsetmacro{\v}{\m-1}
  \foreach \i in {1,...,\X}{
    \foreach \j in {0,...,\v}{
      \foreach \y in {0,...,\i}{
        \pgfmathsetmacro{\w}{\j+\y/(\i+1)}
        \pgfmathsetmacro{\q}{\j+(\y+1)/(\i+1)}
        \draw[gray,thin] (\i-.1,\w+.1)--(\i,\q);
      }
    }
  }
  \foreach \i in {1,...,\X}{
    \foreach \j in {\m,...,\m}{
      \foreach \y in {1,...,\i}{
        \pgfmathsetmacro{\w}{\j+(\y-1)/(\i+1)}
        \pgfmathsetmacro{\q}{\j+\y/(\i+1)}
        \draw[gray,thin] (\i-.1,\w+.1)--(\i,\q);
      }
    }
  }
  \foreach \i in {0,...,0}{
    \foreach \j in {1,...,\m}{
      \draw[gray,thin] (\i-.1,\j-1+.1)--(\i,\j);
    }
  }
  \foreach \i in {2,...,\X}{
    \foreach \j in {0,...,\m}{
      \foreach \y in {2,...,\i}{
        \pgfmathsetmacro{\w}{\j+(\y-1)/\i}
        \pgfmathsetmacro{\z}{\j+(\y-1)/(\i+1)}
        \draw[gray,thin] (\i-1,\w)--(\i-.1,\z+.1);
      }
    }
  }
}

\newtheorem{theorem}{Theorem}
\newtheorem{definition}[theorem]{Definition}
\newtheorem{lemma}[theorem]{Lemma}
\newtheorem{proposition}[theorem]{Proposition}
\newtheorem{corollary}[theorem]{Corollary}
\newtheorem{example}[theorem]{Example}
\newtheorem{assumption}[theorem]{Assumption}
\newtheorem{remark}[theorem]{Remark}

\newcommand\ol{\overline}

\newcommand\RR{{\mathbb R}}

\newcommand\NN{{\mathbb N}}

\newcommand\HH{{\mathbb H}}

\newcommand\YY{{\mathbb {Y}}}

\newcommand\si{\sigma}

\newcommand\q{\quad}

\renewcommand\ell{l}

\newcommand\bm{\mathbf{m}}

\newcounter{mycount}

\numberwithin{equation}{section}
\numberwithin{theorem}{section}
\numberwithin{figure}{section}

\begin{document}
\title{Asymptotics of Bounded Lecture-Hall Tableaux}
\author{Zhongyang Li}
\address{Department of Mathematics, University of Connecticut,
341 Mansfield Road, Unit 1009, Storrs, CT 06269-1009, USA}
\email{zhongyang.li@uconn.edu}
\thanks{Corresponding author: Zhongyang Li.}
\author{David Keating}
\address{Department of Mathematics, University of Massachusetts Boston,
100 William T. Morrissey Boulevard, Boston, MA 02125, USA}
\email{David.Keating@umb.edu}
\author{Istv\'an Prause}
\address{Mathematics, Faculty of Science and Engineering,
{\AA}bo Akademi University, 20500 Turku, Finland}
\email{istvan.prause@abo.fi}



\date{}

\begin{abstract}
We study weighted bounded \(n\)-lecture hall tableaux with bound \(t_n\)
in the joint limit \(n,t_n\to\infty\) with
\(t_n/n\to\alpha\in(0,\infty)\). Using their non-intersecting-path
representation, we construct a Schur generating function adapted to this
model and derive moment formulas for the limiting counting measures.
When the normalized tail sums of the weights have a
\(C^2\) macroscopic profile with negative derivative, the rescaled
height functions converge in probability in
\(L^1_{\mathrm{loc}}(\mathbb R\times(0,\alpha))\) to a deterministic
limit shape whose complex slope
satisfies a
Burgers equation with an explicit source term.
In the uniformly weighted model studied by Corteel,
Keating, and Nicoletti, the source term vanishes, proving their
Conjecture~6.1. Under an additional interval-structure
assumption on the limiting bottom boundary, and only continuity and strict
monotonicity of the tail profile, the unrescaled height fluctuations converge to the Gaussian free field
in the sense of joint moments of horizontal polynomial observables. The method applies although the
particle configurations are not Gelfand--Tsetlin patterns and the associated
dimer model is not doubly periodic.
\end{abstract}

\maketitle

\section{Introduction} 

\begin{introrevision}
Random tilings, tableaux, and non-intersecting paths form a meeting
point of algebraic combinatorics, representation theory, probability,
statistical mechanics, complex analysis, and nonlinear partial
differential equations. A central manifestation of these connections
is the complex Burgers equation: Kenyon and Okounkov showed that it
governs limit shapes in random dimer models \cite{ko07}. At the
fluctuation scale, many exactly solvable random-surface models exhibit
the Gaussian free field \cite{RK01,Li13,bg16}. Two mechanisms have
been particularly successful in turning exact formulas into
asymptotic results: Gelfand--Tsetlin interlacing and Schur-process
methods on the one hand, and the spectral theory of doubly periodic
dimers on the other.

This paper establishes this limit-shape/Burgers/Gaussian-field picture
for weighted bounded lecture-hall tableaux. Although the finite model
admits exact algebraic formulas, the two standard approaches described
above are not directly available: its horizontal slices do not satisfy
the usual interlacing relations, and its associated dimer graph is not
doubly periodic. We develop a Schur generating function tailored to
the lecture-hall branching rule and use it to analyze both the
macroscopic shape and the global fluctuations.

Our results reveal a clear separation between these two scales. At
the macroscopic level, the inhomogeneous weights leave a visible
signature on the limit shape by producing an explicit forcing term in
the complex Burgers equation. At the fluctuation level, we prove a
general multilevel central limit theorem and, under an additional
condition on the bottom boundary, identify its covariance with the
Dirichlet Green kernel after a natural change of coordinates. Thus the
inhomogeneity changes the deterministic geometry of the model, while
the centered fluctuations retain the Gaussian-free-field structure
familiar from more classical random surfaces.

\subsection*{The model and its geometric realization}

Lecture hall tableaux were introduced in \cite{SK20}; they generalize
lecture hall partitions \cite{BE971,BE972} and anti-lecture hall
compositions \cite{SC03}, and recover reverse semistandard Young
tableaux in an appropriate limit. Andrews later gave a proof of the
lecture hall partition theorem using MacMahon's partition analysis
\cite{Andrews98LHP}; see also the survey \cite{SA16}.
\end{introrevision}

\begin{introrevision}
Let $\mu\subset\lambda$ be partitions and let $\lambda/\mu$ denote
the corresponding skew Young diagram. A tableau of shape
$\lambda/\mu$ is a map from its cells to $\NN$, the set of
nonnegative integers. We write $l(\lambda)$ for the length of
$\lambda$, $\YY$ for the set of all partitions, and
$[n]=\{1,\ldots,n\}$.
\end{introrevision}

\begin{definition}Let $n$ be a \compliancerevision{positive }integer and $\mu\subset \lambda$ be partitions such that $l(\lambda)\leq n$.
An $n$-lecture hall tableau of shape $\lambda/\mu$ is a tableau $L$ of shape $\lambda/\mu$ satisfying the following conditions
\begin{align*}
\frac{L(i,j)}{n+c(i,j)}\geq \frac{L(i,j+1)}{n+c(i,j+1)},\qquad
\frac{L(i,j)}{n+c(i,j)}> \frac{L(i+1,j)}{n+c(i+1,j)},
\end{align*}
where $c(i,j)=j-i$ is the content of the cell $(i,j)$. The set of all $n$-lecture hall tableaux is denoted by $LHT_n(\lambda/\mu)$. We write $\lfloor x\rfloor$ for the integer part of a real number $x$. For $L\in LHT_n(\lambda/\mu)$, let $\lfloor L\rfloor$ be the tableau of shape $\lambda/\mu$ whose $(i,j)$th entry is 
   $\left\lfloor\frac{L(i,j)}{n-i+j}\right\rfloor$ .
\end{definition}
See the left panel of Figure~\ref{fig:tpd} for an example of a lecture hall tableau.

\begin{introrevision}
We impose, for every cell $(i,j)$, the additional finite-height condition
\[
    L(i,j)<t(n+j-i),
\]
and say that the tableau is bounded by the positive integer $t$. Bounded lecture hall
tableaux were enumerated in \cite{CK20}. Throughout this paper we
consider the joint scaling regime
\[
    n,t_n\longrightarrow\infty,
    \qquad
    \frac{t_n}{n}\longrightarrow\alpha\in(0,\infty),
\]
so that the number of paths and the height of the associated graph
remain comparable.

The path representation is important both probabilistically and
analytically: it turns a bounded lecture hall tableau into a
non-intersecting path ensemble and, subsequently, into a dimer
configuration. For a prescribed shape $\lambda^{(0)}$ and nonnegative weights
$\mathbf x=(x_0,\ldots,x_{t-1})$, we assign
\[
    \operatorname{wt}_{\mathbf x}(L)
    =
    \prod_{(i,j)\in\lambda^{(0)}}
    x_{\left\lfloor L(i,j)/(n+j-i)\right\rfloor},
\]
and, when the total tableau weight is positive, normalize these weights to obtain the
probability measure studied below. Equal positive weights give the uniformly sampled model of
\cite{SKN21}; more generally, the weights may vary with both the
system size $n$ and the vertical level. We recall the bijection
of \cite{CK20} before
describing the probability measure and its asymptotic behavior.
\end{introrevision}

\begin{definition}\label{dflhg}
\begin{enumerate}
\item Given a positive integer $t$, the lecture hall graph is a graph $\mathcal{G}_t=(V_t,E_t)$. This graph can be described through an embedding in the plane with vertex set $V_t$ given by
\begin{itemize}
    \item $\left(i,\frac{j}{i+1}\right)$ for $i\geq 0$ and $0\leq j<t(i+1)$.
\end{itemize}
and the directed edges given by
\begin{itemize}
    \item from $\left(i,k+\frac{r}{i+1}\right)$
    to $\left(i+1,k+\frac{r}{i+2}\right)$ for $i\geq 0$, $0\leq r\leq i$ and $0\leq k<t$
    \item from $\left(i,k + \frac{r + 1}{i + 1}\right)$ to $\left(i,k + \frac{r}{i + 1}\right)$ for $i\geq 0$ and $0\leq r \leq i$ and
$0 \leq k < t-1$ or for $i\geq 0$ and $0\leq r <i$ and $k = t-1$.
\end{itemize}
\item Given a positive integer $t$ and a partition $\lambda=(\lambda_1,\lambda_2,\ldots,\lambda_n)$ with $\lambda_1\geq \lambda_2\geq\ldots\geq \lambda_n\geq 0$, a non-intersecting path configuration is a system of $n$ paths on the graph $\mathcal{G}_t$. For each integer $i$ satisfying $1\leq i\leq n$, the $i$th path starts at $\left(n-i,t-\frac{1}{n-i+1}\right)$, ends at $(n-i+\lambda_i,0)$ and moves only downwards and rightwards. The paths are said to be non-intersecting if they do not share a vertex.
\item\label{item:finite-horizontal-truncation}
{
Given positive integers $t$ and $N$, let $\mathcal G_{t,N}$ denote the
finite subgraph of $\mathcal G_t$ induced by the vertices whose horizontal
coordinates lie in $[0,N]$.
}
\end{enumerate}
\end{definition}

See the middle graph of Figure \ref{fig:tpd} for an example of $\mathcal{G}_3$ and a configuration of non-intersecting paths on $\mathcal{G}_3$.

\begin{theorem}[\cite{CK20}]
\label{thm:path-bijection}
There is a bijection between bounded lecture hall tableaux of shape
$\lambda$ and bound $t$, and non-intersecting path configurations on
$\mathcal G_t$ whose $i$th path starts at
$\left(n-i,t-\frac{1}{n-i+1}\right)$ and ends at
$(n-i+\lambda_i,0)$, for $i=1,2,\ldots,n$.

More precisely, the path configuration corresponding to a lecture hall
tableau of shape $\lambda$ contains exactly $|\lambda|$ non-vertical
edges. Their left endpoints are
\[
    \left(n+j-i-1,\frac{L(i,j)}{n+j-i}\right),
    \qquad (i,j)\in\lambda.
\]
The corresponding non-intersecting path configuration is obtained uniquely
by adding vertical edges so as to join the prescribed starting and ending
points.
\end{theorem}

For an $n$-lecture hall tableau bounded by $t$, the parameter $t$ is the height of the lecture hall graph $\mathcal G_t$, while $n$ is the number of paths in the corresponding non-intersecting path configuration. See Figure~\ref{fig:tpd}.

\begin{figure}
    \centering
    \begin{tikzpicture}[scale=.9,baseline=(current bounding box.center)]
        \draw (0,0) rectangle (2,2);
        \draw (0,1)--(2,1);
        \draw (1,0)--(1,2);
        \node at (.5,1.5) {$5$};
        \node at (1.5,1.5) {$6$};
        \node at (.5,.5) {$2$};
        \node at (1.5,.5) {$3$};
        \node[font=\small] at (1,2.65) {(a) Lecture-hall tableau};
        \node[font=\scriptsize] at (.5,2.18) {$j=1$};
        \node[font=\scriptsize] at (1.5,2.18) {$j=2$};
        \node[font=\scriptsize,left] at (-.05,1.5) {$i=1$};
        \node[font=\scriptsize,left] at (-.05,.5) {$i=2$};
    \end{tikzpicture}
    \hfill
    \begin{tikzpicture}[scale=.9,baseline=(current bounding box.center)]
        \BLHLL{3}{3}
        \draw[red,very thick,line cap=round,line join=round]
            (0,2)--(1,2)--(1,3/2)--(2,4/3)--(2,0);
        \draw[red,very thick,line cap=round,line join=round]
            (1,5/2)--(2,7/3)--(2,2)--(3,2)--(3,0);
        \foreach \i in {0,...,3}{
            \node[below] at (\i,0) {$\i$};
        }
        \foreach \j in {0,...,3}{
            \node[left] at (0,\j) {$\j$};
        }
        \node[font=\small] at (1.5,3.25) {(b) Paths on $\mathcal G_3$};
        \node[font=\scriptsize,below] at (1.5,-.4) {horizontal coordinate};
        \node[font=\scriptsize,rotate=90] at (-.55,1.5) {height};
    \end{tikzpicture}
    \hfill
    \begin{tikzpicture}[scale=.9,baseline=(current bounding box.center)]
        \DLL{3}{3}
        \tikzset{matching/.style={red,very thick,line cap=round}}
        \draw[matching] (0,2)--(.9,2.1);
        \draw[matching] (1,2)--(.9,1.6);
        \draw[matching] (1,1.5)--(1.9,4/3+.1);
        \draw[matching] (2,4/3)--(1.9,1.1);
        \draw[matching] (2,1)--(1.9,2/3+.1);
        \draw[matching] (2,2/3)--(1.9,1/3+.1);
        \draw[matching] (2,1/3)--(1.9,.1);
        \draw[matching] (1,5/2)--(1.9,7/3+.1);
        \draw[matching] (2,7/3)--(1.9,2+.1);
        \draw[matching] (2,2)--(2.9,2+.1);
        \draw[matching] (3,2)--(2.9,7/4+.1);
        \draw[matching] (3,7/4)--(2.9,6/4+.1);
        \draw[matching] (3,6/4)--(2.9,5/4+.1);
        \draw[matching] (3,5/4)--(2.9,1+.1);
        \draw[matching] (3,1)--(2.9,3/4+.1);
        \draw[matching] (3,3/4)--(2.9,2/4+.1);
        \draw[matching] (3,2/4)--(2.9,1/4+.1);
        \draw[matching] (3,1/4)--(2.9,.1);
        \draw[matching] (0,1)--(-.1,1.1);
        \draw[matching] (0,0)--(-.1,.1);
        \draw[matching] (1,0)--(.9,.1);
        \draw[matching] (1,.5)--(.9,.6);
        \draw[matching] (1,1)--(.9,1.1);
        \draw[matching] (2,5/3)--(1.9,5/3+.1);
        \draw[matching] (2,8/3)--(1.9,8/3+.1);
        \draw[matching] (3,2+1/4)--(2.9,2+1/4+.1);
        \draw[matching] (3,2+2/4)--(2.9,2+2/4+.1);
        \draw[matching] (3,2+3/4)--(2.9,2+3/4+.1);
        \node at (0,3) {$\circ$};
        \node at (1,3) {$\circ$};
        \node at (1.9,-1/3+.1) {$\bullet$};
        \node at (2.9,-1/4+.1) {$\bullet$};
        \draw[matching] (-.1,2.1)--(0,3);
        \draw[matching] (.9,2.6)--(1,3);
        \draw[matching] (1.9,-1/3+.1)--(2,0);
        \draw[matching] (2.9,-1/4+.1)--(3,0);
        \node[font=\small] at (1.5,3.25) {(c) Dimer configuration};
    \end{tikzpicture}
    \caption{Tableau, non-intersecting paths, and dimers. The left panel shows a lecture hall tableau
$L$ of shape $\lambda=(2,2)$ with $n=2$ and entries
$L(1,1)=5$, $L(1,2)=6$, $L(2,1)=2$, and $L(2,2)=3$; the row and column
coordinates are denoted by $(i,j)$. The normalized entries are
$5/2$, $2$, $2$, and $3/2$, respectively, and the tableau is bounded by
$t=3$. The middle panel shows the corresponding non-intersecting path
configuration, with the vertical coordinate recording height, and the
right panel shows the associated dimer configuration on a graph that is
not doubly periodic.}
    \label{fig:tpd}
\end{figure}

\begin{introrevision}
\subsection*{Previous work and structural obstacles}

The uniformly weighted model was studied in \cite{SKN21} using a
heuristic application of the tangent method. That analysis predicted
the frozen boundary but did not determine the full limit shape, and
Conjecture~6.1 of \cite{SKN21} asserted that the slopes of the limiting
height function satisfy the inviscid complex Burgers equation. The
weighted model considered here contains the uniformly sampled model
as a special case and allows the vertical inhomogeneity to vary on a
macroscopic scale.

The two geometric realizations above explain why this problem requires
a different approach. For doubly periodic dimers, the complex Burgers
equation is tied to spectral curves and a variational principle
\cite{ko07}; see also \cite{ADPZ23}. Our argument does not use a
variational principle, and the periodic spectral-curve theory of
\cite{ko07} does not apply directly because the lecture-hall dimer
graph is not doubly periodic. On the algebraic side, Schur generating
functions have been used to study random signatures and uniform dimers
on hexagonal and square-grid domains \cite{bg,bg16,bk}, with extensions
to nonuniform models on rail-yard graphs
\cite{BL17,ZL18,ZL201,ZL20,Li21}. Here the intermediate partitions fail
the usual interlacing inequalities, as shown explicitly in
Figure~\ref{fig:lambdakappaEx}, and the skew lecture-hall weights are
not skew Schur weights. Hence the model cannot be treated simply as an
ordinary Schur process or Gelfand--Tsetlin array.

\subsection*{Main results and significance}

The algebraic starting point is the generating function in
Definition~\ref{df01}. Lemma~\ref{l03} expresses it as a ratio of
ordinary Schur polynomials, thereby making Schur differential
operators available despite the absence of interlacing. This identity
is used in both the law-of-large-numbers and multilevel fluctuation
arguments. Our main conclusions are as follows.

\begin{enumerate}[label=\textup{(\roman*)},leftmargin=*]

\item
\emph{Limit shape and spectral description.}
For sequences of bottom boundaries whose rescaled shapes converge to
a fixed macroscopic profile, and for nonnegative weights that may vary
with both the system size and the vertical level,
Theorem~\ref{t15} proves  convergence of the counting
measure along any sequence of macroscopic horizontal slices for which
the normalized tail-weight ratio converges to a value in $(0,1)$, and
gives the limiting moments explicitly. Proposition~\ref{prop:stieltjes-root}
characterizes its Stieltjes transform by a distinguished root of a
spectral equation. Under the $C^2$ tail-profile hypothesis, these
slice limits assemble into a deterministic height function, and
Proposition~\ref{prop:two-dimensional-height-lln} gives convergence in
probability in
\[
    L^1_{\mathrm{loc}}\bigl(\mathbb R\times(0,\alpha)\bigr).
\]
Under the assumptions of Theorem~\ref{t17}, every interior endpoint of
the liquid region on a horizontal slice corresponds to a multiple root
of the equation determining the complex slope.

\item
\emph{Inhomogeneous Burgers dynamics.}
Let $s(y)$ be the limiting normalized tail sum of the weights.
Under the hypotheses of Theorem~\ref{thm:m31}, that
theorem constructs a complex slope
$u(\chi,y)\in\mathbb H$, recovers both derivatives of the limiting
height function from $u$, and proves
\[
    u_y+u\,u_\chi
    =
    \frac{s''(y)}{s'(y)}u.
\]
Thus the dependence of the forcing term on the weights is compressed
into the single macroscopic profile $s$; The chosen bottom boundary determines the initial data for the equation
and hence affects the resulting limit shape.  For
uniform weights,
\[
    s(y)=1-\frac{y}{\alpha},
\]
so the source term vanishes. Corollary~\ref{cor:uniform-burgers}
therefore proves Conjecture~6.1 of \cite{SKN21}.

\item
\emph{Multilevel Gaussian fluctuations and the GFF.}
Theorem~\ref{t31} proves a multilevel central limit theorem in moments
for the power-sum observables, with an explicit double-contour
covariance and joint Wick moments. This theorem applies under a
general regular bottom boundary and does not require the stronger
geometric assumption used for the GFF identification. Under
Assumption~\ref{ap32} and
conditions~\ref{ass:tail-profile-convergence}--\ref{ass:tail-profile-invertibility}
of Assumption~\ref{ass:tail-profile}, Lemma~\ref{lrd} shows that the
spectral map is a diffeomorphism from the liquid region in
$(\chi,s)$-coordinates to the upper half-plane, and
Lemma~\ref{lem:contour-to-green} transforms
the contour covariance into its Dirichlet Green kernel. Theorem~\ref{thm:gff} shows that the centered discrete height function,
multiplied by $-\sqrt{\pi}$, converges to the pullback of the
zero-boundary Gaussian free field in the sense of joint moments of
polynomially weighted integrals along finitely many horizontal slices.

\item
\emph{An explicit inhomogeneous phase diagram.}
Appendix~\ref{app:linear-bottom-boundary} treats the linear bottom
boundary
\[
    \lambda^{(0)}(n)
    =
    \bigl(2n,2(n-1),\ldots,2\bigr).
\]
The spectral equation becomes a cubic. Subject to the interior-domain
condition stated in the appendix, its discriminant determines the
liquid interval and internal frozen boundary explicitly. For $\alpha=1$,
Figure~\ref{fig:fb1} compares uniform weights, for which $s=1-y$,
with a nonlinear profile satisfying $y=(1-s)^2$. The latter is
realized, for example, by the following explicit choice of
size-dependent weights:
\[
    t_n=n,
    \qquad
    x_r^{(n)}=\sqrt{r+1}-\sqrt r,
    \qquad 0\leq r\leq n-1,
\]
because its normalized tail sum is
\[
    s_n(y)
    =1-\sqrt{\frac{\lfloor ny\rfloor}{n}}
    \longrightarrow 1-\sqrt y.
\]
For this choice, the Burgers source is
$s''(y)/s'(y)=-1/(2y)$, whereas it vanishes for uniform weights.
Both frozen-boundary curves have
\[
    \chi(z)=\frac{3z^3}{z^3+2},
    \qquad z\in(0,1)\cup(1,\infty),
\]
whereas their vertical coordinates are respectively
\[
    y(z)=\frac{3z}{z^3+2}
    \quad\text{and}\quad
    y(z)=\left(\frac{3z}{z^3+2}\right)^2.
\]
Thus this explicit specialization visibly records the weight
profile in the geometry of the random surface.

\end{enumerate}

These results connect an enumerative tableau model, a nonperiodic
random surface, symmetric-function asymptotics, a nonlinear
first-order PDE, and a log-correlated Gaussian field. Thus nonuniform weights change the deterministic limit shape and the
map to spectral coordinates, but not the Gaussian-free-field nature of
the fluctuations. Under Assumption~\ref{ap32}, their covariance in
spectral coordinates is given by the Dirichlet Green kernel. The Gaussian free field
is known to govern height fluctuations in many planar dimer models
\cite{RK01,Li13}; related GFF-type results for standard Young tableaux
appear in \cite{raposo}. The present model shows that this mechanism
can persist when both Gelfand--Tsetlin interlacing and doubly periodic
lattice geometry are absent.
\end{introrevision}
\medskip

\noindent\textbf{Relation to the conference proceedings version.}
Some principal results were announced in the twelve-page conference
proceedings article \cite{LKP24}: a one-level moment formula, a
frozen-boundary criterion, the Burgers equation for uniform weights,
and a Gaussian free field limit. The present article gives complete
proofs and substantially extends these results to nonuniform weights.
It proves convergence of the full two-dimensional height function,
establishes joint Gaussian fluctuations across multiple horizontal
slices, and explains how the weights modify the Burgers equation. It
also gives the precise assumptions and sense of convergence for the
Gaussian free field result. The uniform Burgers equation announced in
\cite{LKP24} is recovered as a special case of
Theorem~\ref{thm:m31}.

\medskip
\noindent\textbf{Organization of the paper.}
Section~\ref{lsti} constructs the lecture-hall Schur generating
function and proves the limiting moment formula
(Theorem~\ref{t15}), together with its spectral and frozen-boundary
consequences. Section~\ref{sect:be} proves the two-dimensional
law of large numbers and the weighted complex Burgers equation
(Theorem~\ref{thm:m31}). Section~\ref{gff} proves the multilevel
Gaussian central limit theorem and identifies, under
Assumption~\ref{ap32}, the limiting height observables with those of
the GFF (Theorem~\ref{thm:gff}). Appendix~\ref{sa} contains technical
Schur-asymptotic results, and
Appendix~\ref{app:linear-bottom-boundary} analyzes the explicit linear
bottom-boundary example.

\section{Limit Shape as $n,t_n\to\infty$}\label{lsti}

In this section, we prove the moment formula for the limit of the counting measure when
$n,t_n\to\infty$ and $t_n/n\to\alpha\in(0,\infty)$ by defining and analyzing a
Schur generating function adapted to lecture hall tableaux. The need for this new
generating function is that the natural particle configurations obtained from the path model
are not Gelfand--Tsetlin patterns: adjacent level partitions need not interlace, as the
explicit example in Figure \ref{fig:lambdakappaEx} shows. This is also distinct from the
doubly periodic dimer setting; the dimer correspondence is recalled in Figure
\ref{fig:tpd}, but the corresponding hexagon--octagon graph is not doubly periodic.
The main result of this section is Theorem~\ref{t15}.

The example in Figure \ref{fig:lambdakappaEx} also explains why the particle configurations
obtained from lecture hall tableaux are not Gelfand--Tsetlin schemes in disguised notation.
Indeed, the path configuration in Figure \ref{fig:lambdakappaEx} comes, through the
bijection of Theorem~\ref{thm:path-bijection}, from the following $5$-lecture hall tableau of shape
$(4,3,1)$ and bound $t=4$:
\[
L=
\begin{array}{cccc}
16 & 8 & 9 & 4\\
12 & 6 & 7\\
7
\end{array}.
\]
The lecture hall inequalities are verified by normalizing each entry by $5+j-i$:
\[
\frac{16}{5}\ge \frac{8}{6}\ge \frac{9}{7}\ge \frac{4}{8},
\qquad
\frac{12}{4}\ge \frac{6}{5}\ge \frac{7}{6},
\]
and
\[
\frac{16}{5}>\frac{12}{4}>\frac{7}{3},\qquad
\frac{8}{6}>\frac{6}{5},\qquad
\frac{9}{7}>\frac{7}{6}.
\]
Moreover, all entries satisfy $L(i,j)<4(5+j-i)$, so the tableau is bounded by $4$.
The partitions read from the horizontal levels of the corresponding path configuration are
\[
\lambda^{(3)}=(1,0,0,0,0),\quad
\lambda^{(2)}=(1,1,1,0,0),\quad
\lambda^{(1)}=(3,3,1,0,0),\quad
\lambda^{(0)}=(4,3,1,0,0).
\]
These adjacent partitions fail the usual Gelfand--Tsetlin interlacing condition; for instance,
interlacing between $\lambda^{(1)}$ and $\lambda^{(2)}$ would require
\[
\lambda^{(1)}_1\ge \lambda^{(2)}_1\ge \lambda^{(1)}_2,
\]
but here $\lambda^{(2)}_1=1<3=\lambda^{(1)}_2$. Thus the sequence of partitions in
Figure \ref{fig:lambdakappaEx} is not a Gelfand--Tsetlin pattern.

Let \(t\) be a positive integer, and let \(\mathcal G_t\) be the lecture hall
graph of height \(t\). Let \(\mathcal M\) be a non-intersecting path configuration
with \(n\) paths and bottom boundary \(\lambda^{(0)}\), meaning that the \(i\)th
path ends at
\[
    (n-i+\lambda^{(0)}_i,0), \qquad 1\le i\le n .
\]
Choose \(N\ge n-1+\lambda^{(0)}_1\), and let \(\mathcal G_{t,N}\) be the
{finite truncation defined in
Definition~\ref{dflhg}(\ref{item:finite-horizontal-truncation})}. Then all paths
of \(\mathcal M\) are contained in
\(\mathcal G_{t,N}\).

Fix an integer \(\kappa\) with \(0\le \kappa<t\). Since \(\mathcal G_{t,N}\)
has only finitely many non-vertical edges, we may choose \(\epsilon>0\) so small
that no non-vertical edge of \(\mathcal G_{t,N}\) intersects the horizontal band
\[
    \kappa < y \le \kappa+\epsilon .
\]
We associate to \(\mathcal M\) a length-\(n\) partition
\[
    \lambda^{(\kappa)}
    =
    \bigl(\lambda^{(\kappa)}_1,\ldots,\lambda^{(\kappa)}_n\bigr)
\]
as follows:
\begin{itemize}
    \item \(\lambda^{(\kappa)}_1\) is the number of absent vertical edges of
    \(\mathcal M\) intersecting the line \(y=\kappa+\epsilon\) and lying to the
    left of the rightmost vertical edge of \(\mathcal M\) present at height
    \(y=\kappa+\epsilon\).

    \item For \(2\le j\le n\), \(\lambda^{(\kappa)}_j\) is the number of absent
    vertical edges of \(\mathcal M\) intersecting the line \(y=\kappa+\epsilon\)
    and lying to the left of the \(j\)th rightmost vertical edge of
    \(\mathcal M\) present at height \(y=\kappa+\epsilon\).
\end{itemize}
See Figure \ref{fig:lambdakappaEx} for an example.

\begin{figure}
  \centering
\begin{tikzpicture}
\LHLL{8}4
\draw [red,very thick] (0,3)--(0,0);
\draw [red,very thick](1,7/2)--(1,0);
\draw [red,very thick](2,11/3)--(2,7/3)--(3,9/4)--(3,0);
\draw [red,very thick](3,15/4)--(3,12/4)--(4,15/5)--(4,13/5)--(4,6/5)--(5,7/6)--(6,8/7)--(6,0);
\draw [red,very thick] (4,19/5) -- (4,16/5) -- (5,19/6) -- (5,16/6) -- (5,8/6) -- (6,9/7) -- (7,10/8)-- (7,4/8) -- (8,4/9) -- (8,0);
\draw[dashed,fill=black,opacity=0.2] (-0.2,2) rectangle (8.2,2+1/12);
\end{tikzpicture}
\caption{Non-intersecting lattice paths on ${\mathcal G}_4$
for $n=5$.  We have $\lambda^{(3)}=(1,0,0,0,0)$,
$\lambda^{(2)}=(1,1,1,0,0)$, $\lambda^{(1)}=(3,3,1,0,0)$ and $\lambda^{(0)}=(4,3,1,0,0)$.
The sequence of partitions $(\lambda^{(0)},\lambda^{(1)},\lambda^{(2)},\lambda^{(3)})$ does not form a Gelfand--Tsetlin scheme.}\label{fig:lambdakappaEx} 
\end{figure}

We use the standard semistandard-tableau convention for Schur and skew Schur
functions; see, for example, \cite[Chapter I, Section 5]{Macdonald95}.  If
\(\mu\subseteq\lambda\) and \(\mathbf z=(z_1,z_2,\ldots)\), then
\[
    s_{\lambda/\mu}(\mathbf z)
    =
    \sum_{T\in \operatorname{SSYT}(\lambda/\mu)}
    \prod_{(i,j)\in\lambda/\mu} z_{T(i,j)},
\]
where the sum is over semistandard Young tableaux of shape \(\lambda/\mu\),
with entries weakly increasing along rows and strictly increasing down columns.
We set \(s_{\lambda/\mu}=0\) unless \(\mu\subseteq\lambda\), and write
\(s_\lambda=s_{\lambda/\varnothing}\).

We shall use the standard branching identity
\[
    s_{\lambda/\mu}(\mathbf z,\mathbf w)
    =
    \sum_{\nu:\,\mu\subseteq\nu\subseteq\lambda}
    s_{\lambda/\nu}(\mathbf z)s_{\nu/\mu}(\mathbf w),
\]
and Weyl's determinantal formula for ordinary Schur polynomials: for
\(\ell(\lambda)\le n\),
\[
    s_\lambda(z_1,\ldots,z_n)
    =
    \frac{\det[z_i^{\lambda_j+n-j}]_{i,j=1}^{n}}{\prod_{1\le i<j\le n}(z_i-z_j)} .
\]

Let
\[
    \mathbf{x}=(x_0,x_1,x_2,\ldots)
\]
be an infinite sequence of variables. For any tableau $T$ of shape
$\lambda/\mu$, define
\begin{align*}
    \mathbf{x}^T:=\prod_{(i,j)\in\lambda/\mu}x_{T(i,j)}.
\end{align*}
For $L\in LHT_n(\lambda/\mu)$, recall that
\[
    \lfloor L\rfloor(i,j)
    =
    \left\lfloor\frac{L(i,j)}{n+j-i}\right\rfloor .
\]
We define the lecture-hall generating function by
\begin{align*}
    L_{\lambda/\mu}^n(\mathbf{x})
    :=
    \sum_{L\in LHT_n(\lambda/\mu)}\mathbf{x}^{\lfloor L\rfloor}.
\end{align*}
Furthermore, we write
\begin{align*}
L_{\lambda}^n(\mathbf{x}):=L_{\lambda/\emptyset}^n(\mathbf{x}).
\end{align*}
Here, the superscript $n$ in $L_{\lambda}^n$ denotes a parameter, not a power.
When $n$ is clear from context, we may omit it from the notation.

\begin{definition}\label{df21}
Fix positive integers $n$ and $t$, and fix a partition
$\lambda^{(0)}\in\YY$ with $\ell(\lambda^{(0)})\leq n$, padded by zeros
to length $n$ when necessary. Let
$\mathcal{T}_{n,t}(\lambda^{(0)})$ denote the set of bounded
$n$-lecture hall tableaux of shape $\lambda^{(0)}$ and bound $t$; thus
$L\in\mathcal{T}_{n,t}(\lambda^{(0)})$ satisfies
\[
    L(i,j)<t(n+j-i),\qquad (i,j)\in\lambda^{(0)}.
\]
Under the path bijection recalled in Section~1, every such tableau corresponds
to a configuration of exactly $n$ non-intersecting paths on $\mathcal{G}_t$
with bottom boundary $\lambda^{(0)}$.

Let
\begin{align}
&\mathbf{x}=(x_0,x_1,\ldots,x_{t-1},0,0,\ldots),\label{dx}\\
&\mathbf{x}_{\kappa}=(x_{\kappa},x_{\kappa+1},\ldots,x_{t-1},0,0,\ldots),\label{dxk}\\
&\mathbf{x}\setminus\mathbf{x}_{\kappa}
=(x_0,x_1,\ldots,x_{\kappa-1},0,0,\ldots),\label{dxk2}
\end{align}
All three sequences in (\ref{dx})--(\ref{dxk2}) are infinite and eventually
zero. Throughout the paper, whenever a finite tuple is used as a
specialization of a symmetric function, it is identified with the infinite
sequence obtained by appending zeros; these trailing zeros may be omitted when
no confusion can arise. Here $x_0,\ldots,x_{t-1}$ are nonnegative, and empty
blocks are interpreted as the zero specialization. In particular,
\[
\mathbf{x}_0=\mathbf{x},\qquad
\mathbf{x}_t=(0,0,\ldots),\qquad
\mathbf{x}\setminus\mathbf{x}_0=(0,0,\ldots),\qquad
\mathbf{x}\setminus\mathbf{x}_t=\mathbf{x}.
\]
For $L\in\mathcal{T}_{n,t}(\lambda^{(0)})$, define
\[
    \operatorname{wt}_{\mathbf{x}}(L)
    :=\mathbf{x}^{\lfloor L\rfloor}
    =\prod_{(i,j)\in\lambda^{(0)}}
      x_{\left\lfloor L(i,j)/(n+j-i)\right\rfloor}.
\]
Since $x_r=0$ for $r\geq t$, the quantity
$L_{\lambda^{(0)}}^n(\mathbf{x})$ is the total weight of
$\mathcal{T}_{n,t}(\lambda^{(0)})$. Assume that it is positive. We define a
probability measure directly on bounded lecture hall tableaux by
\[
    \mathbb{P}_{\mathbf{x}}^{n,t,\lambda^{(0)}}(L)
    :=\frac{\operatorname{wt}_{\mathbf{x}}(L)}{L_{\lambda^{(0)}}^n(\mathbf{x})},
    \qquad L\in\mathcal{T}_{n,t}(\lambda^{(0)}).
\]
Through the path bijection, and through the corresponding dimer bijection when
used, this measure induces the associated probability measures on path and
dimer configurations.

For $0\leq\kappa<t$, let
$\lambda^{(\kappa)}=(\lambda_1^{(\kappa)},\ldots,
\lambda_n^{(\kappa)})$ be the partition read from the path configuration at
the horizontal level $y\in(\kappa,\kappa+\epsilon]$, as defined above. We
define the top boundary by
\[
    \lambda^{(t)}:=\emptyset.
\]
For $0\leq\kappa\leq t$, let
$\rho_{\kappa}^{(\mathbf{x})}$ be the distribution of
$\lambda^{(\kappa)}$ under
$\mathbb{P}_{\mathbf{x}}^{n,t,\lambda^{(0)}}$. Equivalently, for
$\lambda\in\YY$,
\[
\rho_{\kappa}^{(\mathbf{x})}
\bigl(\lambda^{(\kappa)}=\lambda\bigr)
=
\frac{
L_{\lambda}^n(\mathbf{x}_{\kappa})
L_{\lambda^{(0)}/\lambda}^n
   (\mathbf{x}\setminus\mathbf{x}_{\kappa})
}{
L_{\lambda^{(0)}}^n(\mathbf{x})
}.
\]
Here we adopt the convention
\[
L_{\lambda^{(0)}/\lambda}^n(\mathbf{y})=0
\qquad\text{unless}\qquad
\lambda\subseteq\lambda^{(0)}
\]
for every specialization $\mathbf{y}$. In particular, $\rho_0^{(\mathbf{x})}$ is the point mass at
$\lambda^{(0)}$, while $\rho_t^{(\mathbf{x})}$ is the point mass at
$\emptyset$.
When no confusion can arise, we suppress the dependence on $\mathbf{x}$ and
write $\rho_{\kappa}$.
\end{definition}

\begin{definition}\label{df01}
Let $\rho_{\kappa}^{(\mathbf{x})}$ be the probability distribution of
$\lambda^{(\kappa)}$ defined in Definition \ref{df21}. For any scalar $a$
for which the denominators below are nonzero, define its Schur generating
function by
\begin{eqnarray*}
\mathcal{S}_{\rho_{\kappa}^{(\mathbf{x})}}(a,\mathbf{u})
=
\sum_{\lambda\in\YY}
\rho_{\kappa}^{(\mathbf{x})}(\lambda)
\frac{s_{\lambda}(a+\mathbf{u})}{s_{\lambda}(a)},
\end{eqnarray*}
where
\begin{align}
\mathbf{u}=(u_1,u_2,\ldots,u_n,0,0,\ldots),\qquad
|\mathbf{x}|:=\sum_{r=0}^{t-1}x_r=x_0+x_1+\cdots+x_{t-1}\label{dux}
\end{align}
Moreover,
\[
|\mathbf{x}_{\kappa}|:=\sum_{r=\kappa}^{t-1}x_r,
\qquad
|\mathbf{x}\setminus\mathbf{x}_{\kappa}|:=\sum_{r=0}^{\kappa-1}x_r.
\]
For a scalar $a$, we use $a+\mathbf{u}$ as shorthand for the infinite
sequence $(a+u_1,\ldots,a+u_n,0,0,\ldots)$. Likewise, in
$s_{\lambda}(a)$, the scalar $a$ denotes the infinite sequence consisting of
$n$ copies of $a$ followed by zeros. Thus
\begin{eqnarray}
&&s_{\lambda}(a+\mathbf{u})
:=s_{\lambda}(a+u_1,a+u_2,\ldots,a+u_n,0,0,\ldots)\label{sxu}\\
&&s_{\lambda}(a)
:=s_{\lambda}(\underbrace{a,\ldots,a}_{n\text{ copies}},0,0,\ldots).
\end{eqnarray}
The superscript $(\mathbf{x})$ records the vector of weights that determines
$\rho_{\kappa}^{(\mathbf{x})}$, whereas the scalar $a$ is only the evaluation
parameter in the Schur ratio. Thus changing $a$ from $|\mathbf{x}|$ to
$|\mathbf{x}_{\kappa}|$ does not change the probability distribution. When no
confusion can arise, we continue to suppress the superscript $(\mathbf{x})$.
\end{definition}

\begin{lemma}\label{l02}Let $\lambda\in\YY$ with $l(\lambda)\leq n$. Let
\begin{eqnarray*}
\mathbf{a}=(a_0,a_1,\ldots,a_{t-1},0,0,\ldots),\qquad
\mathbf{b}=(b_1,\ldots,b_n).
\end{eqnarray*}
Then 
\begin{eqnarray*}
s_{\lambda}(|\mathbf{a}|+\mathbf{b})=\sum_{\nu\subset\lambda}L_{\lambda/\nu}(\mathbf{a})s_{\nu}(\mathbf{b}),
\end{eqnarray*}
where $|\mathbf{a}|=a_0+a_1+\cdots+a_{t-1}$ and $s_{\lambda}(|\mathbf{a}|+\mathbf{b})$ is as in (\ref{sxu}), and
\begin{eqnarray*}
s_{\nu}(\mathbf{b})=
s_{\nu}(b_1,b_2,\ldots,b_n).
\end{eqnarray*}
\end{lemma}
\begin{proof}The lemma    is the branching identity for the lecture-hall generating functions proved in  Theorem 1.6 of \cite{CK20}   , applied with the outer summation over the intermediate partition $\nu$ .
\end{proof}

We shall use the following immediate specialization of Lemma \ref{l02}:
for every partition $\eta$ with $\ell(\eta)\leq n$ and every specialization $\mathbf{a}$ of the form used in Lemma \ref{l02},
\[
    L_{\eta}^n(\mathbf{a})=s_{\eta}(|\mathbf{a}|).
\]
Indeed, set $\mathbf{b}=\mathbf{0}$ in Lemma \ref{l02}. Then
$s_{\nu}(\mathbf{0})=0$ for every nonempty partition $\nu$, whereas
$s_{\emptyset}(\mathbf{0})=1$, so only the term $\nu=\emptyset$ remains.

\begin{lemma}\label{l03}Assume the partition on the bottom boundary $\lambda^{(0)}$ is fixed, and the partition on the top boundary is given by the empty partition. Then for any $0<\kappa<t$,
\begin{eqnarray*}
\mathcal{S}_{\rho_{\kappa}^{(\mathbf{x})}}
(|\mathbf{x}_{\kappa}|,\mathbf{u})=
\frac{s_{\lambda^{(0)}}(|\mathbf{x}|+\mathbf{u})}{s_{
\lambda^{(0)}}(|\mathbf{x}|)},
\end{eqnarray*}
where    $\mathbf{x}_{\kappa}$  is defined by (\ref{dxk}). Here the measure on the left-hand side is the same
$\rho_{\kappa}^{(\mathbf{x})}$ determined by the full weight vector $\mathbf{x}$
in Definition \ref{df21}; only the scalar evaluation parameter of its Schur
generating function is $|\mathbf{x}_{\kappa}|$.
\end{lemma}

\begin{proof}Recall that 
$\mathbf{x}\setminus \mathbf{x}_{\kappa}=(x_0,x_1,\ldots,x_{\kappa-1},0,0,\ldots)$ is defined by (\ref{dxk2}).

By Definition \ref{df01}, we have
\begin{align}
\mathcal{S}_{\rho_{\kappa}^{(\mathbf{x})}}
(|\mathbf{x}_{\kappa}|,\mathbf{u})
&=\sum_{\lambda\in \YY}\rho_{\kappa}^{(\mathbf{x})}(\lambda)\frac{s_{\lambda}(|\mathbf{x}_{\kappa}|+\mathbf{u})}{s_{\lambda}(|\mathbf{x}_{\kappa}|)}\label{srk}\\
&=\sum_{\lambda\in \YY}\frac{L_{\lambda}(\mathbf{x}_{\kappa})L_{\lambda^{(0)}/ \lambda}(\mathbf{x}\setminus \mathbf{x}_{\kappa} )}{L_{\lambda^{(0)}}(\mathbf{x})}\frac{s_{\lambda}(|\mathbf{x}_{\kappa}|+\mathbf{u})}{s_{\lambda}(|\mathbf{x}_{\kappa}|)}\notag\\
&=\sum_{\lambda\in \YY}\frac{L_{\lambda^{(0)}/ \lambda}(\mathbf{x}\setminus \mathbf{x}_{\kappa} )s_{\lambda}(|\mathbf{x}_{\kappa}|+\mathbf{u})}{L_{\lambda^{(0)}}(\mathbf{x})}\notag\\
&=\frac{s_{\lambda^{(0)}}(|\mathbf{x}|+\mathbf{u})}{s_{
\lambda^{(0)}}(|\mathbf{x}|)},\notag
\end{align}
In passing from the second line to the third line, we used the specialization
\[
    L_{\lambda}^n(\mathbf{x}_{\kappa})
    =s_{\lambda}(|\mathbf{x}_{\kappa}|)
\]
displayed above, which cancels the corresponding denominator. For the last
line, apply Lemma \ref{l02} with
\[
    \mathbf{a}=\mathbf{x}\setminus\mathbf{x}_{\kappa},
    \qquad
    \mathbf{b}=|\mathbf{x}_{\kappa}|+\mathbf{u}.
\]
Since
\[
    |\mathbf{x}\setminus\mathbf{x}_{\kappa}|
    +|\mathbf{x}_{\kappa}|=|\mathbf{x}|,
\]
the numerator becomes $s_{\lambda^{(0)}}(|\mathbf{x}|+\mathbf{u})$.
Finally, the same specialization identity gives
$L_{\lambda^{(0)}}^n(\mathbf{x})=s_{\lambda^{(0)}}(|\mathbf{x}|)$,
which yields the displayed ratio.
\end{proof}

Define the following differential operator $\mathcal{D}_{j,\kappa}$ 
 
  \begin{eqnarray*}
\mathcal{D}_{j,\kappa} :=\frac{1}{V(\mathbf{u})}
\left[\sum_{i=1}^{n}\left((|\mathbf{x}_{\kappa}|+u_i)\frac{\partial}{\partial u_i}\right)^j\right]
V(\mathbf{u}),
\end{eqnarray*} 
where
\begin{eqnarray*}
V(\mathbf{u})=\prod_{i<j}(u_i-u_j).
\end{eqnarray*}

We introduce the following definition to study the distribution of random partitions.
\begin{definition}\label{dfcm}Let $\lambda$ be a length-$n$  partition. We define the counting measure $m(\lambda)$ as    the  probability measure on $\RR$  
  given by
\begin{eqnarray*}
m(\lambda)=\frac{1}{n}\sum_{i=1}^{n}\delta_{(\lambda_i+n-i)/n}.
\end{eqnarray*} 
Here $\delta_a$ denotes the Dirac mass at $a$.  If $\lambda$ is random, then    $m(\lambda)$ is  the corresponding random counting measure.
\end{definition}

\noindent\textbf{Schur-operator identity.}
For later use, we record the following identity. For every $j,m\in\NN$,
\begin{align}
&\mathcal{D}_{j,\kappa}^{m}
\mathcal{S}_{\rho_{\kappa}}(|\mathbf{x}_{\kappa}|,\mathbf{u})\notag\\
&\qquad=
\sum_{\lambda\in\YY}\rho_{\kappa}(\lambda)
\left[\sum_{i=1}^{n}(\lambda_i+n-i)^j\right]^m
\frac{s_{\lambda}(|\mathbf{x}_{\kappa}|+\mathbf{u})}
{s_{\lambda}(|\mathbf{x}_{\kappa}|)}.
\label{eq:schur-operator-action}
\end{align}
Indeed, set $a:=|\mathbf{x}_{\kappa}|$ and
$d_r:=\lambda_r+n-r$ for $1\leq r\leq n$. Since
$V(\mathbf{u})=\prod_{p<q}((a+u_p)-(a+u_q))$, Weyl's determinantal
formula gives
\begin{align*}
V(\mathbf{u})s_{\lambda}(a+\mathbf{u})
&=\det\!\left[(a+u_i)^{d_r}\right]_{i,r=1}^{n}\\
&=\sum_{\sigma\in\mathfrak S_n}\operatorname{sgn}(\sigma)
  \prod_{i=1}^{n}(a+u_i)^{d_{\sigma(i)}}.
\end{align*}
If $E_i:=(a+u_i)\frac{\partial}{\partial u_i}$, then every monomial in
this alternant is an eigenfunction of $\sum_{i=1}^{n}E_i^j$ with the
same eigenvalue $\sum_{r=1}^{n}d_r^j$. Hence
\begin{align*}
\frac{1}{V(\mathbf{u})}
\left[\sum_{i=1}^{n}
\left((|\mathbf{x}_{\kappa}|+u_i)
\frac{\partial}{\partial u_i}\right)^j\right]^m
V(\mathbf{u})s_{\lambda}(|\mathbf{x}_{\kappa}|+\mathbf{u})
=
\left[\sum_{i=1}^{n}(\lambda_i+n-i)^j\right]^m
s_{\lambda}(|\mathbf{x}_{\kappa}|+\mathbf{u}).
\end{align*}
Substituting this identity into Definition~\ref{df01} proves
\eqref{eq:schur-operator-action}.

The following lemma relates the moments of the random counting measure to the operator $\mathcal{D}_{j,\kappa}$. 
\begin{lemma}\label{l04}Let $j,m\in\NN$. Then
 
  \begin{eqnarray*}
\left.\frac{1}{n^{(j+1)m}}\mathcal{D}_{j,\kappa}^m \mathcal{S}_{\rho_{\kappa}}(|\mathbf{x}_{\kappa}|,\mathbf{u})\right|_{\mathbf{u}=0}=\mathbb{E}\left[\left(\int_{\RR} x^{j}\,d\mathbf{m}_{\rho_{\kappa}}(x)\right)^m\right]
\end{eqnarray*} 
where $\mathbf{m}_{\rho_{\kappa}}$ is the random counting measure for the random partition $\lambda^{(\kappa)}$ defined in Definition \ref{df21}. 
\end{lemma}
\begin{proof}
Evaluating \eqref{eq:schur-operator-action} at $\mathbf{u}=\mathbf{0}$ gives
\begin{align*}
\left.\mathcal{D}_{j,\kappa}^{m}
\mathcal{S}_{\rho_{\kappa}}(|\mathbf{x}_{\kappa}|,\mathbf{u})
\right|_{\mathbf{u}=\mathbf{0}}
=
\sum_{\lambda\in\YY_n}\rho_{\kappa}(\lambda)
\left[\sum_{i=1}^{n}(\lambda_i+n-i)^j\right]^m.
\end{align*}
By the definition of the counting measure,
\[
\int_{\RR}x^j\,d m(\lambda)(x)
=
\frac{1}{n^{j+1}}\sum_{i=1}^{n}(\lambda_i+n-i)^j.
\]
Dividing the preceding identity by $n^{(j+1)m}$ proves the lemma.
\end{proof}

\begin{definition}[Regular sequences of partitions]
\label{def:regular-bottom}
For every $n\geq1$, let
\[
    \lambda(n)
    =
    \bigl(\lambda_1(n)\geq\lambda_2(n)\geq\cdots\geq
    \lambda_n(n)\geq0\bigr)
    \in\YY_n,
\]
where partitions of length less than $n$ are padded by zeros. The
sequence $\{\lambda(n)\}_{n\geq1}$ is called \emph{regular} if there
exist a bounded, piecewise-continuous, nonincreasing function
$f:[0,1]\to[0,\infty)$ and a constant $C>0$ such that
\begin{align}
    \frac1n\sum_{i=1}^n
    \left|
        \frac{\lambda_i(n)}n-f\left(\frac{i}{n}\right)
    \right|
    &\longrightarrow0,
    \label{eq:regular-bottom-L1}\\
    \left|
        \frac{\lambda_i(n)}n-f\left(\frac{i}{n}\right)
    \right|
    &\leq C,
    \qquad 1\leq i\leq n,\quad n\geq1.
    \label{eq:regular-bottom-uniform}
\end{align}
The function $f$ is called the limiting profile of the sequence, and
the associated probability measure is
\begin{equation}
    \mathbf m_f
    :=
    \bigl(u\mapsto f(u)+1-u\bigr)_*
    \operatorname{Leb}_{[0,1]}.
    \label{eq:measure-from-regular-profile}
\end{equation}
Here $T_*\mu$ denotes the pushforward of $\mu$ under the map $T$, and
$\operatorname{Leb}_{[0,1]}$ denotes Lebesgue measure on $[0,1]$.
When $\lambda(n)=\lambda^{(0)}(n)$ is the bottom-boundary sequence, we
write $f_0$ for its limiting profile and $\mathbf m_0$ for the measure
in \eqref{eq:measure-from-regular-profile}, and call
$\{\lambda^{(0)}(n)\}_{n\geq1}$ a regular bottom-boundary sequence.
\end{definition}

\begin{lemma}
\label{lem:regular-implies-moment-convergence}
If $\{\lambda(n)\}_{n\geq1}$ is regular with limiting profile $f$, then
the counting measures defined in Definition~\ref{dfcm} converge weakly to
$\mathbf m_f$; that is,
\[
    m(\lambda(n))\Longrightarrow\mathbf m_f.
\]
The convergence also holds in the sense of moments. Moreover,
$\mathbf m_f$ and all the measures $m(\lambda(n))$, $n\geq1$, are
supported in a common compact interval.
\end{lemma}

\begin{proof}
For $1\leq i\leq n$, set
\[
    q_{i,n}:=\frac{\lambda_i(n)+n-i}{n}.
\]
By \eqref{eq:regular-bottom-uniform},
\[
    0\leq q_{i,n}\leq \lVert f\rVert_\infty+C+1,
\]
so the counting measures have uniformly bounded supports. Also,
\[
    \operatorname{supp}(\mathbf m_f)
    \subseteq[0,\lVert f\rVert_\infty+1].
\]
Moreover,
\eqref{eq:regular-bottom-L1} gives
\[
    \frac1n\sum_{i=1}^n
    \left|
        q_{i,n}-\left(f\left(\frac{i}{n}\right)+1-\frac{i}{n}\right)
    \right|
    \longrightarrow0.
\]
The assertion follows by Riemann-sum approximation, first for
polynomials and then for continuous functions on the common compact
interval containing the supports.
\end{proof}

\begin{theorem}\label{t15}
For each $n\geq 1$, let $t_n$ be a positive integer, let
$\kappa_n\in\{0,1,\ldots,t_n\}$, and let $\lambda^{(0)}(n)$ be a
bottom-boundary partition with $\ell(\lambda^{(0)}(n))\leq n$ for a system of
$n$ non-intersecting paths on $\mathcal G_{t_n}$. Let
\[
\mathbf{x}^{(n)}
=
\bigl(x_0^{(n)},x_1^{(n)},\ldots,x_{t_n-1}^{(n)},0,0,\ldots\bigr)
\]
be a nonnegative weight specialization such that
$L_{\lambda^{(0)}(n)}^n(\mathbf{x}^{(n)})>0$. Define
\[
\rho_{\kappa_n}^{(n)}
:=
\rho_{\kappa_n}^{(\mathbf{x}^{(n)})},
\]
where the distribution on the right-hand side is the one in Definition
\ref{df21}, with parameters $n$, $t_n$, and $\lambda^{(0)}(n)$. Let
$\lambda^{(\kappa_n)}(n)$ be a random partition with this distribution and let
\[
\mathbf m_n:=m\bigl(\lambda^{(\kappa_n)}(n)\bigr)
\]
be its random counting measure. Set
\[
|\mathbf{x}^{(n)}|
:=\sum_{r=0}^{t_n-1}x_r^{(n)},
\qquad
|\mathbf{x}_{\kappa_n}^{(n)}|
:=\sum_{r=\kappa_n}^{t_n-1}x_r^{(n)}.
\]
Assume that
\begin{equation}
 y:=\lim_{n\to\infty}\frac{\kappa_n}{n},
 \qquad
 s:=\lim_{n\to\infty}
 \frac{|\mathbf{x}_{\kappa_n}^{(n)}|}{|\mathbf{x}^{(n)}|},
 \qquad
 \alpha:=\lim_{n\to\infty}\frac{t_n}{n},
 \label{lms1}
\end{equation}
where
\[
\alpha\in(0,\infty),
\qquad
s\in(0,1),
\qquad
0<y<\alpha.
\]

Assume also that $\{\lambda^{(0)}(n)\}_{n\geq1}$ is a regular
bottom-boundary sequence in the sense of
Definition~\ref{def:regular-bottom}, and let $\mathbf m_0$ be its
associated compactly supported probability measure defined by
\eqref{eq:measure-from-regular-profile}.

Then there is a deterministic probability measure $\mathbf m_y$ on $\mathbb R$
such that, for every fixed $j\geq 0$,
\[
X_{j,n}:=\int_{\mathbb R}x^j\,\mathbf m_n(dx)
\longrightarrow
\int_{\mathbb R}x^j\,\mathbf m_y(dx)
\quad\text{in }L^2,
\]
and hence in probability. Equivalently, the random measures $\mathbf m_n$
converge in probability to $\mathbf m_y$ in the sense of moments.

Let \(f_0\) and \(C\) be the profile and the constant in
Definition~\ref{def:regular-bottom}, and set
\[
    B:=\lVert f_0\rVert_\infty+C+1.
\]
Then
\[
    \operatorname{supp}(\mathbf m_n)\subset[0,B]
    \quad\text{almost surely},
    \qquad
    \operatorname{supp}(\mathbf m_y)\subset[0,B],
\]
and
\[
    \mathbf m_n\Longrightarrow\mathbf m_y
    \qquad\text{weakly in probability}.
\]

The limiting moments are
\[
\int_{\mathbb R}x^j\,\mathbf m_y(dx)
=
\frac{1}{2(j+1)\pi\mathbf i}
\oint_{\gamma_1}
\frac{1}{z-1+s}
\left(
 (z-1+s)H_{\mathbf m_0}'(z)
 +\frac{z-1+s}{z-1}
\right)^{j+1}
\,dz.
\]
Here $\gamma_1=\{z:|z-1|=r\}$ is oriented counterclockwise, where $r>0$ is
chosen so that $r<s$, $H_{\mathbf m_0}'$ is analytic on
$\{z:|z-1|\leq r\}$, and $z=1$ is the only singularity of the integrand inside
$\gamma_1$. The function $H_{\mathbf m_0}$ is defined in (\ref{hmz}). The
notation $\mathbf m_y$ suppresses the dependence on the limiting weights and on
$\alpha$: the moment formula depends on the macroscopic level and on the aspect
ratio through the tail-weight ratio $s$. In particular, for uniform weights,
\[
s=1-\frac{y}{\alpha}.
\]
\end{theorem}

\begingroup

\begin{proof}
Set
\[
    \rho_n:=\rho_{\kappa_n}^{(n)},
    \qquad
    a_n:=|\mathbf x_{\kappa_n}^{(n)}|,
    \qquad
    b_n:=|\mathbf x^{(n)}|,
    \qquad
    s_n:=\frac{a_n}{b_n}.
\]
By \eqref{lms1}, \(s_n\to s\in(0,1)\). Moreover,
\(0<y<\alpha\) implies \(0<\kappa_n<t_n\) for all sufficiently large
\(n\). The hypotheses also give \(a_n,b_n>0\) for all sufficiently
large \(n\). Discarding finitely many indices does not affect the
conclusion, so we henceforth restrict to such \(n\).

For
\[
    \mathbf z=(z_1,\ldots,z_n),
    \qquad
    \mathbf 1_n=(1,\ldots,1),
\]
define the standard Schur generating function of \(\rho_n\) by
\[
    \widehat{\mathcal S}_n(\mathbf z)
    :=
    \sum_{\lambda\in\YY_n}
    \rho_n(\lambda)
    \frac{s_\lambda(\mathbf z)}{s_\lambda(1^n)}.
\]
Definition~\ref{df01} and homogeneity of Schur polynomials give the
first equality below, Lemma~\ref{l03} gives the second, and homogeneity
together with \(s_n=a_n/b_n\) gives the last:
\begin{align}
    \widehat{\mathcal S}_n(\mathbf z)
    &=
    \mathcal S_{\rho_n}
    \bigl(a_n,a_n(\mathbf z-\mathbf 1_n)\bigr)
    \notag\\
    &=
    \frac{
        s_{\lambda^{(0)}(n)}
        \bigl(b_n\mathbf 1_n+a_n(\mathbf z-\mathbf 1_n)\bigr)
    }{
        s_{\lambda^{(0)}(n)}(b_n\mathbf 1_n)
    }
    \notag\\
    &=
    \frac{
        s_{\lambda^{(0)}(n)}
        \bigl(\mathbf 1_n+s_n(\mathbf z-\mathbf 1_n)\bigr)
    }{
        s_{\lambda^{(0)}(n)}(1^n)
    }.
    \label{eq:standard-level-schur}
\end{align}

Fix \(k\geq1\). Since the convergence in Lemma~\ref{la1} is locally
uniform and \(s_n\to s\), specializing
\(z_{k+1}=\cdots=z_n=1\) in
\eqref{eq:standard-level-schur} yields
\begin{equation}
    \frac1n
    \log
    \widehat{\mathcal S}_n
    (z_1,\ldots,z_k,1^{n-k})
    \longrightarrow
    \sum_{i=1}^k
    H_{\mathbf m_0}\bigl(1+s(z_i-1)\bigr)
    \label{eq:level-schur-lln}
\end{equation}
locally uniformly in a fixed complex neighborhood of
\((1,\ldots,1)\). Here and below, the logarithm is taken in the
holomorphic branch normalized to vanish at \((1,\ldots,1)\).

Cauchy's integral formula permits differentiation of
\eqref{eq:level-schur-lln}; hence, locally uniformly,
\[
    \frac1n
    \frac{\partial}{\partial z_1}
    \log
    \widehat{\mathcal S}_n
    (z_1,\ldots,z_k,1^{n-k})
    \longrightarrow
    F_\rho(z_1),
\]
where
\[
    F_\rho(z)
    :=
    sH_{\mathbf m_0}'\bigl(1+s(z-1)\bigr).
\]
Thus the criterion of
\cite[General Example~2.3]{bg16} is satisfied, and
\(\{\rho_n\}_{n\geq1}\) is LLN-appropriate.

It now follows from \cite[Theorem~2.4]{bg16} that the random counting
measures \(\mathbf m_n\) converge in probability, in the sense of
moments, to a deterministic probability measure, which we denote by
\(\mathbf m_y\). For every \(j\geq1\), writing
\[
    M_j:=\int_{\mathbb R}x^j\,\mathbf m_y(dx),
\]
the same theorem gives
\[
    M_j
    =
    \frac{1}{2(j+1)\pi\mathbf i}
    \oint_{|\zeta|=\epsilon}
    \frac{1}{1+\zeta}
    \left(
        \frac1\zeta+1+
        s(1+\zeta)
        H_{\mathbf m_0}'(1+s\zeta)
    \right)^{j+1}
    d\zeta,
\]
where \(\epsilon>0\) is sufficiently small. Under the change of
variables \(z=1+s\zeta\),
\[
    \frac{d\zeta}{1+\zeta}
    =
    \frac{dz}{z-1+s},
\]
and
\[
    \frac1\zeta+1+
    s(1+\zeta)H_{\mathbf m_0}'(1+s\zeta)
    =
    \frac{z-1+s}{z-1}
    +(z-1+s)H_{\mathbf m_0}'(z).
\]
Consequently,
\[
    M_j
    =
    \frac{1}{2(j+1)\pi\mathbf i}
    \oint_{\gamma_1}
    \frac{1}{z-1+s}
    \left(
        (z-1+s)H_{\mathbf m_0}'(z)
        +\frac{z-1+s}{z-1}
    \right)^{j+1}
    dz.
\]
The image of \(|\zeta|=\epsilon\) is a sufficiently small positively
oriented circle about \(z=1\). By the analyticity assumptions in the
statement, this circle may be deformed to \(\gamma_1\) without
crossing any singularity.

For \(j=0\), the convergence assertion follows from
\[
    \int_{\mathbb R}1\,\mathbf m_n(dx)=1.
\]
The corresponding contour formula also equals \(1\), since its
integrand reduces to
\[
    H_{\mathbf m_0}'(z)+\frac1{z-1},
\]
whose residue at \(z=1\) is \(1\).

It remains to strengthen convergence in probability to \(L^2\).
By Definition~\ref{df21}, every partition in the support of \(\rho_n\)
is contained in \(\lambda^{(0)}(n)\). Let \(f_0\) and \(C\) be the
profile and constant from Definition~\ref{def:regular-bottom}, and set
\[
    B:=\lVert f_0\rVert_\infty+C+1.
\]
For every \(\lambda\in\operatorname{supp}\rho_n\) and
\(1\leq i\leq n\), condition
\eqref{eq:regular-bottom-uniform} gives
\[
    0
    \leq
    \frac{\lambda_i+n-i}{n}
    \leq
    \frac{\lambda_i^{(0)}(n)+n-i}{n}
    \leq
    f_0\left(\frac{i}{n}\right)+C+1
    \leq B.
\]
Therefore, almost surely,
\[
    0\leq
    X_{j,n}
    =
    \frac1n\sum_{i=1}^n
    \left(
        \frac{\lambda_i^{(\kappa_n)}(n)+n-i}{n}
    \right)^j
    \leq B^j.
\]
Since \(X_{j,n}\to M_j\) in probability, necessarily
\(0\leq M_j\leq B^j\). Hence, for every \(\delta>0\),
\[
    \mathbb E|X_{j,n}-M_j|^2
    \leq
    \delta^2+
    B^{2j}
    \mathbb P\bigl(|X_{j,n}-M_j|>\delta\bigr).
\]
Letting \(n\to\infty\) and then \(\delta\downarrow0\) proves that
\(X_{j,n}\to M_j\) in \(L^2\).

Taking \(j=2r\) in the preceding bound gives
\[
    \int_{\mathbb R}|x|^{2r}\,\mathbf m_y(dx)
    \leq B^{2r},
    \qquad r\geq1,
\]
and hence
\[
    \operatorname{supp}(\mathbf m_y)\subset[-B,B].
\]
Thus \(\mathbf m_n\) and \(\mathbf m_y\) are supported in a common
compact interval. The convergence of all moments proved above and
polynomial approximation therefore imply
\[
    \mathbf m_n\Longrightarrow\mathbf m_y
    \qquad\text{weakly in probability}.
\]
Indeed, from every subsequence one may extract a further subsequence
along which all moments converge almost surely simultaneously;
polynomial approximation then gives almost-sure weak convergence
along that further subsequence.

Finally, if
\(\psi\in C_c(\mathbb R\setminus[0,B])\) is nonnegative, then
\[
    \int_{\mathbb R}\psi\,d\mathbf m_n=0
    \qquad\text{almost surely}.
\]
The weak convergence and the fact that \(\mathbf m_y\) is
deterministic give
\[
    \int_{\mathbb R}\psi\,d\mathbf m_y=0.
\]
Consequently,
\[
    \operatorname{supp}(\mathbf m_y)\subset[0,B].
\]

\end{proof}

\endgroup

\begingroup

\begin{assumption}[Limiting scaled path domain]
\label{ass:scaled-domain}
For each \(n\), set
\[
    N_n
    :=
    n-1+\lambda^{(0)}_1(n)
    =
    \max_{1\leq i\leq n}
    \bigl\{n-i+\lambda_i^{(0)}(n)\bigr\},
\]
and let \(\mathcal G_{t_n,N_n}\) be the finite truncation defined in
Definition~\ref{dflhg}(\ref{item:finite-horizontal-truncation}).
Denote by \(\lvert\mathcal G_{t_n,N_n}\rvert\) the geometric
realization of the embedded graph, and let
\(\mathsf F_b(\mathcal G_{t_n,N_n})\) be the collection of its
bounded faces. Define
\[
    \mathcal D_n
    :=
    \lvert\mathcal G_{t_n,N_n}\rvert
    \cup
    \bigcup_{F\in\mathsf F_b(\mathcal G_{t_n,N_n})}\overline F.
\]
Thus \(\mathcal D_n\) is obtained by filling all bounded faces of the
finite embedded graph. Since every admissible path uses only edges of
\(\mathcal G_{t_n,N_n}\), the embedded image of every admissible path
configuration is contained in
\(\lvert\mathcal G_{t_n,N_n}\rvert\subset\mathcal D_n\). Indeed, its
paths move only to the right and their terminal horizontal coordinates
are bounded above by \(N_n\).

Define the compact rescaled domain
\[
    \mathcal R_n
    :=
    \frac1n\mathcal D_n
    =
    \left\{
        \left(\frac{x}{n},\frac{y}{n}\right):(x,y)\in\mathcal D_n
      \right\}.
\]
Assume that there exists a compact set
\(\mathcal R=\overline{\mathcal R^{\circ}}\subset\mathbb R^2\), where
\(\mathcal R^{\circ}\) is a bounded, simply connected domain, such that
\[
    d_{\mathrm H}(\mathcal R_n,\mathcal R)
    \longrightarrow0,
\]
where \(d_{\mathrm H}\) is the Hausdorff metric induced by the Euclidean
distance. We call \(\mathcal R\) the limiting scaled path domain.
All subsequent statements involving \(\mathcal R\) or \(\mathcal L\)
are understood under this standing assumption.

Let \(C_0\) be a constant for which
\eqref{eq:regular-bottom-uniform} holds for the bottom-boundary
sequence. Then
\[
    \frac{N_n}{n}
    \leq
    \lVert f_0\rVert_\infty+C_0+1.
\]
Together with \(t_n/n\to\alpha\), this shows that the compact sets
\(\mathcal R_n\) lie in a common bounded rectangle.
\end{assumption}
\endgroup

\begin{definition}\label{def:transforms}
Let $\eta$ be a compactly supported probability measure on $\mathbb R$.  Its
Stieltjes transform is
\[
\mathrm{St}_\eta(w):=\int_{\mathbb R}\frac{\eta(d\xi)}{w-\xi},
\qquad w\in\mathbb C\setminus\operatorname{supp}(\eta).
\]
Its ordinary moment generating function (or moment series) is
\[
S_\eta(z):=\sum_{k=0}^{\infty}M_k(\eta)z^{k+1}
=z+\sum_{k=1}^{\infty}M_k(\eta)z^{k+1},
\qquad
M_k(\eta):=\int_{\mathbb R}\xi^k\,\eta(d\xi).
\]
For $|w|$ sufficiently large, equivalently for $z$ sufficiently close to
$0$, the two transforms satisfy
\[
\mathrm{St}_\eta(w)=S_\eta(w^{-1}),
\qquad
S_\eta(z)=\mathrm{St}_\eta(z^{-1}).
\]
Since $S_\eta(0)=0$ and $S_\eta'(0)=1$, it has a local inverse under
composition near $0$, denoted by $S_\eta^{(-1)}$.
\end{definition}

  For fixed $y$ and the corresponding value of $s$ in Theorem \ref{t15}, set
\begin{equation}
U_y(z):=(z-1+s)H'_{\mathbf m_0}(z)+\frac{z-1+s}{z-1}.
\label{duyz}
\end{equation}

\begin{proposition}[Stieltjes transform of the limit measure]
\label{prop:stieltjes-root}
For $w$ of sufficiently large modulus, let $\mathbf z_1(w,y)$ be the unique
solution of
 \[
U_y(z)=w
\]
that    satisfies $\mathbf z_1(w,y)\to1$ as $w\to\infty$. Then
\begin{equation}
\mathrm{St}_{\mathbf m_y}(w)
=
\operatorname{Log}\!\left(
\frac{\mathbf z_1(w,y)-1+s}{s}
\right),
\label{eq:st-root-large-w}
\end{equation}
where the logarithm is normalized by requiring the right-hand side to tend to
$0$ as $w\to\infty$.  Moreover, the identity extends by analytic
continuation to every connected domain that meets the large-$w$ domain and on
which the physical branch $\mathbf z_1(\,\cdot\,,y)$ and the normalized
logarithm are holomorphic.
\end{proposition}

\begin{proof}
Set
\[
    Q_y(z):=\frac{1}{U_y(z)}.
\]
Since \(H'_{\mathbf m_0}\) is analytic in a neighborhood of \(z=1\),
equation \eqref{duyz} gives
\[
    U_y(z)
    =
    \frac{s}{z-1}+O(1),
    \qquad z\to1.
\]
Hence
\[
    Q_y(z)
    =
    \frac{z-1}{s}
    +O\bigl((z-1)^2\bigr),
    \qquad z\to1.
\]
Thus the apparent singularity of \(Q_y\) at \(z=1\) is removable, and
the holomorphic extension satisfies
\[
    Q_y(1)=0,
    \qquad
    Q_y'(1)=\frac1s\neq0.
\]
By the complex inverse-function theorem, \(Q_y\) has a local
holomorphic inverse \(Z_y\) near \(q=0\), characterized by
\[
    Q_y(Z_y(q))=q,
    \qquad
    Z_y(0)=1.
\]
For nonzero \(q\) sufficiently close to \(0\), we have
\[
    U_y(Z_y(q))=\frac1q.
\]
By the defining asymptotic condition for the branch
\(\mathbf z_1(w,y)\), this implies
\[
    Z_y(q)=\mathbf z_1(q^{-1},y).
\]

Define
\[
    \Phi(z)
    :=
    \operatorname{Log}\!
    \left(
        \frac{z-1+s}{s}
    \right),
\]
where the branch of \(\operatorname{Log}\) is chosen near \(z=1\) so
that \(\Phi(1)=0\). Then
\[
    \Phi'(z)=\frac{1}{z-1+s}.
\]
After shrinking \(\gamma_1\), if necessary, we may assume that
\(Q_y\) is one-to-one in a neighborhood of the closed region bounded
by \(\gamma_1\), that \(z=1\) is its only zero inside \(\gamma_1\),
and that \(\Phi\) is holomorphic there.

For every integer \(\ell\geq1\), Cauchy's coefficient formula,
followed by the change of variables \(q=Q_y(z)\), gives
\begin{align*}
[q^\ell]\Phi(Z_y(q))
&=
\frac{1}{2\pi\mathbf i}
\oint_{Q_y(\gamma_1)}
\frac{\Phi(Z_y(\zeta))}{\zeta^{\ell+1}}\,d\zeta\\
&=
\frac{1}{2\pi\mathbf i}
\oint_{\gamma_1}
\Phi(z)
\frac{Q_y'(z)}{Q_y(z)^{\ell+1}}\,dz.
\end{align*}
Since
\[
    \frac{Q_y'(z)}{Q_y(z)^{\ell+1}}
    =
    -\frac1\ell
    \frac{d}{dz}\bigl(Q_y(z)^{-\ell}\bigr),
\]
integration by parts around the closed contour yields
\[
[q^\ell]\Phi(Z_y(q))
=
\frac{1}{2\pi\mathbf i\,\ell}
\oint_{\gamma_1}
\Phi'(z)Q_y(z)^{-\ell}\,dz.
\]
Using
\[
    \Phi'(z)=\frac{1}{z-1+s},
    \qquad
    Q_y(z)^{-1}=U_y(z),
\]
we obtain
\begin{equation}
[q^\ell]\Phi(Z_y(q))
=
\frac{1}{2\pi\mathbf i\,\ell}
\oint_{\gamma_1}
\frac{U_y(z)^\ell}{z-1+s}\,dz.
\label{eq:LB-coefficient}
\end{equation}

Applying the moment formula of Theorem~\ref{t15} with
\[
    j=\ell-1
\]
shows that the right-hand side of
\eqref{eq:LB-coefficient} is
\[
    M_{\ell-1}(\mathbf m_y)
    =
    \int_{\mathbb R}
    \xi^{\ell-1}\,\mathbf m_y(d\xi).
\]
Therefore,
\[
    [q^\ell]\Phi(Z_y(q))
    =
    M_{\ell-1}(\mathbf m_y),
    \qquad \ell\geq1.
\]
Moreover,
\[
    \Phi(Z_y(0))
    =
    \Phi(1)
    =
    0.
\]
Consequently,
\begin{align*}
\Phi(Z_y(q))
&=
\sum_{\ell=1}^{\infty}
M_{\ell-1}(\mathbf m_y)q^\ell=
\sum_{k=0}^{\infty}
M_k(\mathbf m_y)q^{k+1}=
S_{\mathbf m_y}(q).
\end{align*}
The equality first holds as a power-series identity near \(q=0\);
since both sides are holomorphic there, it is also an identity of
holomorphic functions for sufficiently small \(|q|\).

Finally, let \(q=w^{-1}\). For \(|w|\) sufficiently large,
\[
    Z_y(w^{-1})
    =
    \mathbf z_1(w,y).
\]
Using the relation
\[
    \mathrm{St}_{\mathbf m_y}(w)
    =
    S_{\mathbf m_y}(w^{-1})
\]
from Definition~\ref{def:transforms}, we obtain
\begin{align*}
\mathrm{St}_{\mathbf m_y}(w)
&=
S_{\mathbf m_y}(w^{-1})=
\Phi(Z_y(w^{-1}))=
\operatorname{Log}\!
\left(
    \frac{\mathbf z_1(w,y)-1+s}{s}
\right).
\end{align*}
This proves \eqref{eq:st-root-large-w}. The identity extends to every
connected domain on which the branch
\(\mathbf z_1(\,\cdot\,,y)\) and the normalized logarithm admit
holomorphic continuation, by the identity theorem.
\end{proof}

\begin{proposition}[Local existence of the density]
\label{prop:density-existence}
Fix \(y\), and let \(I\subset\mathbb R\) be an open interval.
Let \(\mathbf z_1(w,y)\) be the solution branch introduced in
Proposition~\ref{prop:stieltjes-root}. For \(w\) in the large-\(w\)
domain, set
\[
    G_y(w)
    :=
    \operatorname{Log}\!\left(
        \frac{\mathbf z_1(w,y)-1+s}{s}
    \right),
\]
where the branch of \(\operatorname{Log}\) is the one used in
Proposition~\ref{prop:stieltjes-root}.

Assume that \(G_y\) admits a holomorphic continuation to an open
complex neighborhood \(\Omega\) of \(I\). Then, by
Proposition~\ref{prop:stieltjes-root} and the uniqueness of analytic
continuation,
\[
    \mathrm{St}_{\mathbf m_y}(w)=G_y(w),
    \qquad w\in\Omega\cap\mathbb H.
\]
Then the restriction of \(\mathbf m_y\) to \(I\) is absolutely
continuous and has the continuous density
\[
    f_y(\chi)
    =
    -\frac{1}{\pi}\Im G_y(\chi),
    \qquad \chi\in I.
\]
More precisely,
\[
    -\frac{1}{\pi}
    \Im\mathrm{St}_{\mathbf m_y}(\chi+\mathbf i\delta)
    \longrightarrow
    f_y(\chi)
\]
locally uniformly for \(\chi\in I\) as \(\delta\downarrow0\).
\end{proposition}

\begin{proof}
Let \(K\Subset I\) be compact. Since \(G_y\) is holomorphic in a
neighborhood of \(K\),
\[
    G_y(\chi+\mathbf i\delta)
    \longrightarrow
    G_y(\chi)
\]
uniformly for \(\chi\in K\) as \(\delta\downarrow0\). Hence
\[
    -\frac{1}{\pi}
    \Im\mathrm{St}_{\mathbf m_y}(\chi+\mathbf i\delta)
    \longrightarrow
    -\frac{1}{\pi}\Im G_y(\chi)
\]
locally uniformly on \(I\). Set
\[
    f_y(\chi):=-\frac{1}{\pi}\Im G_y(\chi).
\]
Then \(f_y\) is continuous on \(I\).

For \(\delta>0\), define
\[
    f_{y,\delta}(\chi)
    :=
    -\frac{1}{\pi}
    \Im\mathrm{St}_{\mathbf m_y}(\chi+\mathbf i\delta).
\]
By the definition of the Stieltjes transform,
\[
    f_{y,\delta}(\chi)
    =
    \int_{\mathbb R}
    P_\delta(\chi-\xi)\,\mathbf m_y(d\xi),
\]
where
\[
    P_\delta(u)
    =
    \frac{1}{\pi}\frac{\delta}{u^2+\delta^2}
\]
is the Poisson kernel.

Let \(\varphi\in C_c(I)\). By Fubini's theorem,
\[
    \int_{\mathbb R}
    \varphi(\chi)f_{y,\delta}(\chi)\,d\chi
    =
    \int_{\mathbb R}
    (P_\delta*\varphi)(\xi)\,\mathbf m_y(d\xi).
\]
Since \(P_\delta*\varphi\to\varphi\) uniformly,
\[
    \lim_{\delta\downarrow0}
    \int_{\mathbb R}
    \varphi(\chi)f_{y,\delta}(\chi)\,d\chi
    =
    \int_{\mathbb R}\varphi(\xi)\,\mathbf m_y(d\xi).
\]
On the other hand, the locally uniform convergence
\(f_{y,\delta}\to f_y\) gives
\[
    \lim_{\delta\downarrow0}
    \int_{\mathbb R}
    \varphi(\chi)f_{y,\delta}(\chi)\,d\chi
    =
    \int_I\varphi(\chi)f_y(\chi)\,d\chi.
\]
Therefore,
\[
    \int_I\varphi(\chi)\,\mathbf m_y(d\chi)
    =
    \int_I\varphi(\chi)f_y(\chi)\,d\chi
\]
for every \(\varphi\in C_c(I)\). Hence
\[
    \mathbf m_y|_I=f_y(\chi)\,d\chi.
\]
\end{proof}

\begin{definition}\label{df28}
A point $(\chi,y)\in\mathcal R$ is called a \emph{regular liquid
point} if there is an open interval $I\ni\chi$ and a connected complex
neighborhood $\Omega$ of $I$ satisfying the assumptions of
Proposition~\ref{prop:density-existence}, and if the continuous density
supplied by that proposition satisfies
\[
    0<f_y(\chi)<1.
\]
Let
\[
    \mathcal L
    :=
    \{(\chi,y)\in\mathcal R:(\chi,y)\text{ is a regular liquid point}\}.
\]
By Proposition~\ref{prop:density-existence}, the limiting measure
$\mathbf m_y$ has a continuous density on every horizontal interval
contained in $\mathcal L$. The particle-spacing estimate gives only the
non-strict bound $0\leq f_y\leq1$; the strict inequalities above are part
of the definition. We call $\mathcal L$ the regular liquid region and
call its boundary $\partial\mathcal L$ the frozen boundary.

\end{definition}

\begin{lemma}[Nonreal roots in the regular liquid region]
\label{lem:liquid-nonreal-root}
Let \((\chi,y)\in\mathcal L\), and let
\(\mathbf z_1(\chi,y)\) denote the boundary value at \(w=\chi\) of the
solution branch introduced in
Proposition~\ref{prop:stieltjes-root}. Then
\[
    \Im \mathbf z_1(\chi,y)<0.
\]
Consequently,
\[
    \mathbf z_+(\chi,y)
    :=
    \overline{\mathbf z_1(\chi,y)}
\]
lies in the upper half-plane and satisfies
\[
    U_y(\mathbf z_+(\chi,y))=\chi.
\]
If the equation \(U_y(z)=\chi\) has at most one nonreal
complex-conjugate pair, then
\[
    \bigl\{
        \mathbf z_1(\chi,y),
        \mathbf z_+(\chi,y)
    \bigr\}
\]
is that unique pair.
\end{lemma}

\begin{proof}
By Definition~\ref{df28}, the assumptions of
Proposition~\ref{prop:density-existence} hold in a neighborhood of
\(\chi\), and
\[
    0<f_y(\chi)<1.
\]
Let
\[
    G_y(w)
    :=
    \operatorname{Log}\!\left(
        \frac{\mathbf z_1(w,y)-1+s}{s}
    \right)
\]
be the holomorphic continuation appearing in
Proposition~\ref{prop:density-existence}. The density formula gives
\[
    0
    <
    -\frac{1}{\pi}\Im G_y(\chi)
    <
    1,
\]
and hence
\[
    -\pi<\Im G_y(\chi)<0.
\]
Since
\[
    \frac{\mathbf z_1(\chi,y)-1+s}{s}
    =
    \exp(G_y(\chi)),
\]
we have
\[
    \Im\mathbf z_1(\chi,y)
    =
    s\exp\bigl(\Re G_y(\chi)\bigr)
    \sin\bigl(\Im G_y(\chi)\bigr)
    <0.
\]
Thus \(\mathbf z_1(\chi,y)\) is a lower-half-plane root of
\(U_y(z)=\chi\).

Since \(\mathbf m_0\) is a probability measure supported on
\(\mathbb R\), the function \(H'_{\mathbf m_0}\) satisfies the
Schwarz symmetry
\[
    H'_{\mathbf m_0}(\overline z)
    =
    \overline{H'_{\mathbf m_0}(z)}
\]
on the conjugation-invariant domain under consideration. Since
\(s\in\mathbb R\), it follows that
\[
    U_y(\overline z)
    =
    \overline{U_y(z)}.
\]

Therefore,
\[
\begin{aligned}
    U_y(\mathbf z_+(\chi,y))
    &=
    U_y\bigl(\overline{\mathbf z_1(\chi,y)}\bigr) \\
    &=
    \overline{U_y(\mathbf z_1(\chi,y))} \\
    &=
    \chi.
\end{aligned}
\]
The last assertion follows immediately from the one-pair assumption.
\end{proof}

\begin{theorem}[Multiple-root condition at a horizontal frozen boundary]
\label{t17}
Fix \(y\), and set
\[
    \mathcal L_y
    :=
    \{x\in\mathbb R:(x,y)\in\mathcal L\},
    \qquad
    \mathcal R_y
    :=
    \{x\in\mathbb R:(x,y)\in\mathcal R\}.
\]
Let \(I\) be a connected component of \(\mathcal L_y\), and let
\(\chi\) be a finite endpoint of \(I\) lying in
\(\operatorname{int}(\mathcal R_y)\). Choose a sequence
\[
    \chi_k\in I,
    \qquad
    \chi_k\longrightarrow\chi,
\]
and set
\[
    z_k:=\mathbf z_+(\chi_k,y)\in\mathbb H,
\]
where \(\mathbf z_+\) is defined in
Lemma~\ref{lem:liquid-nonreal-root}.

Assume that, after passing to a subsequence,
\[
    z_k\longrightarrow z_*
\]
for some finite \(z_*\in\mathbb C\), and that \(U_y\) is holomorphic on
a conjugation-invariant open neighborhood of
\(\{z_*,\overline{z_*}\}\). Assume also that, for every real \(x\)
sufficiently close to \(\chi\), the equation
\[
    U_y(z)=x
\]
has at most one nonreal complex-conjugate pair of roots.

Then \(z_*\in\mathbb R\) and
\[
    U_y(z_*)=\chi,
    \qquad
    U_y'(z_*)=0.
\]
Thus \(U_y(z)-\chi\) has a real zero of multiplicity at least two. If,
in addition,
\[
    U_y''(z_*)\neq0,
\]
then the multiplicity is exactly two.
\end{theorem}

\begin{proof}
For every \(k\),
\[
    U_y(z_k)=\chi_k.
\]
Since \(z_k\to z_*\), \(\chi_k\to\chi\), and \(U_y\) is holomorphic near
\(z_*\), we obtain
\[
    U_y(z_*)=\chi.
\]
Moreover, \(z_k\in\mathbb H\), so \(z_*\) lies in the closed upper
half-plane. We first prove that \(z_*\) is real.

Suppose, to the contrary, that \(z_*\in\mathbb H\). We distinguish two
cases.

\medskip
\noindent\emph{Case 1: \(U_y'(z_*)\neq0\).}
By the complex inverse-function theorem, there are a disk \(V\) centered at
\(\chi\), a neighborhood \(W_+\) of \(z_*\), and a holomorphic inverse branch
\[
    \zeta_+:V\longrightarrow W_+
\]
such that
\[
    U_y(\zeta_+(w))=w,
    \qquad
    \zeta_+(\chi)=z_*.
\]
After shrinking \(V\), we may assume that \(V\) is invariant under complex
conjugation. By the Schwarz symmetry assumed in the statement, the function
\[
    \zeta_-(w):=\overline{\zeta_+(\overline w)},
    \qquad w\in V,
\]
is holomorphic and satisfies
\[
    U_y(\zeta_-(w))=w,
    \qquad
    \zeta_-(\chi)=\overline{z_*}.
\]
Shrinking \(V\) once more, we may assume that
\[
    \zeta_+(V\cap\mathbb R)\subset\mathbb H,
    \qquad
    \zeta_-(V\cap\mathbb R)
    \subset\{z\in\mathbb C:\Im z<0\},
\]
and that
\[
    \frac{\zeta_-(w)-1+s}{s}\neq0,
    \qquad w\in V.
\]
For all sufficiently large \(k\), local uniqueness of the inverse branches
gives
\[
    \zeta_+(\chi_k)=z_k,
    \qquad
    \zeta_-(\chi_k)=\overline{z_k}
    =\mathbf z_1(\chi_k,y),
\]
where the last equality follows from
Lemma~\ref{lem:liquid-nonreal-root}.

Fix one such index \(k_0\). Since \((\chi_{k_0},y)\) is a regular liquid
point, Definition~\ref{df28} provides an open interval
\[
    J_{k_0}\ni\chi_{k_0}
\]
and a connected open complex neighborhood \(\Omega_{k_0}\) of \(J_{k_0}\)
to which the logarithmic branch
\[
    G_y(w)
    =
    \operatorname{Log}\!\left(
        \frac{\mathbf z_1(w,y)-1+s}{s}
    \right),
\]
initially defined for sufficiently large \(\lvert w\rvert\), admits a
single-valued holomorphic continuation. Denote this continuation by
\[
    G_{y,k_0}:\Omega_{k_0}\longrightarrow\mathbb C.
\]
By Proposition~\ref{prop:stieltjes-root} and uniqueness of analytic
continuation,
\[
    G_{y,k_0}(w)=\mathrm{St}_{\mathbf m_y}(w),
    \qquad w\in\Omega_{k_0}\cap\mathbb H.
\]

Set
\[
    \widehat\zeta_{y,k_0}(w)
    :=
    1-s+s\exp\bigl(G_{y,k_0}(w)\bigr).
\]
After shrinking \(\Omega_{k_0}\) around \(\chi_{k_0}\), if necessary, the
image of \(\widehat\zeta_{y,k_0}\) lies in a domain on which \(U_y\) is
holomorphic. On a nonempty open subset of
\(\Omega_{k_0}\cap\mathbb H\), Proposition~\ref{prop:stieltjes-root}
shows that
\[
    U_y\bigl(\widehat\zeta_{y,k_0}(w)\bigr)=w.
\]
The identity theorem therefore gives the same equality throughout the
shrunken neighborhood. At \(w=\chi_{k_0}\),
\[
    \widehat\zeta_{y,k_0}(\chi_{k_0})
    =\mathbf z_1(\chi_{k_0},y)
    =\zeta_-(\chi_{k_0}).
\]
Both \(\widehat\zeta_{y,k_0}\) and \(\zeta_-\) are thus local inverse
branches of \(U_y\) taking the same value at \(\chi_{k_0}\). Since this root
is simple, local uniqueness implies that they agree on a connected open
neighborhood of \(\chi_{k_0}\) contained in
\(V\cap\Omega_{k_0}\).

Because \(V\) is simply connected and
\((\zeta_-(w)-1+s)/s\) does not vanish on \(V\), there exists a holomorphic
function
\[
    \widetilde G_y:V\longrightarrow\mathbb C
\]
such that
\[
    \exp\bigl(\widetilde G_y(w)\bigr)
    =
    \frac{\zeta_-(w)-1+s}{s},
    \qquad w\in V.
\]
Choose this logarithmic branch so that it agrees with \(G_{y,k_0}\) on the
open neighborhood of \(\chi_{k_0}\) constructed above. Hence
\(\widetilde G_y\) and \(\mathrm{St}_{\mathbf m_y}\) agree on a nonempty
open subset of \(V\cap\mathbb H\). Since both functions are holomorphic on
the connected set \(V\cap\mathbb H\), the identity theorem yields
\[
    \mathrm{St}_{\mathbf m_y}(w)=\widetilde G_y(w),
    \qquad w\in V\cap\mathbb H.
\]
Thus the Stieltjes-transform branch extends through \(\chi\) to a full
complex neighborhood, not merely from the liquid side.

Since \(\chi\in\operatorname{int}_{\mathbb R}(\mathcal R_y)\), choose an
open interval \(J_0\ni\chi\) such that
\[
    \overline{J_0}\subset V\cap\mathbb R,
    \qquad
    J_0\subset\mathcal R_y.
\]
Proposition~\ref{prop:density-existence}, applied to
\(\widetilde G_y\), shows that \(\mathbf m_y\) is absolutely continuous on
\(J_0\) with continuous density
\[
    f_y(x)
    =
    -\frac{1}{\pi}\Im\widetilde G_y(x),
    \qquad x\in J_0.
\]
For \(x\in I\cap J_0\), we have \((x,y)\in\mathcal L\), and hence
\[
    0<f_y(x)<1.
\]
Therefore,
\[
    -\pi<\Im\widetilde G_y(x)<0,
    \qquad x\in I\cap J_0.
\]
Since \(\chi\) is an endpoint of \(I\), points of \(I\cap J_0\) approach
\(\chi\) from the liquid side. By continuity,
\[
    -\pi\leq\Im\widetilde G_y(\chi)\leq0.
\]
On the other hand,
\[
    \exp\bigl(\widetilde G_y(\chi)\bigr)
    =
    \frac{\zeta_-(\chi)-1+s}{s}
    =
    \frac{\overline{z_*}-1+s}{s}
\]
is nonreal. Hence \(\Im\widetilde G_y(\chi)\) is neither \(0\) nor
\(-\pi\), and thus
\[
    -\pi<\Im\widetilde G_y(\chi)<0.
\]
By continuity, after shrinking to a smaller open interval
\(J\ni\chi\) with \(J\subset J_0\), we have
\[
    0<-\frac{1}{\pi}\Im\widetilde G_y(x)<1,
    \qquad x\in J.
\]
Moreover, \(\zeta_-\) is a holomorphic inverse branch on \(V\), so
\[
    U_y'(\zeta_-(w))\,\zeta_-'(w)=1
\]
and therefore \(U_y'(\zeta_-(w))\neq0\) on \(V\). Thus the same root and
logarithmic branches verify all the regularity conditions in
Definition~\ref{df28} for every \(x\in J\). Consequently,
\[
    J\subset\mathcal L_y.
\]
The interval \(J\) intersects \(I\) and contains the endpoint
\(\chi\notin I\). Hence \(I\cup J\) is a connected subset of
\(\mathcal L_y\) strictly larger than the connected component \(I\), a
contradiction.

\medskip
\noindent\emph{Case 2: \(U_y'(z_*)=0\).}
Let \(m\geq2\) be the multiplicity of \(z_*\) as a zero of
\[
    F(z):=U_y(z)-\chi.
\]
Choose a closed disk \(D_+\Subset\mathbb H\), centered at \(z_*\) and
contained in the holomorphy domain of \(U_y\), such that \(F\) has no zero
on \(\partial D_+\) and no zero in \(D_+\) other than \(z_*\), counted with
multiplicity \(m\). Set
\[
    c:=\min_{z\in\partial D_+}|F(z)|>0.
\]
For every real \(\varepsilon\) with \(|\varepsilon|<c\), Rouch\'e's theorem
shows that
\[
    U_y(z)-(\chi+\varepsilon)=F(z)-\varepsilon
\]
has exactly \(m\) zeros in \(D_+\), counted with multiplicity.

Since \(U_y\) is nonconstant, \(U_y'\) is not identically zero and has only
finitely many zeros in \(\overline{D_+}\). Hence the set
\[
    \mathcal C
    :=
    \bigl\{
        U_y(c_0)-\chi:
        c_0\in\overline{D_+},\ U_y'(c_0)=0
    \bigr\}
    \cap\mathbb R
\]
is finite. We may therefore choose a nonzero real \(\varepsilon\),
arbitrarily small and within the neighborhood in which the one-pair
assumption holds, such that
\[
    |\varepsilon|<c,
    \qquad
    \varepsilon\notin\mathcal C.
\]
For this choice, all \(m\) zeros in \(D_+\) are simple and hence distinct.
Since \(m\geq2\), the equation
\[
    U_y(z)=\chi+\varepsilon
\]
has at least two distinct roots in the upper half-plane. By Schwarz
symmetry, their conjugates produce at least two nonreal complex-conjugate
pairs, contradicting the one-pair assumption.

The two cases exclude \(z_*\in\mathbb H\). Since \(z_*\) lies in the closed
upper half-plane, we conclude that
\[
    z_*\in\mathbb R.
\]

It remains to prove that the limiting real root is multiple. Since
\(z_*\in\mathbb R\),
\[
    \overline{z_k}\longrightarrow z_*.
\]
For each \(k\), the two distinct points \(z_k\) and \(\overline{z_k}\) are
zeros of
\[
    U_y(z)-\chi_k.
\]
Suppose that \(z_*\) were a simple zero of \(U_y(z)-\chi\). Choose a closed
disk \(D\), centered at \(z_*\) and contained in the holomorphy domain of
\(U_y\), such that \(U_y(z)-\chi\) has no zero on \(\partial D\) and exactly
one zero in \(D\), counted with multiplicity. Set
\[
    c_D:=\min_{z\in\partial D}|U_y(z)-\chi|>0.
\]
For all sufficiently large \(k\),
\[
    |\chi_k-\chi|<c_D.
\]
Rouch\'e's theorem then implies that \(U_y(z)-\chi_k\) has exactly one zero
in \(D\), counted with multiplicity. However, for all sufficiently large
\(k\), both \(z_k\) and \(\overline{z_k}\) lie in \(D\), and they are
distinct. This contradiction shows that \(z_*\) is a zero of
\(U_y(z)-\chi\) of multiplicity at least two. Therefore,
\[
    U_y'(z_*)=0.
\]
If \(U_y''(z_*)\neq0\), then the multiplicity is exactly two.
\end{proof}

The linear bottom-boundary model provides an explicit illustration of
Theorem~\ref{t17}.  A complete analysis, including the identification of
the Stieltjes-transform branch and the proof that the resulting
double-root locus is the internal frozen boundary, is given in
Appendix~\ref{app:linear-bottom-boundary}. The corresponding frozen boundary curves for the two weight regimes considered 
are shown in Figure~\ref{fig:fb1}.

\begin{assumption}[Macroscopic tail-weight profile]
\label{ass:tail-profile}
Let $t_n$, $\mathbf x^{(n)}$, and $\alpha$ be as in
Theorem~\ref{t15}. For $y\in[0,\alpha]$, set
\begin{equation}
    s_n(y)
    :=
    \frac{
        \displaystyle\sum_{r=\lfloor ny\rfloor}^{t_n-1}x_r^{(n)}
    }{
        \displaystyle\sum_{r=0}^{t_n-1}x_r^{(n)}
    },
    \label{eq:discrete-tail-profile}
\end{equation}
where an empty sum is interpreted as zero. We shall use the following
nested conditions, listed from weakest to strongest; each condition is
imposed together with all preceding ones.
\begin{enumerate}[label=\textup{(\arabic*)},ref=\textup{(\arabic*)}]
\item\label{ass:tail-profile-convergence}
There exists a continuous function
\[
    s:(0,\alpha)\longrightarrow[0,1]
\]
such that
\[
    s_n(y)\longrightarrow s(y),
    \qquad y\in(0,\alpha),
\]
and
\[
    \lim_{y\downarrow0}s(y)=1,
    \qquad
    \lim_{y\uparrow\alpha}s(y)=0.
\]
\item\label{ass:tail-profile-invertibility}
The function $s$ is strictly decreasing.
\item\label{ass:tail-profile-regularity}
The function $s$ belongs to $C^2((0,\alpha))$ and satisfies
\[
    s'(y)<0,
    \qquad y\in(0,\alpha).
\]
\end{enumerate}
Under condition~\ref{ass:tail-profile-invertibility}, the function $s$ is a
homeomorphism from $(0,\alpha)$ onto $(0,1)$. Denote its inverse by
$y=y(s)$ and set
\begin{equation}
    \vartheta(s):=\alpha-y(s),
    \qquad s\in(0,1).
    \label{eq:theta-tail-profile}
\end{equation}
\end{assumption}

Since the weights are nonnegative, every $s_n$ is nonincreasing.
Consequently, condition~\ref{ass:tail-profile-convergence} and the
continuity of $s$ imply that $s_n\to s$ locally uniformly on
$(0,\alpha)$. In particular, if $y_n\to y\in(0,\alpha)$, then
\[
    s_n(y_n)\longrightarrow s(y).
\]

Under Assumption~\ref{ass:tail-profile}\ref{ass:tail-profile-invertibility}, define
\begin{align}
    \widetilde{\mathcal L}
    :=
    \{(\chi,s(y)):(\chi,y)\in\mathcal L\}.\label{dll}
\end{align}
The $s$-coordinate is used in Section~\ref{gff} because the spectral
equation is naturally written as $V_s(\zeta)=\chi$. Since
$s:(0,\alpha)\to(0,1)$ is a homeomorphism,
\[
    (\chi,y)\longmapsto(\chi,s(y))
\]
is a homeomorphism from $\mathcal L$ onto
$\widetilde{\mathcal L}$. Thus the two regions contain the same liquid
points in different vertical coordinates.

\begin{lemma}[Support in the limiting path domain]
\label{lem:limit-measure-domain-support}
Assume condition~\ref{ass:tail-profile-convergence} of
Assumption~\ref{ass:tail-profile}. If \(y\in(0,\alpha)\) and
\(s(y)\in(0,1)\), then
\[
    \operatorname{supp}(\mathbf m_y)\subset\mathcal R_y,
    \qquad
    \mathcal R_y:=\{\chi\in\mathbb R:(\chi,y)\in\mathcal R\}.
\]
\end{lemma}

\begin{proof}
Choose \(\kappa_n/n\to y\), set
\(\mathbf m_n:=m(\lambda^{(\kappa_n)}(n))\), and choose admissible
slices \(y_n=(\kappa_n+\epsilon_n)/n\to y\). Every atom \(\xi\) of
\(\mathbf m_n\) satisfies
\[
    (\xi,y_n)\in\mathcal R_n.
\]
If \(\chi\notin\mathcal R_y\), choose open intervals \(U\ni\chi\)
and \(V\ni y\) such that
\[
    (\overline U\times\overline V)\cap\mathcal R=\varnothing.
\]
Assumption~\ref{ass:scaled-domain} gives
\(\mathcal R_n\cap(U\times V)=\varnothing\) for all sufficiently
large \(n\), and hence \(\mathbf m_n(U)=0\) almost surely.
The weak convergence in Theorem~\ref{t15} gives
\(\mathbf m_n\Rightarrow\mathbf m_y\) weakly in probability.
Testing against nonnegative functions in \(C_c(U)\) yields
\(\mathbf m_y(U)=0\). This proves the assertion.
\end{proof}

\begin{lemma}[Joint regularity in the liquid region]
\label{lem:joint-liquid-regularity}
Assume
Assumption~\ref{ass:tail-profile}\ref{ass:tail-profile-regularity}.
Then \(\mathcal L\) is open. More precisely, for every
\((\chi_0,y_0)\in\mathcal L\), there exist open intervals
\(J\ni\chi_0\) and \(K\ni y_0\), and a complex neighborhood
\(\Omega\) of \(\overline J\), such that
\[
    J\times K\Subset\mathcal L.
\]
For each \(v\in K\), let
\(\mathbf z_1(\,\cdot\,,v)\) be the physical branch from
Proposition~\ref{prop:stieltjes-root}. These branches and their
normalized logarithms admit compatible holomorphic continuations to
\(\Omega\) such that
\[
    (w,v)\longmapsto\mathbf z_1(w,v)
\]
and
\[
    (w,v)\longmapsto
    \operatorname{Log}\left(
        \frac{\mathbf z_1(w,v)-1+s(v)}{s(v)}
    \right)
\]
are holomorphic in \(w\in\Omega\) and \(C^2\) in \(v\in K\), locally
uniformly in \(w\). Here \(\operatorname{Log}\) is normalized as in
Proposition~\ref{prop:stieltjes-root}.
\end{lemma}

\begin{proof}
Fix \((\chi_0,y_0)\in\mathcal L\) and set \(s_0=s(y_0)\).
By Definition~\ref{df28}, the normalized logarithm from
Proposition~\ref{prop:stieltjes-root} admits analytic continuation
from the large-\(w\) domain to a neighborhood of \(\chi_0\).
Denote this continuation by \(G_0\) and set
\[
    z_0(w):=1-s_0+s_0e^{G_0(w)}.
\]
Proposition~\ref{prop:stieltjes-root} and the identity theorem give
\[
    U_{y_0}(z_0(w))=w.
\]
Differentiating yields
\[
    U_{y_0}'(z_0(w))z_0'(w)=1,
\]
and hence \(U_{y_0}'(z_0(w))\neq0\).

Cover a continuation path from the large-\(w\) domain to
\(\chi_0\) by finitely many overlapping disks. Since the expression
in \eqref{duyz} is holomorphic in \(z\) and in the parameter \(s\),
the parameter-dependent holomorphic implicit-function theorem,
applied successively along these disks, extends the branch for
\(s\) near \(s_0\). Uniqueness on the overlaps and agreement on the
large-\(w\) disk identify the resulting branch with
\(\mathbf z_1\) from Proposition~\ref{prop:stieltjes-root}. The
normalized logarithm extends simultaneously. Since
\(s\in C^2((0,\alpha))\), these continuations are holomorphic in
\(w\) and \(C^2\) in \(v\), locally uniformly in \(w\).

At \((\chi_0,y_0)\), Definition~\ref{df28} gives
\[
    0<
    -\frac1\pi
    \Im\operatorname{Log}\left(
        \frac{\mathbf z_1(\chi_0,y_0)-1+s(y_0)}{s(y_0)}
    \right)
    <1.
\]
By continuity, there exist open intervals \(J\ni\chi_0\) and
\(K\ni y_0\), and a complex neighborhood \(\Omega\) of
\(\overline J\), such that
\begin{equation}
    0<
    -\frac1\pi
    \Im\operatorname{Log}\left(
        \frac{\mathbf z_1(x,v)-1+s(v)}{s(v)}
    \right)
    <1,
    \qquad (x,v)\in J\times K.
    \label{eq:joint-liquid-density-bounds}
\end{equation}
Proposition~\ref{prop:density-existence} and
Lemma~\ref{lem:limit-measure-domain-support} then give
\[
    J\subset\operatorname{supp}(\mathbf m_v)
    \subset\mathcal R_v,
    \qquad v\in K.
\]
Thus Definition~\ref{df28} implies
\(J\times K\subset\mathcal L\). Hence \(\mathcal L\) is open, and
shrinking \(J\) and \(K\) once more gives
\(J\times K\Subset\mathcal L\).
\end{proof}

\section{Rescaled Height Function and Complex Burgers Equation}\label{sect:be}

In this section, we prove that the slopes of the rescaled height
function are encoded by a complex slope variable. Under
Assumption~\ref{ass:tail-profile}\ref{ass:tail-profile-regularity}, this variable satisfies a weighted
complex Burgers equation whose source term is determined by the
macroscopic tail-weight profile. For uniform weights the profile is
affine, the source term vanishes, and the resulting corollary gives the
standard inviscid complex Burgers equation conjectured in
Conjecture~6.1 of \cite{SKN21}.

On the lecture hall graph $\mathcal{G}_t$, define  the discrete height function $h^{\mathrm{disc}}$ associated to a  non-intersecting path configuration as follows. The height at the lower left corner is 0, and the height increases by 1 whenever crossing a path from the left to the right. Define the rescaled height function by 
 
  \begin{equation}
    h_n(\chi,y)
    :=
    \frac{1}{n}h^{\mathrm{disc}}(n\chi,ny).
    \label{eq:rescaled-height}
\end{equation}

\begingroup
Fix \(y\in(0,\alpha)\), and choose integers
\(\kappa_n\in\{0,1,\ldots,t_n\}\) such that the hypotheses of
Theorem~\ref{t15} are satisfied. In particular,
\[
    \frac{\kappa_n}{n}\longrightarrow y,
    \qquad
    \frac{|\mathbf x_{\kappa_n}^{(n)}|}
         {|\mathbf x^{(n)}|}
    \longrightarrow s\in(0,1).
\]
Let
\[
    \mathbf m_n
    :=
    m\bigl(\lambda^{(\kappa_n)}(n)\bigr)
\]
be the random counting measure in Theorem~\ref{t15}. Choose
\(\epsilon_n\in(0,1)\) as in the construction of
\(\lambda^{(\kappa_n)}(n)\), and set
\[
    y_n:=\frac{\kappa_n+\epsilon_n}{n}.
\]
Then \(y_n\to y\). For all sufficiently large \(n\), the horizontal
coordinates of the \(n\) paths at this slice are
\[
    \lambda_i^{(\kappa_n)}(n)+n-i,
    \qquad 1\leq i\leq n.
\]
Since the discrete height is normalized to be zero to the left of all
paths and increases by one when a path is crossed from left to right,
we have, with the natural constant extension outside the path domain,
\begin{equation}
    h_n(\chi,y_n)
    =F_n(\chi)
    :=\mathbf m_n(({-}\infty,\chi])
    \quad\text{for Lebesgue-a.e. }\chi\in\mathbb R.
    \label{eq:height-as-cdf}
\end{equation}
The values at the finitely many jump points depend only on the chosen
step-function convention and play no role below.

\begin{proposition}[Slice-wise law of large numbers for the height]
\label{prop:slice-height-lln}
Let \(\mathbf m_y\) be the deterministic probability measure supplied
by Theorem~\ref{t15}, and define
\begin{equation}
    \mathbf h(\chi,y)
    :=F_y(\chi)
    :=\mathbf m_y(({-}\infty,\chi]),
    \qquad \chi\in\mathbb R.
    \label{eq:canonical-slice-height}
\end{equation}
Then
\begin{equation}
    h_n(\,\cdot\,,y_n)
    \longrightarrow
    \mathbf h(\,\cdot\,,y)
    \quad\text{in probability in }L^1(\mathbb R).
    \label{eq:limiting-height-profile}
\end{equation}
 Moreover,
\[
    D_\chi\mathbf h(\,\cdot\,,y)=\mathbf m_y
\]
in the sense of distributions.
\end{proposition}

\begin{proof}

By Theorem~\ref{t15}, there exists \(B<\infty\) such that
\(\mathbf m_n\) and \(\mathbf m_y\) are supported in \([0,B]\), and
\[
    \mathbf m_n\Longrightarrow\mathbf m_y
    \qquad\text{weakly in probability}.
\]

Since all the measures are supported in \([0,B]\), the weak
convergence implies
\[
    \int_{\mathbb R}|F_n(\chi)-F_y(\chi)|\,d\chi
    \longrightarrow0
    \qquad\text{in probability}.
\]
Indeed, along every almost-surely weakly convergent subsequence,
\(F_n(\chi)\to F_y(\chi)\) at every continuity point of \(F_y\), and
the assertion follows from dominated convergence.
Together with \eqref{eq:height-as-cdf} and
\eqref{eq:canonical-slice-height}, this proves
\eqref{eq:limiting-height-profile}.

Finally, for every
\(\varphi\in C_c^\infty(\mathbb R)\), the definition of
\(\mathbf h\) as the cumulative distribution function of
\(\mathbf m_y\) gives
\[
    -\int_{\mathbb R}
        \mathbf h(\chi,y)\varphi'(\chi)\,d\chi
    =
    \int_{\mathbb R}
        \varphi(\chi)\,\mathbf m_y(d\chi).
\]
This is precisely
\[
    D_\chi\mathbf h(\,\cdot\,,y)=\mathbf m_y
\]
in the sense of distributions.
\end{proof}

Under
Assumption~\ref{ass:tail-profile}\ref{ass:tail-profile-convergence},
the local uniform convergence of the tail profiles gives, whenever
\(\kappa_n/n\to y\),
\[
    \frac{|\mathbf x_{\kappa_n}^{(n)}|}{|\mathbf x^{(n)}|}
    =s_n\left(\frac{\kappa_n}{n}\right)
    \longrightarrow s(y).
\]
Thus the proposition applies at every \(y\) for which
\(s(y)\in(0,1)\). Its limit is independent of the chosen sequence of
levels because the moment formula in Theorem~\ref{t15} depends on the
level sequence only through \(s(y)\), and compactly supported
probability measures are determined by their moments. In particular,
under
Assumption~\ref{ass:tail-profile}\ref{ass:tail-profile-regularity},
the profile \(\mathbf h(\chi,y)\) is defined for every
\(y\in(0,\alpha)\).  We call \[
    \mathbf h:\mathbb R\times(0,\alpha)\longrightarrow[0,1],
\]
the limiting rescaled height profile.
\endgroup

\begin{proposition}[Two-dimensional law of large numbers for the height]
\label{prop:two-dimensional-height-lln}
Assume
Assumption~\ref{ass:tail-profile}\ref{ass:tail-profile-regularity}.
Then, for every compact interval \(K\Subset(0,\alpha)\),
\begin{equation}
    \lVert h_n-\mathbf h\rVert_{L^1(\mathbb R\times K)}
    \longrightarrow0
    \qquad\text{in probability}.
    \label{eq:two-dimensional-height-lln}
\end{equation}
Consequently,
\[
    h_n\longrightarrow\mathbf h
    \qquad\text{in probability in }
    L^1_{\mathrm{loc}}
    \bigl(\mathbb R\times(0,\alpha)\bigr).
\]
We call \(\mathbf h\) the deterministic limit shape.
\end{proposition}

\begin{proof}
Fix \(y\in(0,\alpha)\). For all sufficiently large \(n\), set
\[
    \kappa_n^\pm:=\lfloor ny\rfloor\pm1,
    \qquad
    y_n^\pm
    :=
    \frac{\kappa_n^\pm+\epsilon_n^\pm}{n},
    \qquad
    0<\epsilon_n^\pm<1,
\]
where \(y_n^\pm\) are admissible slices. Then
\[
    y_n^-<y<y_n^+,
    \qquad
    y_n^\pm\longrightarrow y.
\]
Since every path moves only downward and to the right,
\[
    h_n(\,\cdot\,,y_n^-)
    \leq h_n(\,\cdot\,,y)
    \leq h_n(\,\cdot\,,y_n^+)
    \qquad\text{a.e. on }\mathbb R.
\]

For each sign,
\[
    \frac{\kappa_n^\pm}{n}\longrightarrow y,
    \qquad
    \frac{|\mathbf x_{\kappa_n^\pm}^{(n)}|}
         {|\mathbf x^{(n)}|}
    \longrightarrow s(y),
\]
where the second convergence follows from local uniform convergence
of the tail profiles. Applying
Proposition~\ref{prop:slice-height-lln} separately to the two level
sequences gives
\[
    h_n(\,\cdot\,,y_n^\pm)
    \longrightarrow
    \mathbf h(\,\cdot\,,y)
    \quad\text{in probability in }L^1(\mathbb R).
\]
Therefore, with
\[
    X_n(y)
    :=
    \lVert
        h_n(\,\cdot\,,y)-\mathbf h(\,\cdot\,,y)
    \rVert_{L^1(\mathbb R)},
\]
the preceding inequalities yield
\[
\begin{aligned}
    X_n(y)
    \leq{}&
    \lVert
        h_n(\,\cdot\,,y_n^-)-\mathbf h(\,\cdot\,,y)
    \rVert_{L^1(\mathbb R)}
    \\
    &+
    \lVert
        h_n(\,\cdot\,,y_n^+)-\mathbf h(\,\cdot\,,y)
    \rVert_{L^1(\mathbb R)}
    \longrightarrow0
\end{aligned}
\]
in probability.

With \(B\) as in Theorem~\ref{t15}, both height functions take values
in \([0,1]\) and agree outside \([0,B]\). Hence
\[
    0\leq X_n(y)\leq B
\]
uniformly in \(n\) and \(y\), and consequently
\[
    \mathbb E X_n(y)\longrightarrow0.
\]

It remains to justify integration in \(y\). If \(y_r\to y\), the
continuity of \(s\) and the moment formula in Theorem~\ref{t15} give
convergence of every moment of \(\mathbf m_{y_r}\) to the corresponding
moment of \(\mathbf m_y\). Their common compact support then implies
\[
    \mathbf m_{y_r}\Longrightarrow\mathbf m_y.
\]
Thus \(y\mapsto\mathbf m_y\) is weakly continuous and
\(\mathbf h\) is Borel measurable. For every compact
\(K\Subset(0,\alpha)\), Tonelli's theorem and dominated convergence
now give
\[
\begin{aligned}
    \mathbb E
    \lVert h_n-\mathbf h\rVert_{L^1(\mathbb R\times K)}
    &=
    \int_K\mathbb E X_n(y)\,dy
    \longrightarrow0.
\end{aligned}
\]
The assertion follows from Markov's inequality.
\end{proof}

Let \((\chi,y)\) be a regular liquid point, and let \(I\ni\chi\)
be an interval supplied by Definition~\ref{df28}.  The preceding
distributional identity and Proposition~\ref{prop:density-existence}
give
\[
    D_\chi\mathbf h(\,\cdot\,,y)
    =
    \mathbf m_y
    =
    f_y(\chi)\,d\chi
    \qquad\text{on }I.
\]
Since \(f_y\) is continuous on \(I\), the function
\(\mathbf h(\,\cdot\,,y)\) has a \(C^1\) representative on \(I\).
For this representative,
\[
    \frac{\partial\mathbf h(\chi,y)}{\partial\chi}
    =
    f_y(\chi).
\]
Moreover, Proposition~\ref{prop:density-existence} gives
\[
    f_y(\chi)
    =
    -\lim_{\delta\downarrow0}
    \frac{1}{\pi}
    \Im\mathrm{St}_{\mathbf m_y}
    (\chi+\mathbf i\delta).
\]
By Lemma~\ref{lem:liquid-nonreal-root},
\[
    \mathbf z_+(\chi,y)
    =
    \overline{\mathbf z_1(\chi,y)}
\]
is the corresponding upper-half-plane root.  Hence
\begin{equation}
    \frac{\partial\mathbf h(\chi,y)}{\partial\chi}
    =
    f_y(\chi)
    =
    -\lim_{\delta\downarrow0}
    \frac{1}{\pi}
    \Im\mathrm{St}_{\mathbf m_y}
    (\chi+\mathbf i\delta)
    =
    \frac{1}{\pi}
    \operatorname{Arg}
    \bigl(\mathbf z_+(\chi,y)-1+s\bigr),
    \label{hnx}
\end{equation}
where \(\operatorname{Arg}\) takes values in \((0,\pi)\) on the upper
half-plane.

For the remainder of this section, assume
Assumption~\ref{ass:tail-profile}\ref{ass:tail-profile-regularity} and write
\[
    s=s(y).
\]
In particular, \(s\in C^2((0,\alpha))\) and \(s'(y)<0\).

\begingroup
Assume
\begin{eqnarray*}
\lim_{n\rightarrow\infty}
\frac{\lambda_1^{(0)}(n)+n-1}{n}
=\beta\in[1,\infty).
\end{eqnarray*}
For every partition in the support of \(\rho_{\kappa_n}^{(n)}\),
\[
    0\leq
    \frac{\lambda_i^{(\kappa_n)}(n)+n-i}{n}
    \leq
    \frac{\lambda_1^{(0)}(n)+n-1}{n},
    \qquad 1\leq i\leq n.
\]

The weak convergence in Theorem~\ref{t15} and the endpoint assumption
therefore imply
\endgroup
\[
    \operatorname{supp}(\mathbf m_y)\subset[0,\beta],
    \qquad y\in(0,\alpha).
\]

For $j\geq0$, write
\[
    M_j(y)
    :=
    \int_{\mathbb R}u^j\,\mathbf m_y(du).
\]
By Theorem~\ref{t15} and
Assumption~\ref{ass:tail-profile}\ref{ass:tail-profile-regularity},
each function $M_j$ belongs to $C^1((0,\alpha))$. Since all the
measures $\mathbf m_y$ are supported in the common compact interval
$[0,\beta]$, polynomial approximation shows that
\[
    y\longmapsto\mathbf m_y
\]
is weakly continuous. By
Proposition~\ref{prop:slice-height-lln} and the support inclusion
\(\operatorname{supp}(\mathbf m_y)\subset[0,\beta]\), the canonical
limiting height profile in \eqref{eq:canonical-slice-height} can be
written as
\begin{equation}
    \mathbf h(\chi,y)=
    \begin{cases}
        0, & \chi\leq0,\\
        \mathbf m_y([0,\chi]), & 0<\chi<\beta,\\
        1, & \chi\geq\beta.
    \end{cases}
    \label{eq:canonical-height-profile}
\end{equation}
The values at the jump points are immaterial. By
Proposition~\ref{prop:slice-height-lln}, this representative is the
slice-wise \(L^1(\mathbb R)\)-limit in probability in
\eqref{eq:limiting-height-profile} and satisfies
\[
    D_\chi\mathbf h(\,\cdot\,,y)=\mathbf m_y.
\]
It is also a bounded Borel function on
\(\mathbb R\times(0,\alpha)\). Consequently, its distributional derivative
\(D_y\mathbf h\) is well defined and is supported in \([0,\beta]\) in the
horizontal variable.

For every integer $j\geq1$, Stieltjes integration by parts gives
\begin{align}
M_j(y)
&=
\int_0^\beta u^j\,d_\chi\mathbf h(u,y)\notag\\
&=
\beta^j-j\int_0^\beta
\mathbf h(u,y)u^{j-1}\,du.
\label{pdr}
\end{align}
Taking the distributional derivative of \eqref{pdr} with respect to
$y$ yields
\begin{equation}
    \int_0^\beta
    u^{j-1}D_y\mathbf h(u,y)\,du
    =
    -\frac1j M_j'(y)
    \qquad\text{in }\mathcal D'((0,\alpha)).
    \label{eq:weak-vertical-moments}
\end{equation}
The right-hand side is a continuous function of $y$.

Fix an open interval $K$ with
\[
    \overline K\Subset(0,\alpha).
\]
Choose $r>0$ such that
\[
    r<\min_{y\in\overline K}s(y)
\]
and such that $H'_{\mathbf m_0}$ is analytic on
\[
    \{z\in\mathbb C:|z-1|\leq r\}.
\]
In the contour formula below, take
\[
    \gamma_1=\{z\in\mathbb C:|z-1|=r\}.
\]
All distributional identities in the following contour computation are
understood in $\mathcal D'(K)$.

Applying Theorem~\ref{t15} and differentiating its contour formula gives
\begin{align*}
\int_0^\beta
u^{j-1}D_y\mathbf h(u,y)\,du
&=
-\frac{1}{2j(j+1)\pi\mathbf i}
\frac{d}{dy}
\oint_{\gamma_1}
\frac{U_y(z)^{j+1}}{z-1+s}
\,dz.
\end{align*}
Making the change of variables
\[
    w=z-1+s,
\]
let
\[
    \gamma_s:=\{w\in\mathbb C:|w-s|=r\},
\]
the positively oriented image of $\gamma_1$. Set
\[
    A_s(w):=wH_{\mathbf m_0}'(w+1-s)+\frac{w}{w-s},
    \qquad
    B_s(w):=H_{\mathbf m_0}'(w+1-s)+\frac{1}{w-s}.
\]
Then $A_s(w)=wB_s(w)$ and
$\partial_sA_s(w)=-w\partial_wB_s(w)$. Hence, in
$\mathcal D'(K)$,
\begin{align}
\int_0^\beta
u^{j-1}D_y\mathbf h(u,y)\,du
&=
-\frac{1}{2j(j+1)\pi\mathbf i}
s'(y)
\oint_{\gamma_s}
\frac{1}{w}
\frac{\partial}{\partial s}A_s(w)^{j+1}
\,dw\notag\\
&=
\frac{s'(y)}{2j\pi\mathbf i}
\oint_{\gamma_s}A_s(w)^jB_s'(w)\,dw\notag\\
&=
\frac{s'(y)}{2j\pi\mathbf i}
\oint_{\gamma_s}A_s(w)^j\,dB_s(w).
\label{eq:weak-vertical-contour-moments}
\end{align}

For $x\in\mathbb C\setminus[0,\beta]$, define the distributional
Stieltjes transform of $D_y\mathbf h$ by
\begin{equation}
    \left\langle
        \mathcal S(x,\,\cdot\,),\psi
    \right\rangle
    :=
    \left\langle
        D_y\mathbf h(u,y),
        \frac{\psi(y)}{x-u}
    \right\rangle,
    \qquad
    \psi\in C_c^\infty((0,\alpha)).
    \label{eq:distributional-stieltjes-transform}
\end{equation}
We now restrict this distribution to $K$. In the remainder of this
computation, $x$ denotes the complex Stieltjes variable, $u$ denotes the
real horizontal variable, and $w$ denotes the contour variable. Recall
the function $A_s$ defined above. Since $\overline K$ is compact and the
family of contours varies continuously with $y$,
\[
    M_K
    :=
    \sup_{\substack{y\in\overline K\\ w\in\gamma_{s(y)}}}
    |A_{s(y)}(w)|
    <\infty.
\]
Assume
\[
    |x|>\max\{\beta,2M_K\}.
\]
The expansion
\[
    \frac{1}{x-u}
    =
    \sum_{j=1}^{\infty}\frac{u^{j-1}}{x^j}
\]
converges in $C^\infty$ on a neighborhood of $[0,\beta]$. Therefore,
by \eqref{eq:weak-vertical-moments},
\begin{align}
\mathcal S(x,y)
&=
\sum_{j=1}^{\infty}x^{-j}
\int_0^\beta u^{j-1}D_y\mathbf h(u,y)\,du
\label{stx}\\
&=
\sum_{j=1}^{\infty}
\frac{s'(y)}{2j\pi\mathbf i}
\oint_{\gamma_s}
\left(\frac{A_s(w)}{x}\right)^j
 d\left(
 H_{\mathbf m_0}'(w+1-s)
 +\frac{1}{w-s}
\right)
\notag\\
&=
-\frac{s'(y)}{2\pi\mathbf i}
\oint_{\gamma_s}
\log\left(1-\frac{A_s(w)}{x}\right)
 d\left(
 H_{\mathbf m_0}'(w+1-s)
 +\frac{1}{w-s}
\right)
\notag\\
&=
-\frac{s'(y)}{2\pi\mathbf i}
\oint_{\gamma_s}
 d\left[
 \left(
 H_{\mathbf m_0}'(w+1-s)
 +\frac{1}{w-s}
 \right)
 \log\left(1-\frac{A_s(w)}{x}\right)
 \right]
\notag\\
&\quad+
\frac{s'(y)}{2\pi\mathbf i}
\oint_{\gamma_s}
\frac{
 H_{\mathbf m_0}'(w+1-s)+\dfrac{1}{w-s}
}{
 1-\dfrac{A_s(w)}{x}
}
\frac{\partial}{\partial w}
\left(1-\frac{A_s(w)}{x}\right)
\,dw.
\notag
\end{align}
The equalities initially hold in $\mathcal D'(K)$; the final contour
expression is continuous in $y\in K$, so it represents
$\mathcal S(x,\,\cdot\,)$ as an ordinary function on $K$.
By (\ref{hmd}) we obtain
\begin{eqnarray}
H_{\mathbf{m}_0}'(w+1-s)+\frac{1}{w-s}=\frac{1}{(w+1-s) S_{\bm_0}^{(-1)}(\ln (w+1-s))}.\label{hpe}
\end{eqnarray}

\begin{lemma}\label{l31}
With $K$ and $r$ fixed as above, there exists $R_K>0$ such that, for
every $y\in K$, with $s=s(y)$, and every $x$ satisfying $|x|>R_K$, the
equation
\[
    A_s(w)=x
\]
has exactly one solution inside $\gamma_s$, and $A_s$ has exactly one
pole inside $\gamma_s$.
\end{lemma}

\begin{proof}
Recall that whenever a formal power series of the form
\[
    f(z)=\sum_{i=1}^{\infty}f_i z^i\in\RR[[z]]
\]
has $f_1\neq0$, there exists a unique composition inverse
\[
    g(u)=\sum_{i=1}^{\infty}g_i u^i\in\RR[[u]].
\]
From the explanations leading to \eqref{hmd},
\[
    S_{\mathbf m_0}^{(-1)}(u)
    =
    u+\sum_{i=2}^{\infty}s_i u^i.
\]
It follows that
\begin{align}
\frac{1}{
(w+1-s)S_{\bm_0}^{(-1)}(\ln(w+1-s))
}
=
\frac{1}{w-s}
+\sum_{i=0}^{\infty}c_i(w-s)^i.
\label{sip}
\end{align}
By \eqref{hpe}, $A_s$ has exactly one pole, at $w=s$, inside the disk
bounded by $\gamma_s$. Since $s(\overline K)$ is a compact subset of
$(0,1)$, the Laurent expansion above and the corresponding Rouch\'e
estimates are uniform for $y\in K$. Hence there exists $R_K>0$ such
that, whenever $y\in K$ and $|x|>R_K$, the equation $A_s(w)=x$ has
exactly one solution in that disk. This proves the lemma.
\end{proof}

To compute $\mathcal S(x,y)$, we first assume that
\[
    |x|>\max\{\beta,2M_K,R_K\}.
\]
It follows from Lemma~\ref{l31} that
\begin{align}
-\frac{1}{2\pi\mathbf i}
\oint_{\gamma_s}
 d\left[
 \left(
 H_{\mathbf m_0}'(w+1-s)
 +\frac{1}{w-s}
 \right)
 \log\left(1-\frac{A_s(w)}{x}\right)
 \right]
=0.
\label{vpt}
\end{align}
Under Assumption~\ref{ass:tail-profile}\ref{ass:tail-profile-regularity}, the parameter in the moment formula
is $s=s(y)$. It follows from \eqref{stx} and \eqref{vpt} that
\begin{align}
\mathcal S(x,y)
=
\frac{s'(y)}{2\pi\mathbf i}
\oint_{\gamma_s}
\frac{
 H_{\mathbf m_0}'(w+1-s)+\dfrac{1}{w-s}
}{
 1-\dfrac{A_s(w)}{x}
}
\frac{\partial}{\partial w}
\left(1-\frac{A_s(w)}{x}\right)
\,dw.
\label{str}
\end{align}
Recall also that $B_s(w)=H'_{\mathbf m_0}(w+1-s)+(w-s)^{-1}$ and
$A_s(w)=wB_s(w)$.
Since
\[
    -\frac{B_s(w)}{1-A_s(w)/x}
    \frac{\partial}{\partial w}
    \left(
        1-\frac{A_s(w)}{x}
    \right)
    =
    \frac{B_s(w)A_s'(w)}{x-A_s(w)},
\]
equation \eqref{str} becomes
\begin{equation}
    \mathcal S(x,y)
    =
    -\frac{s'(y)}{2\pi\mathbf i}
    \oint_{\gamma_s}
    \frac{B_s(w)A_s'(w)}
         {x-A_s(w)}
    \,dw.
    \label{eq:st-xi-residue-form}
\end{equation}

For sufficiently large \(\lvert x\rvert\), let
\[
    w_x
    :=
    \mathbf z_1(x,y)+s-1.
\]
Then
\[
    A_s(w_x)=x,
\]
and \(w_x\) is the unique simple solution of this equation inside
\(\gamma_s\).

We first compute the residue at \(w=s\). By \eqref{sip}, writing
\(t=w-s\), we have
\[
    B_s(w)
    =
    \frac1t+c_0+c_1t+O(t^2)
\]
for some constants \(c_0,c_1\). Since
\[
    A_s(w)=(s+t)B_s(w),
\]
we obtain
\[
    A_s(w)
    =
    \frac{s}{t}
    +
    1+sc_0
    +
    (c_0+sc_1)t
    +
    O(t^2),
\]
and hence
\[
    B_s(w)A_s'(w)
    =
    -\frac{s}{t^3}
    -
    \frac{sc_0}{t^2}
    +
    O(t^{-1}).
\]
Moreover,
\[
    \frac1{x-A_s(w)}
    =
    -\frac{t}{s}
    -
    \frac{x-1-sc_0}{s^2}t^2
    +
    O(t^3).
\]
Multiplying these expansions gives
\[
    \frac{B_s(w)A_s'(w)}
         {x-A_s(w)}
    =
    \frac1{t^2}
    +
    \frac{x-1}{s}\frac1t
    +
    O(1).
\]
Therefore,
\[
    \operatorname*{Res}_{w=s}
    \frac{B_s(w)A_s'(w)}
         {x-A_s(w)}
    =
    \frac{x-1}{s}.
\]

At the simple zero \(w=w_x\), we have
\[
    x-A_s(w)
    =
    -A_s'(w_x)(w-w_x)
    +
    O\bigl((w-w_x)^2\bigr),
\]
and therefore
\[
    \operatorname*{Res}_{w=w_x}
    \frac{B_s(w)A_s'(w)}
         {x-A_s(w)}
    =
    -B_s(w_x).
\]
The residue theorem now yields
\begin{align}
\mathcal S(x,y)
&=
-s'(y)
\left(
    \frac{x-1}{s}
    -
    B_s(w_x)
\right)\notag\\
&=
-s'(y)\frac{x-1}{s}
+
 s'(y)
\left(
    H'_{\mathbf m_0}
    \bigl(\mathbf z_1(x,y)\bigr)
    +
    \frac{1}{\mathbf z_1(x,y)-1}
\right)\notag\\
&=
-s'(y)\frac{x-1}{s}
+
\frac{s'(y)}{
    \mathbf z_1(x,y)
    S_{\mathbf m_0}^{(-1)}
    \bigl(\log\mathbf z_1(x,y)\bigr)
}.
\label{eq:st-xi-root}
\end{align}
The last equality follows from \eqref{hpe}.

The identity is first obtained for sufficiently large
$\lvert x\rvert$. Let $\psi\in C_c^\infty(K)$. On every connected
$x$-domain in
$\mathbb C\setminus[0,\beta]$ to which the branch
$\mathbf z_1(\,\cdot\,,y)$ admits holomorphic continuation for
$y\in K$, both the pairing
$\langle\mathcal S(x,\,\cdot\,),\psi\rangle$ and the pairing of the
right-hand side of \eqref{eq:st-xi-root} with $\psi$ are holomorphic
in $x$. They agree for large $\lvert x\rvert$ and hence agree throughout
the domain by the identity theorem. Thus \eqref{eq:st-xi-root} holds as
an identity of distributions on $K$; since its right-hand side is
continuous in $y$, it is the ordinary-function representative of
$\mathcal S(x,\,\cdot\,)$ on $K$. Since the open interval $K$ with
$\overline K\Subset(0,\alpha)$ was arbitrary,
\eqref{eq:st-xi-root} holds in $\mathcal D'((0,\alpha))$.

Let $(\chi_0,y_0)\in\mathcal L$. By
Lemma~\ref{lem:joint-liquid-regularity}, after shrinking the
neighborhoods if necessary, there exist open intervals
$J\ni\chi_0$ and $K_0\ni y_0$ such that
\[
    J\times K_0\Subset\mathcal L,
\]
and the physical branch $\mathbf z_1(x,y)$ is holomorphic in $x$ in
a complex neighborhood of $J$ and $C^2$ in $y\in K_0$.

Equation~\eqref{eq:st-xi-root} therefore gives the locally uniform
boundary limit
\begin{align}
\lim_{\epsilon\downarrow0}
\mathcal S(\chi+\mathbf i\epsilon,y)
&=
-s'(y)\frac{\chi-1}{s}
+
\frac{s'(y)}{
    \mathbf z_1(\chi,y)
    S_{\mathbf m_0}^{(-1)}
    \bigl(\log\mathbf z_1(\chi,y)\bigr)
},
\qquad (\chi,y)\in J\times K_0.
\label{eq:st-xi-boundary-weighted}
\end{align}
For $\epsilon>0$, set
\[
    g_\epsilon(\chi,y)
    :=
    -\frac{1}{\pi}
    \Im\mathcal S(\chi+\mathbf i\epsilon,y).
\]
By the definition \eqref{eq:distributional-stieltjes-transform},
\[
    g_\epsilon
    =
    P_\epsilon*_{\chi}D_y\mathbf h
\]
in the sense of distributions, where
\[
    P_\epsilon(v)
    :=
    \frac{1}{\pi}
    \frac{\epsilon}{v^2+\epsilon^2}
\]
is the Poisson kernel and the convolution is taken only in the
horizontal variable. Hence
\[
    g_\epsilon
    \longrightarrow
    D_y\mathbf h
\]
in $\mathcal D'(J\times K_0)$ as $\epsilon\downarrow0$. On the other
hand, \eqref{eq:st-xi-boundary-weighted} gives locally uniform
convergence to
\[
    g(\chi,y)
    :=
    -\frac{1}{\pi}
    \Im\left[
        \frac{s'(y)}{
            \mathbf z_1(\chi,y)
            S_{\mathbf m_0}^{(-1)}
            \bigl(\log\mathbf z_1(\chi,y)\bigr)
        }
    \right].
\]
Therefore
\begin{equation}
    D_y\mathbf h(\chi,y)
    =
    -\frac{1}{\pi}
    \Im\left[
        \frac{s'(y)}{
            \mathbf z_1(\chi,y)
            S_{\mathbf m_0}^{(-1)}
            \bigl(\log\mathbf z_1(\chi,y)\bigr)
        }
    \right]
    \qquad\text{in }\mathcal D'(\mathcal L).
    \label{eq:weak-vertical-density}
\end{equation}

Hence we have the following theorem.
\begin{theorem}
\label{thm:m31}
Suppose that the data $t_n$, $\mathbf x^{(n)}$, and
$\lambda^{(0)}(n)$ satisfy the
standing hypotheses of Theorem~\ref{t15}; in particular, the
bottom-boundary sequence is regular in the sense of
Definition~\ref{def:regular-bottom}. Suppose
Assumption~\ref{ass:tail-profile}\ref{ass:tail-profile-regularity} holds. Assume also that
\[
    \lim_{n\to\infty}
    \frac{\lambda_1^{(0)}(n)+n-1}{n}
    =
    \beta
    \in[1,\infty),
\]
and that the regular liquid region $\mathcal L$ is nonempty.

Under these assumptions, let $\mathbf h$ be the deterministic limit
shape from Proposition~\ref{prop:two-dimensional-height-lln}; by
\eqref{eq:canonical-height-profile}, it is represented by the
cumulative distribution functions of the measures \(\mathbf m_y\).

For each $(\chi,y)\in\mathcal L$, let
$\mathbf z_+(\chi,y)\in\mathbb H$ be the upper-half-plane root
defined in Lemma~\ref{lem:liquid-nonreal-root}, and set
\begin{equation}
    u(\chi,y)
    :=
    \frac{s'(y)}{
        \mathbf z_+(\chi,y)
        S_{\mathbf m_0}^{(-1)}
        \bigl(\log\mathbf z_+(\chi,y)\bigr)
    }.
    \label{eq:weighted-complex-slope}
\end{equation}
Here $\log$ denotes the restriction to $\mathbb H$ of the principal
logarithm on $\mathbb C\setminus(-\infty,0]$, so that
$0<\Im\log z<\pi$ for $z\in\mathbb H$. The notation
$S_{\mathbf m_0}^{(-1)}$ refers to the local compositional inverse
introduced in Definition~\ref{def:transforms}, analytically continued
from its local definition near $0$ so that \eqref{hpe} holds.

Then
\[
    u(\chi,y)\in\mathbb H,
    \qquad
    (\chi,y)\in\mathcal L.
\]
Moreover, $\mathbf h$ has a $C^1$ representative on $\mathcal L$,
which we continue to denote by $\mathbf h$, and
\begin{equation}
    \frac{\partial\mathbf h}{\partial\chi}
    =
    1-\frac{1}{\pi}\operatorname{Arg}(u),
    \qquad
    \frac{\partial\mathbf h}{\partial y}
    =
    \frac{1}{\pi}\Im u,
    \label{s21}
\end{equation}
where $\operatorname{Arg}$ takes values in $(0,\pi)$ on
$\mathbb H$.

Finally, $u$ is locally $C^1$ on $\mathcal L$ and satisfies the
weighted complex Burgers equation
\begin{equation}
    u_y+u\,u_\chi
    =
    \frac{s''(y)}{s'(y)}\,u.
    \label{be-weighted}
\end{equation}
\end{theorem}

\begin{proof}
Fix $(\chi,y)\in\mathcal L$, write $s=s(y)$, and set
\[
    z:=\mathbf z_+(\chi,y).
\]
By Lemma~\ref{lem:liquid-nonreal-root},
\[
    z\in\mathbb H,
    \qquad
    U_y(z)=\chi.
\]

Lemma~\ref{lem:joint-liquid-regularity} shows that
$\mathbf z_+(\chi,y)$ is locally $C^2$ as a function of
$(\chi,y)$. Consequently, the function $u$ in
\eqref{eq:weighted-complex-slope} is locally $C^1$.

We first identify the horizontal slope. Equation \eqref{hnx} gives
\begin{equation}
    \frac{\partial\mathbf h}{\partial\chi}
    =
    \frac{1}{\pi}
    \operatorname{Arg}(z-1+s).
    \label{eq:horizontal-slope-weighted}
\end{equation}

Define
\begin{equation}
    Q(z)
    :=
    \frac{1}{
        zS_{\mathbf m_0}^{(-1)}(\log z)
    }
    =
    H'_{\mathbf m_0}(z)+\frac{1}{z-1},
    \label{qzh}
\end{equation}
where the second equality follows from \eqref{hmd}. Since
$\mathbf m_0$ is supported on $\mathbb R$, the relevant analytic
branches satisfy Schwarz symmetry:
\[
    Q(\overline z)=\overline{Q(z)}.
\]
By the definition of $u$,
\begin{equation}
    u=s'(y)Q(z).
    \label{qza}
\end{equation}
Moreover, Lemma~\ref{lem:liquid-nonreal-root} gives
$\mathbf z_1(\chi,y)=\overline z$. Therefore,
\[
    s'(y)Q\bigl(\mathbf z_1(\chi,y)\bigr)
    =
    s'(y)Q(\overline z)
    =
    \overline u.
\]
Equation~\eqref{eq:weak-vertical-density} and
$\mathbf z_1(\chi,y)=\overline z$ give
\begin{equation}
    D_y\mathbf h
    =
    -\frac{1}{\pi}\Im\overline u
    =
    \frac{1}{\pi}\Im u
    \qquad\text{in }\mathcal D'(\mathcal L).
    \label{eq:weak-vertical-slope}
\end{equation}

We next determine the half-plane containing $u$. By
\eqref{duyz}, \eqref{qzh}, and $U_y(z)=\chi$,
\[
    \chi=(z-1+s)Q(z).
\]
Together with \eqref{qza}, this yields
\begin{equation}
    (z-1+s)u=s'(y)\chi.
    \label{eq:zu-chi-weighted}
\end{equation}
The measure $\mathbf m_y$ is supported in $[0,\beta]$. Since
$(\chi,y)$ is a regular liquid point, the density of
$\mathbf m_y$ is positive on an open interval containing $\chi$;
hence $0<\chi<\beta$. Also, $z-1+s\in\mathbb H$, while
$s'(y)\chi<0$. It follows from
\eqref{eq:zu-chi-weighted} that $u\in\mathbb H$ and
\[
    \operatorname{Arg}(u)
    =
    \pi-\operatorname{Arg}(z-1+s).
\]
Combining this identity with
\eqref{eq:horizontal-slope-weighted} proves the first identity in
\eqref{s21}.

Equations~\eqref{hnx} and \eqref{eq:weak-vertical-slope} show that
the two distributional derivatives of $\mathbf h$ on $\mathcal L$ are
represented by the continuous functions
\[
    1-\frac{1}{\pi}\operatorname{Arg}(u)
    \qquad\text{and}\qquad
    \frac{1}{\pi}\Im u.
\]
A standard distributional regularity argument, obtained for example by
local mollification, now implies that $\mathbf h$ has a $C^1$
representative on $\mathcal L$ with these partial derivatives. Thus the
two identities in \eqref{s21} hold pointwise for this representative.

It remains to prove \eqref{be-weighted}. Put
\[
    q(\chi,y):=Q(z(\chi,y))=\frac{u(\chi,y)}{s'(y)}.
\]
Since $q$ is a function of $z$ alone,
\begin{equation}
    q_\chi z_y-q_yz_\chi=0.
    \label{eq:q-parallel-gradients}
\end{equation}
The root equation can be written as
\begin{equation}
    (z-1+s)q=\chi.
    \label{eq:zq-root}
\end{equation}
Differentiating \eqref{eq:zq-root} with respect to $\chi$ and $y$
gives
\[
    qz_\chi+(z-1+s)q_\chi=1,
\]
and
\[
    q\bigl(z_y+s'(y)\bigr)+(z-1+s)q_y=0.
\]
Multiplying the first identity by $q_y$, the second by $q_\chi$,
and subtracting the first resulting identity from the second, then using
\eqref{eq:q-parallel-gradients}, gives
\begin{equation}
    q_y+s'(y)q\,q_\chi=0.
    \label{eq:q-burgers}
\end{equation}
Since $u=s'q$, equation \eqref{eq:q-burgers} implies
\begin{align*}
    u_y+u\,u_\chi
    &=
    s''q+s'q_y+(s'q)(s'q_\chi)\\
    &=
    s''q+s'\bigl(q_y+s'q\,q_\chi\bigr)\\
    &=
    s''q
    =
    \frac{s''}{s'}u.
\end{align*}
This proves \eqref{be-weighted}.
\end{proof}

\begin{corollary}[Uniform weights]
\label{cor:uniform-burgers}
Assume all hypotheses of Theorem~\ref{thm:m31} except
Assumption~\ref{ass:tail-profile}\ref{ass:tail-profile-regularity}, and
suppose that, for every $n$,
\begin{equation}
    x_r^{(n)}=1,
    \qquad
    0\leq r\leq t_n-1.
    \label{uw}
\end{equation}
Then
\[
    s(y)=1-\frac{y}{\alpha},
    \qquad
    s'(y)=-\frac{1}{\alpha},
    \qquad
    s''(y)=0.
\]
Consequently,
\[
    u(\chi,y)
    =
    -\frac{1}{
        \alpha\,
        \mathbf z_+(\chi,y)
        S_{\mathbf m_0}^{(-1)}
        \bigl(\log\mathbf z_+(\chi,y)\bigr)
    }
\]
satisfies the slope identities \eqref{s21} and the standard complex
Burgers equation
\begin{equation}
    u_y+u\,u_\chi=0.
    \label{be}
\end{equation}
\end{corollary}

\begin{proof}
Under \eqref{uw}, the discrete tail profile in
\eqref{eq:discrete-tail-profile} is
\[
    s_n(y)
    =
    \frac{(t_n-\lfloor ny\rfloor)_+}{t_n},
\]
and it converges uniformly on $[0,\alpha]$ to
$1-y/\alpha$. The conclusion follows from
Theorem~\ref{thm:m31}.
\end{proof}

\begin{remark}
With the height-function and complex-slope conventions used here,
Corollary~\ref{cor:uniform-burgers} proves Conjecture~6.1 of
\cite{SKN21}. The specialization \eqref{uw} is precisely the uniformly
sampled model considered in \cite{SKN21}, rather than an additional
restriction on that conjecture.
\end{remark}

\section{Height fluctuations and the Gaussian free field (GFF) in weighted proportional regimes}\label{gff}

In this section, we prove convergence of the unrescaled height
fluctuations to the Gaussian free field in the weighted proportional
regime under Assumption~\ref{ass:tail-profile}\ref{ass:tail-profile-invertibility} and Assumption~\ref{ap32}. The
main idea is as follows. First, we apply the Schur difference operator
from Section~\ref{lsti} to the Schur generating functions, obtain the
moments of the height observables, and verify Wick's formula in the
scaling limit. Second, we construct an explicit diffeomorphism from the
liquid region in $(\chi,s)$-coordinates to the upper half-plane and
show that the limiting covariance is the Green kernel there. The
uniformly weighted proportional regime is recorded separately as a
corollary.

The main theorem in Section \ref{gff} is Theorem \ref{thm:gff}.

Let $C_c^{\infty}(\mathbb H)$ denote the space of smooth real-valued
functions with compact support in the upper half-plane, and let $dA$ denote
planar Lebesgue measure. The \textbf{zero-boundary Gaussian free field}
(GFF) on $\mathbb H$ is the centered Gaussian random distribution $\Xi$
characterized by
\begin{equation}
    \mathbb E[\Xi(f_1)\Xi(f_2)]
    =
    \int_{\mathbb H}\int_{\mathbb H}
    f_1(z)G_{\mathbb H}(z,w)f_2(w)\,dA(z)\,dA(w),
    \label{eq:gff-covariance}
\end{equation}
for $f_1,f_2\in C_c^{\infty}(\mathbb H)$, where
\begin{equation}
    G_{\mathbb H}(z,w)
    :=
    -\frac{1}{2\pi}
    \log\left|\frac{z-w}{z-\overline w}\right|,
    \qquad z,w\in\mathbb H,
    \label{eq:green-H}
\end{equation}
is the Green function of the Dirichlet Laplacian on $\mathbb H$.

We shall also use the standard extension of the GFF to finite signed Borel
measures of finite Green energy. If
\[
    \int_{\mathbb H}\int_{\mathbb H}
    |G_{\mathbb H}(z,w)|\,|\nu|(dz)\,|\nu|(dw)
    <\infty,
\]
then $\nu$ has finite Green energy and $\Xi(\nu)$ is defined by
completion from smooth compactly supported densities. For two such
measures,
\begin{equation}
    \mathbb E[\Xi(\nu_1)\Xi(\nu_2)]
    =
    \int_{\mathbb H}\int_{\mathbb H}
    G_{\mathbb H}(z,w)\,\nu_1(dz)\,\nu_2(dw).
    \label{eq:gff-measure-covariance}
\end{equation}

The GFF is conformally invariant in law. More precisely, if
$\mathcal D\subsetneq\mathbb C$ is simply connected and
$\phi:\mathcal D\to\mathbb H$ is conformal, then the random distribution
on $\mathcal D$ defined by
\[
    (\phi^*\Xi)(f)
    :=
    \Xi\left(
        (f\circ\phi^{-1})\left|(\phi^{-1})'\right|^2
    \right),
    \qquad f\in C_c^\infty(\mathcal D),
\]
has the law of the zero-boundary GFF on $\mathcal D$. Its joint moments
are determined by Wick's rule and the covariance kernel
$G_{\mathbb H}$. See \cite{SS07} for further background.

 Let $f$ be a function of $r$ variables. Define the symmetrization of $f$ as
 \begin{align}
     \mathrm{Sym}_{v_1,\ldots,v_r}f(v_1,\ldots,v_r):=\frac{1}{r!}\sum_{\sigma\in S_r}f(v_{\si(1)},\ldots,v_{\si(r)}).
     \label{dsym}
 \end{align}

For $\theta\in(0,\alpha)$, $k\geq1$, and all sufficiently large $n$,
set
\begin{equation}
    p_k^{(\lfloor\theta n\rfloor)}
    :=
    \sum_{i=1}^{n}
    \left(
        \lambda_i^{(t_n-\lfloor\theta n\rfloor)}+n-i
    \right)^k.
    \label{defpk}
\end{equation}
For fixed $g\geq1$ and
$\alpha_1,\ldots,\alpha_g\in(0,\alpha)$, we use
\eqref{defpk} with $\theta=\alpha_r$, $r=1,\ldots,g$. Powers such as
$[p_k^{(\lfloor\theta n\rfloor)}]^m$ denote ordinary powers of this
random variable.

  The following lemma is a multilevel extension of Lemma~\ref{l04}.
\begin{lemma}Let $g\geq1$, let $\alpha_1,\ldots,\alpha_g\in(0,\alpha)$, and let $k_1,\ldots,k_g$ and $m_1,\ldots,m_g$ be positive integers. For $1\leq r\leq g$, define $\kappa_r$ by
\begin{align}
\lfloor\alpha_r n\rfloor+\kappa_r=t_n.\label{kal}
\end{align}
Assume that
\begin{align*}
\kappa_1\geq\kappa_2\geq\ldots\geq\kappa_g.
\end{align*}
Then
  
  \begin{align}
\mathbb{E}\prod_{r=1}^{g}\left[p_{k_r}^{(\lfloor\alpha_r n\rfloor)}\right]^{m_r}=\left.\mathcal{D}_{k_1,\kappa_1}^{m_1}\cdot\ldots\cdot\mathcal{D}_{k_g,\kappa_g}^{m_g}\mathcal{S}_{\rho_{\kappa_g}}(|\mathbf{x}_{\kappa_g}|,\mathbf{u})\right|_{\mathbf{u}=0}.\label{mmt}
\end{align} 
 \end{lemma}

 \begin{proof} When $g=1$, (\ref{mmt}) follows from Lemma \ref{l04}. We now prove the $g=2$ case. The case $g\geq3$ is obtained by iterating the same conditioning argument over the nested levels $\kappa_1\geq\cdots\geq\kappa_g$. If two adjacent levels coincide, the corresponding conditional distribution is a point mass, so the same argument applies without change. Applying \eqref{eq:schur-operator-action} at level $\kappa_2$, we have
\begin{align}
&\mathcal{D}_{k_1,\kappa_1}^{m_1}\mathcal{D}_{k_2,\kappa_2}^{m_2}\mathcal{S}_{\rho_{\kappa_2}}(|\mathbf{x}_{\kappa_2}|,\mathbf{u})=\mathcal{D}_{k_1,\kappa_1}^{m_1}\left[\sum_{\lambda\in \YY}\rho_{\kappa_2}(\lambda)
\left[\sum_{i=1}^n(\lambda_i+n-i)^{k_2}\right]^{m_2}\frac{s_{\lambda}(|\mathbf{x}_{\kappa_2}|+\mathbf{u})}{s_{\lambda}(|\mathbf{x}_{\kappa_2}|)}\right].\label{dd1}
 \end{align}
 By (\ref{srk}), one has
  
  \begin{align}
\frac{s_{\lambda}(|\mathbf{x}_{\kappa_2}|+\mathbf{u})}{s_{
\lambda}(|\mathbf{x}_{\kappa_2}|)}
&=\sum_{\mu\in \YY}\rho_{\kappa_1}(\mu\mid\lambda^{(\kappa_2)}=\lambda)\frac{s_{\mu}(|\mathbf{x}_{\kappa_1}|+\mathbf{u})}{s_{\mu}(|\mathbf{x}_{\kappa_1}|)}.\label{cpm}
\end{align} 
 Here $\rho_{\kappa_1}(\mu\mid\lambda^{(\kappa_2)}=\lambda)$ denotes the conditional distribution of the partition at level $\kappa_1$ given the partition at level $\kappa_2$.  It follows from (\ref{dd1}) and (\ref{cpm}) that
 
  \begin{align*}
&\mathcal{D}_{k_1,\kappa_1}^{m_1}\mathcal{D}_{k_2,\kappa_2}^{m_2}\mathcal{S}_{\rho_{\kappa_2}}(|\mathbf{x}_{\kappa_2}|,\mathbf{u})\\
&=\mathcal{D}_{k_1,\kappa_1}^{m_1}\left[\sum_{\lambda\in \YY}\rho_{\kappa_2}(\lambda)
\left[\sum_{i=1}^n(\lambda_i+n-i)^{k_2}\right]^{m_2}\sum_{\mu\in \YY}\rho_{\kappa_1}(\mu\mid\lambda^{(\kappa_2)}=\lambda)\frac{s_{\mu}(|\mathbf{x}_{\kappa_1}|+\mathbf{u})}{s_{\mu}(|\mathbf{x}_{\kappa_1}|)}\right]\\
&=\sum_{\lambda,\mu\in\YY}
\mathbb P\bigl(
    \lambda^{(\kappa_2)}=\lambda,
    \lambda^{(\kappa_1)}=\mu
\bigr)
\left[\sum_{i=1}^n(\lambda_i+n-i)^{k_2}\right]^{m_2}\\
&\qquad\times
\left[\sum_{j=1}^n(\mu_j+n-j)^{k_1}\right]^{m_1}
\frac{s_{\mu}(|\mathbf{x}_{\kappa_1}|+\mathbf{u})}
     {s_{\mu}(|\mathbf{x}_{\kappa_1}|)}
\end{align*} 
where the last identity follows by applying \eqref{eq:schur-operator-action} at level $\kappa_1$ to the conditional distribution in \eqref{cpm}. Then, when $g=2$, \eqref{mmt} follows by letting $\mathbf{u}=\mathbf{0}$.
 \end{proof}

\begin{theorem}
\label{t31}
Suppose that the data $t_n$, $\mathbf x^{(n)}$, and
$\lambda^{(0)}(n)$ satisfy the standing hypotheses of
Theorem~\ref{t15}.  In particular, the bottom-boundary sequence is
regular in the sense of Definition~\ref{def:regular-bottom}.  Fix $g\geq1$ and
$\alpha_1,\ldots,\alpha_g\in(0,\alpha)$, and assume that the limits
\begin{equation}
    s_r
    :=
    \lim_{n\to\infty}
    \frac{
        |\mathbf x^{(n)}_{t_n-\lfloor\alpha_r n\rfloor}|
    }{
        |\mathbf x^{(n)}|
    }
    \in(0,1),
    \qquad r=1,\ldots,g,
    \label{eq:section4-tail-limits}
\end{equation}
exist. For $r\in[g]$ and $k\geq1$, let
$p_k^{(\lfloor\alpha_rn\rfloor)}$ be defined by \eqref{defpk}.

Every finite subcollection of
\begin{equation}
    \left\{
        n^{-k}\left(
            p_k^{(\lfloor\alpha_rn\rfloor)}
            -\mathbb E p_k^{(\lfloor\alpha_rn\rfloor)}
        \right)
    \right\}_{r\in[g],\,k\geq1}
    \label{crv}
\end{equation}
converges, in the sense of moments, to a centered Gaussian vector. Its
covariance is given as follows. For $r_1,r_2\in[g]$ and
$k_1,k_2\geq1$,
{
interchange the two pairs if necessary so that
$\alpha_{r_1}\leq\alpha_{r_2}$.  Because the weights are nonnegative,
this ordering also gives $s_{r_1}\leq s_{r_2}$. }
\begin{align}
&\Sigma_{(r_1,k_1),(r_2,k_2)}
:=
\lim_{n\to\infty}
\frac{
    \operatorname{Cov}\left(
        p_{k_1}^{(\lfloor\alpha_{r_1}n\rfloor)},
        p_{k_2}^{(\lfloor\alpha_{r_2}n\rfloor)}
    \right)
}{n^{k_1+k_2}}
\label{covf}\\
&\quad=
\frac{1}{(2\pi\mathbf i)^2}
\oint_{|z-1|=\epsilon}
\oint_{|w-1|=2\epsilon}
\left(
    (z-1+s_{r_1})H_{\mathbf m_0}'(z)
    +\frac{z-1+s_{r_1}}{z-1}
\right)^{k_1}\notag\\
&\qquad\qquad\times
\left(
    (w-1+s_{r_2})H_{\mathbf m_0}'(w)
    +\frac{w-1+s_{r_2}}{w-1}
\right)^{k_2}
Q(z,w)\,dz\,dw,\notag
\end{align}
where
\begin{align*}
Q(z,w)
&=
\frac{1}{(z-w)^2}\\
&\quad+
\frac{\partial^2}{\partial z\,\partial w}
\log\left(
    1-
    \frac{(z-1)(w-1)}{z-w}
    \left[zH_{\mathbf m_0}'(z)-wH_{\mathbf m_0}'(w)\right]
\right).
\end{align*}
Here $\epsilon>0$ is chosen so that the closed disk
$|z-1|\leq2\epsilon$ lies in the common domain of analyticity of the
functions in the integrand and contains no singularity other than the
displayed pole at $1$. Both contours are positively oriented; their
radii make them disjoint and nested.
\end{theorem}

\begin{definition}\label{def:global-F}
For use throughout this section, write
\[
S_n(\mathbf u):=\frac{s_{\lambda^{(0)}}\left(1^n+\frac{\mathbf u}{|\mathbf x|}\right)}{s_{\lambda^{(0)}}(1^n)},
\qquad
u_{j,\kappa}:=\frac{u_j}{{|\mathbf x_\kappa|}},
\]
and define
\begin{align}
\mathcal F_{\kappa,k}(\mathbf u)
:=\frac{1}{V(\mathbf u)S_n(\mathbf u)}
\left[\sum_j\left((1+u_{j,\kappa})\frac{\partial}{\partial u_{j,\kappa}}\right)^k\right]
V(\mathbf u)S_n(\mathbf u).
\label{def:global-F-formula}
\end{align}
This definition is global in the present section and is not confined to a proof.
\end{definition}

\begingroup

\begin{definition}[Derivative-stable $n$-degree]
\label{def:derivative-stable-degree}
Let \(A_n(\mathbf z)\) be analytic in a neighborhood of
\(\mathbf z=\mathbf1^n\). For $q\geq1$, set
\[
    J_q(A_n)
    :=
    \max_{i_1,\ldots,i_q\in[n]}
    \left|
        \left.
        \partial_{z_{i_1}}\cdots\partial_{z_{i_q}}
        A_n(\mathbf z)
        \right|_{\mathbf z=1^n}
    \right|,
    \qquad
    J_0(A_n):=|A_n(\mathbf1^n)|,
\]
where the coordinate indices \(i_1,\ldots,i_q\) need not be distinct. For an integer $D$, we write
\[
    \deg_n A_n\leq D
    \quad\text{if}\quad
    J_q(A_n)=O(n^D)
\]
for every fixed $q$, and
\[
    \deg_n A_n<D
    \quad\text{if}\quad
    J_q(A_n)=o(n^D)
\]
for every fixed $q$. For every fixed \(q\), the \(O\)- and \(o\)-estimates are uniform over
the coordinate indices \(i_1,\ldots,i_q\). When the displayed function itself carries
coordinate indices, such as $\partial_i\mathcal F_{a,n}$, the estimate
is also required to hold uniformly over all permitted values of those
indices. 
\end{definition}

By definition, applying any fixed number of coordinate derivatives
preserves the degree bound. By the product
rule, degrees add under multiplication, and the product has strict
degree if one factor does. A sum of $O(n^r)$ terms satisfying an
$O(n^D)$, respectively $o(n^D)$, bound uniformly in the summation
indices has degree at most $D+r$, respectively less than $D+r$.
\endgroup

\begingroup

For use in the one- and two-operator estimates below, put
\begin{align}
    b_n:=|\mathbf x^{(n)}|,
    \qquad
    \widehat S_n(\mathbf z)
    :=
    \frac{s_{\lambda^{(0)}(n)}(z_1,\ldots,z_n)}
         {s_{\lambda^{(0)}(n)}(1^n)},
    \qquad
    L_n:=\log\widehat S_n,\label{dsnh}
\end{align}
and
\[
    V_n(\mathbf z):=\prod_{i<j}(z_i-z_j),
    \qquad
    z_i:=1+\frac{u_i}{b_n},
    \qquad
    \partial_i:=\frac{\partial}{\partial z_i}.
\]
For an observable indexed by \(a\), set
\begin{align}
    s_{a,n}:=\frac{|\mathbf x^{(n)}_{\kappa_a}|}{b_n},
    \qquad
    g_{a,n}(z):=z-1+s_{a,n},
    \qquad
    \mathscr D_{a,n}
    :=
    \sum_{i=1}^n
    \bigl(g_{a,n}(z_i)\partial_i\bigr)^{k_a}.\label{dsag}
\end{align}
Then
\begin{equation}
    (1+u_{i,\kappa_a})
    \frac{\partial}{\partial u_{i,\kappa_a}}
    =g_{a,n}(z_i)\partial_i,
    \label{eq:affine-euler-change}
\end{equation}
and the normalized one-operator function is
\begin{equation}
    \mathcal F_{a,n}(\mathbf z)
    :=
    \frac{
        \mathscr D_{a,n}
        \bigl(V_n(\mathbf z)\widehat S_n(\mathbf z)\bigr)
    }{
        V_n(\mathbf z)\widehat S_n(\mathbf z)
    }
    =
    \mathcal F_{\kappa_a,k_a}
    \bigl(b_n(\mathbf z-\mathbf1)\bigr).
    \label{eq:affine-one-operator}
\end{equation}
Introduce the notation
\begin{align}
\mathcal F^{\mathrm{lead}}_{a,n}
&:=
\sum_{q=0}^{k_a}\binom{k_a}{q}(q+1)!
\sum_{\substack{B\subset[n]\\ |B|=q+1}}
\operatorname{Sym}_{b_1,\ldots,b_{q+1}}
\frac{
    g_{a,n}(z_{b_1})^{k_a}
    (\partial_{b_1}L_n)^{k_a-q}
}{
    \prod_{j=2}^{q+1}(z_{b_1}-z_{b_j})
}.
\label{eq:correct-inner-leading}
\end{align}
Here and below an empty product is one.

\begin{lemma}[One-operator expansion and bounds]\label{adfd}
Suppose \(s_{a,n}\to s_a\in(0,1)\). In the sense of
Definition~\ref{def:derivative-stable-degree},
\[
\begin{array}{lll}
\deg_n\mathcal F_{a,n}\leq k_a+1,
&\deg_n(\partial_i\mathcal F_{a,n})\leq k_a,
&\deg_n(\partial_i\partial_j\mathcal F_{a,n})<k_a
\quad(i\ne j),
\end{array}
\]
while
\[
    \deg_n\mathcal F^{\mathrm{lead}}_{a,n}\leq k_a+1,
    \qquad
    \deg_n(\partial_i\mathcal F^{\mathrm{lead}}_{a,n})\leq k_a.
\]
Moreover, uniformly in \(i\in[n]\),
\begin{equation}
    \deg_n\partial_i\left(
        \mathcal F_{a,n}-\mathcal F^{\mathrm{lead}}_{a,n}
    \right)<k_a.
    \label{eq:inner-leading-remainder}
\end{equation}
Each displayed estimate controls all further fixed-order derivatives
of the displayed function.
\end{lemma}

\begin{proof}
Lemmas~\ref{la1} and~\ref{la3}, their local uniformity, and Cauchy's
formula give the following bounds for every fixed mixed derivative of
\(L_n\): a derivative supported on one coordinate has degree at most
\(1\), one supported on exactly two distinct coordinates has degree at
most \(0\), and one supported on at least three distinct coordinates
has degree less than \(0\). Repeated indices are allowed. Because the
functions are symmetric, these bounds are uniform over all choices of
the coordinate indices.

Apply the product rule in \eqref{eq:affine-one-operator} and collect
terms involving the same coordinate indices. By the argument of
\cite[Lemma~5.2]{bg16}, symmetrization cancels the apparent poles at
\(z_i=z_j\). The resulting functions are analytic near
\(\mathbf z=\mathbf1^n\), and the same Taylor-coefficient argument
preserves the uniform \(O/o\) degree bounds used here.

We first estimate the expression
\(\mathcal F^{\mathrm{lead}}_{a,n}\) in
\eqref{eq:correct-inner-leading}. Fix \(q\). For each set
\(B\subset[n]\) with \(|B|=q+1\), the factor
\[
    (\partial_{b_1}L_n)^{k_a-q}
\]
has degree at most \(k_a-q\), while \(g_{a,n}^{k_a}\) has degree
zero. The symmetrized quotient does not increase this degree. Since
there are
\[
    \binom{n}{q+1}=O(n^{q+1})
\]
possible choices of \(B\), the \(q\)-th sum in
\eqref{eq:correct-inner-leading} has degree at most
\[
    (k_a-q)+(q+1)=k_a+1.
\]
After applying \(\partial_i\), the \(q\)-th sum has degree at most
\(k_a\). Indeed, if \(i\in B\), fixing \(i\) reduces the number of
choices of \(B\) from \(O(n^{q+1})\) to \(O(n^q)\). If \(i\notin B\),
the derivative can act only through \(L_n\), producing a mixed
derivative and again lowering the degree by one.

These are the only terms that can attain degree \(k_a\). A term in the
full product-rule expansion that is not included in
\(\mathcal F^{\mathrm{lead}}_{a,n}\) either uses one of the
differentiations on \(g_{a,n}\), producing no logarithmic-derivative
factor, or combines two or more differentiations into a higher
derivative of \(L_n\). In either case it loses at least one power of
\(n\). Here we use
\[
    g_{a,n}'=1,
    \qquad
    g_{a,n}^{(q)}=0\quad(q\geq2),
\]
together with the fact that \(s_{a,n}\) remains in a compact subset of
\((0,1)\). This proves \eqref{eq:inner-leading-remainder}.

The same counting gives
\[
    \deg_n\mathcal F_{a,n}\leq k_a+1,
    \qquad
    \deg_n(\partial_i\mathcal F_{a,n})\leq k_a.
\]
If \(i\ne j\), applying both \(\partial_i\) and \(\partial_j\) forces
one further loss, and hence
\[
    \deg_n(\partial_i\partial_j\mathcal F_{a,n})<k_a.
\]

Finally, the argument remains valid after any fixed number of further
mixed derivatives: for each fixed derivative order, the same index
count is uniform over all choices of the derivative indices. Thus all
the estimates hold in the sense of
Definition~\ref{def:derivative-stable-degree}.
\end{proof}
\endgroup

\begingroup

For a symmetric analytic germ \(H_n\) at \(\mathbf z=\mathbf1^n\), set
\begin{align}
    \mathscr T_{a,n}[H_n]
    :=
    \frac{
        \mathscr D_{a,n}
        \bigl(V_n\widehat S_n H_n\bigr)
    }{
        V_n\widehat S_n
    },
    \qquad
    \mathscr I_{a,n}[H_n]
    :=
    \mathscr T_{a,n}[H_n]-\mathcal F_{a,n}H_n.\label{dti}
\end{align}
Also define
\begin{align}
\mathscr J_{a,n}[H_n]
&:=
k_a\sum_{p=0}^{k_a-1}
\binom{k_a-1}{p}(p+1)!
\sum_{\substack{A\subset[n]\\ |A|=p+1}}
\operatorname{Sym}_{a_1,\ldots,a_{p+1}}
\frac{
    g_{a,n}(z_{a_1})^{k_a}
    (\partial_{a_1}L_n)^{k_a-1-p}
    \partial_{a_1}H_n
}{
    \prod_{j=2}^{p+1}(z_{a_1}-z_{a_j})
}.
\label{eq:affine-interaction-leading}
\end{align}

\begin{lemma}[Affine interaction expansion]
\label{lem:affine-interaction}
Suppose that \(s_{a,n}\to s_a\in(0,1)\) and, for some integer \(d\),
\[
    \deg_n(\partial_iH_n)\leq d
\]
uniformly in \(i\in[n]\). Then
\[
    \deg_n\mathscr I_{a,n}[H_n]\leq k_a+d,
    \qquad
    \deg_n\mathscr J_{a,n}[H_n]\leq k_a+d,
\]
and
\begin{equation}
    \deg_n\bigl(
        \mathscr I_{a,n}[H_n]
        -\mathscr J_{a,n}[H_n]
    \bigr)<k_a+d.
    \label{eq:affine-interaction-remainder}
\end{equation}
If \(\widetilde H_n\) is another symmetric analytic germ satisfying
\[
    \deg_n\partial_i(H_n-\widetilde H_n)<d
\]
uniformly in \(i\), then
\begin{equation}
    \deg_n\bigl(
        \mathscr J_{a,n}[H_n]
        -\mathscr J_{a,n}[\widetilde H_n]
    \bigr)<k_a+d.
    \label{eq:affine-interaction-replacement}
\end{equation}
All the displayed estimates are derivative-stable in the sense of
Definition~\ref{def:derivative-stable-degree}.
\end{lemma}

\begin{proof}
In the product-rule expansion of \(\mathscr T_{a,n}[H_n]\), the terms
on which no derivative falls on \(H_n\) sum exactly to
\(\mathcal F_{a,n}H_n\) and therefore cancel from
\(\mathscr I_{a,n}[H_n]\). After grouping equal coordinate sets and
symmetrizing, every remaining term is a constant multiple of
\[
\sum_{\substack{A\subset[n]\\ |A|=p+1}}
\operatorname{Sym}_{a_1,\ldots,a_{p+1}}
\frac{
    g_{a,n}(z_{a_1})^{k_a-m_0}
    \partial_{a_1}^{m_1}H_n
    \prod_{\ell\geq1}
        (\partial_{a_1}^{\ell}L_n)^{d_\ell}
}{
    \prod_{j=2}^{p+1}(z_{a_1}-z_{a_j})
},
\]
where \(m_0,m_1,p,d_1,d_2,\ldots\) are nonnegative integers, only
finitely many \(d_\ell\) are nonzero, \(m_1\geq1\), and
\begin{equation}
    m_0+m_1+p+\sum_{\ell\geq1}\ell d_\ell=k_a.
    \label{eq:affine-interaction-budget}
\end{equation}
Since \(m_1\geq1\), the assumption
\(\deg_n(\partial_iH_n)\leq d\), together with
Definition~\ref{def:derivative-stable-degree}, gives
\[
    \deg_n\bigl(\partial_{a_1}^{m_1}H_n\bigr)\leq d.
\] Hence its degree is at most
\[
    d+\sum_{\ell\geq1}d_\ell+(p+1).
\]
Here the last term comes from the \(O(n^{p+1})\) choices of the
coordinate set \(A\).
By \eqref{eq:affine-interaction-budget}, the difference between
\(k_a+d\) and this bound is
\[
    m_0+(m_1-1)
    +\sum_{\ell\geq2}(\ell-1)d_\ell.
\]
It can vanish only when
\[
    m_0=0,
    \qquad m_1=1,
    \qquad d_\ell=0\quad(\ell\geq2),
\]
in which case \(d_1=k_a-1-p\). These are precisely the terms in
\eqref{eq:affine-interaction-leading}; their coefficient is
\begin{align}
    k_a\binom{k_a-1}{p}(p+1)!.\label{dce}
\end{align}
Indeed, choosing the derivative acting on \(H_n\) and the \(p\)
derivatives acting on \(V_n\) gives
\(k_a\binom{k_a-1}{p}\). By the normalization in \eqref{dsym}, the
factor \(p!\) from the \(p\)-th Vandermonde derivative, together with
the choice of the distinguished element of \(A\), becomes the factor
\((p+1)!\) multiplying the symmetrized summand. Hence the coefficient
in \eqref{eq:affine-interaction-leading} is
(\ref{dce}).

All other terms have degree less than \(k_a+d\). Since \(k_a\) is
fixed, only finitely many parameter choices occur, and the same count
applies uniformly after any fixed number of further mixed derivatives.
Thus \eqref{eq:affine-interaction-remainder} follows. 

Since \(\mathscr J_{a,n}\) is linear in the derivative of its
argument, the difference
\(\mathscr J_{a,n}[H_n]-\mathscr J_{a,n}[\widetilde H_n]\) is obtained
from \eqref{eq:affine-interaction-leading} by replacing
\(\partial_{a_1}H_n\) with
\(\partial_{a_1}(H_n-\widetilde H_n)\). The assumed strict degree bound
and the preceding count then give
\eqref{eq:affine-interaction-replacement}.
\end{proof}
\endgroup

We now apply these preparations to the covariance calculation.

\begin{lemma}\label{le43}The covariance formula \eqref{covf} holds.
\end{lemma}
\begin{proof}
Fix $r_1,r_2\in[g]$ and set
\begin{equation}
    \kappa_i:=t_n-\lfloor\alpha_{r_i}n\rfloor,
    \qquad i=1,2.
    \label{dkr}
\end{equation}
By symmetry, after interchanging the two pairs $(r_i,k_i)$ if necessary,
we may assume $\kappa_1\geq\kappa_2$.  When the two levels coincide, the
same computation below is read with both operators acting on the common
level.  Equation \eqref{mmt} gives
\begin{align*}
&\frac{\operatorname{Cov}\left(
    p_{k_1}^{(\lfloor\alpha_{r_1}n\rfloor)},
    p_{k_2}^{(\lfloor\alpha_{r_2}n\rfloor)}
\right)}{n^{k_1+k_2}}\\
&=\left.\frac{1}{n^{k_1+k_2}}\left[
\mathcal{D}_{k_1,\kappa_1}\mathcal{D}_{k_2,\kappa_2}
\mathcal{S}_{\rho_{\kappa_2}}(|\mathbf{x}_{\kappa_2}|,\mathbf{u})
-\mathcal{D}_{k_1,\kappa_1}
\mathcal{S}_{\rho_{\kappa_1}}(|\mathbf{x}_{\kappa_1}|,\mathbf{u})
\mathcal{D}_{k_2,\kappa_2}
\mathcal{S}_{\rho_{\kappa_2}}(|\mathbf{x}_{\kappa_2}|,\mathbf{u})
\right]\right|_{\mathbf{u}=0}.
\end{align*}
\begingroup

Use the preceding notation with
\((\kappa_a,k_a)=(\kappa_i,k_i)\), \(i=1,2\), and write
\[
    s_i:=\lim_{n\to\infty}s_{i,n}=s_{r_i}.
\]
By Definition~\ref{def:global-F} and (\ref{dsnh}), and the change of variables
\(z_i=1+u_i/b_n\), we have
\[
    S_n(\mathbf u)=\widehat S_n(\mathbf z),
\] and
\eqref{eq:affine-euler-change} gives the exact change of the normalized
operators. The constant relating the Vandermonde determinants in the
\(\mathbf u\)- and \(\mathbf z\)-variables cancels from each normalized
operator, so the level dependence remains entirely in \(g_{i,n}\).

Consequently, the preceding operator identity becomes
\begin{equation}
\operatorname{Cov}\left(
    p_{k_1}^{(\lfloor\alpha_{r_1}n\rfloor)},
    p_{k_2}^{(\lfloor\alpha_{r_2}n\rfloor)}
\right)
=
\left.
\mathscr I_{1,n}[\mathcal F_{2,n}]
\right|_{\mathbf z=\mathbf1}.
\label{eq:covariance-as-interaction}
\end{equation}
Lemma~\ref{adfd} gives, uniformly in \(i\),
\[
    \deg_n(\partial_i\mathcal F_{2,n})\leq k_2,
    \qquad
    \deg_n\partial_i\left(
        \mathcal F_{2,n}-\mathcal F^{\mathrm{lead}}_{2,n}
    \right)<k_2.
\]
Applying Lemma~\ref{lem:affine-interaction} with
\[
    H_n=\mathcal F_{2,n},
    \qquad
    \widetilde H_n=\mathcal F^{\mathrm{lead}}_{2,n},
    \qquad d=k_2,
\]
we obtain
\begin{equation}
\operatorname{Cov}\left(
    p_{k_1}^{(\lfloor\alpha_{r_1}n\rfloor)},
    p_{k_2}^{(\lfloor\alpha_{r_2}n\rfloor)}
\right)
=
\left.
\mathscr J_{1,n}[\mathcal F^{\mathrm{lead}}_{2,n}]
\right|_{\mathbf z=\mathbf1}
+o(n^{k_1+k_2}).
\label{eq:covariance-leading-interaction}
\end{equation}
Set
\begin{align}
\mathcal C_n
&:=\mathscr J_{1,n}[\mathcal F^{\mathrm{lead}}_{2,n}]
\notag\\
&=k_1\sum_{p=0}^{k_1-1}
\binom{k_1-1}{p}(p+1)!
\sum_{\substack{A\subset[n]\\ |A|=p+1}}
\operatorname{Sym}_{a_1,\ldots,a_{p+1}}
\frac{
    g_{1,n}(z_{a_1})^{k_1}
    (\partial_{a_1}L_n)^{k_1-1-p}
    \partial_{a_1}\mathcal F^{\mathrm{lead}}_{2,n}
}{
    \prod_{j=2}^{p+1}(z_{a_1}-z_{a_j})
}.
\label{cc1}
\end{align}

Put
\[
    F(z):=H_{\mathbf m_0}'(z)
\]
and
\begin{equation}
    G(z,w)
    :=
    \frac{\partial^2}{\partial z\,\partial w}
    \log\left(
        1-
        \frac{(z-1)(w-1)}{z-w}
        \bigl[zF(z)-wF(w)\bigr]
    \right).
    \label{eq:covariance-mixed-kernel}
\end{equation}
Lemmas~\ref{la1} and~\ref{la3} give, locally uniformly near
$(1,1)$,
\[
    \frac1n\partial_iL_n\longrightarrow F(z_i),
    \qquad
    \partial_i\partial_jL_n\longrightarrow G(z_i,z_j)
    \quad(i\ne j).
\]

Substitute the expansion
\eqref{eq:correct-inner-leading} of
\(\mathcal F^{\mathrm{lead}}_{2,n}\) into \eqref{cc1}.
The resulting double sum is indexed by
\[
    A=\{a_1,\ldots,a_{p+1}\},
    \qquad
    B=\{b_1,\ldots,b_{q+1}\}.
\]
We group its summands according to the size of \(A\cap B\).
If $A\cap B=\varnothing$, then $\partial_{a_1}$ acts on the inner
summand only through $L_n$.  For $0\le q\le k_2-1$,
\begin{equation}
    \partial_{a_1}(\partial_{b_1}L_n)^{k_2-q}
    =
    (k_2-q)(\partial_{b_1}L_n)^{k_2-q-1}
    \partial_{a_1}\partial_{b_1}L_n.
    \label{eq:correct-chain-rule}
\end{equation}
For $q=k_2$, the power $(\partial_{b_1}L_n)^0$ is identically one,
so its derivative is zero.  Thus the term $q=k_2$ has no
disjoint-index contribution. 

For fixed \(p\) and \(q\), each disjoint-index summand has order
\(n^{k_1+k_2-p-q-2}\), whereas the number of disjoint pairs
\((A,B)\) is
\[
    \binom{n}{p+1}\binom{n-p-1}{q+1}
    \sim
    \frac{n^{p+q+2}}{(p+1)!(q+1)!}.
\]
Lemma~\ref{la2} converts the two symmetrized quotients into the
\(p\)-th and \(q\)-th derivatives at \(1\). By Lemma~\ref{la2}, when the variables in the two symmetrized
quotients tend to \(1\), the quotients contribute the factors
\(1/(p+1)!\) and \(1/(q+1)!\), respectively. These cancel the
prefactors \((p+1)!\) and \((q+1)!\) already present in the two
operator expansions. The factorials remaining in the formula below
come from the number of choices of the disjoint sets \(A\) and \(B\). Using
\[
    (k_2-q)\binom{k_2}{q}
    =
    k_2\binom{k_2-1}{q},
\]
together with the locally uniform limits above, we conclude that the
disjoint-index contribution, divided by \(n^{k_1+k_2}\), converges to
\begin{align}
&\sum_{p=0}^{k_1-1}\sum_{q=0}^{k_2-1}
\frac{k_1k_2}{(p+1)!(q+1)!}
\binom{k_1-1}{p}\binom{k_2-1}{q}
\notag\\[-2mm]
&\qquad{}\times
\left.
\partial_z^p\partial_w^q
\left\{
    g_1(z)^{k_1}F(z)^{k_1-1-p}
    g_2(w)^{k_2}F(w)^{k_2-1-q}
    G(z,w)
\right\}
\right|_{z=w=1},
\label{eq:disjoint-finite-sum}
\end{align}
where \(g_i(z):=z-1+s_i\).  Define
\[
    \Phi_i(z)
    :=g_i(z)\left(F(z)+\frac1{z-1}\right).
\]
For every function $\Psi$ analytic near $1$, Cauchy's formula and the
binomial theorem give
\begin{align}
&\frac1{2\pi\mathbf i}
\oint_{|z-1|=\epsilon}\Phi_i(z)^{k_i}\Psi(z)\,dz
\notag\\
&\quad=
\sum_{p=0}^{k_i-1}
\frac{k_i}{(p+1)!}\binom{k_i-1}{p}
\left.
\partial_z^p
\bigl[g_i(z)^{k_i}F(z)^{k_i-1-p}\Psi(z)\bigr]
\right|_{z=1}.
\label{eq:one-variable-residue-identity}
\end{align}
Applying \eqref{eq:one-variable-residue-identity} in both variables
shows that \eqref{eq:disjoint-finite-sum} equals
\begin{equation}
\frac1{(2\pi\mathbf i)^2}
\oint_{|z-1|=\epsilon}\oint_{|w-1|=2\epsilon}
    \Phi_1(z)^{k_1}\Phi_2(w)^{k_2}G(z,w)\,dz\,dw.
\label{eq:disjoint-contour}
\end{equation}

It remains to consider $|A\cap B|=1$.  For $0\le q\le k_2$, write
\[
    g_2(z)^{k_2}F(z)^{k_2-q}
    =\sum_{j\ge0}c_j^{(q)}(z-1)^j
\]
and set
\[
    \mathcal B_q(z)
    :=
    \frac{
        g_2(z)^{k_2}F(z)^{k_2-q}
        -\sum_{j=0}^{q-1}c_j^{(q)}(z-1)^j
    }{(z-1)^q},
\]
where the sum in the numerator is empty for $q=0$; by construction,
$\mathcal B_q$ is analytic near $1$.

By symmetry, assume
\(A\cap B=\{a_1\}=\{b_1\}\) and set
\(z=z_{a_1}=z_{b_1}\).
Keeping \(z\) free while letting the remaining \(B\)-coordinates tend
to \(1\), Lemma~\ref{lem:partial-confluent-divided-difference}
reduces the inner divided difference to \(\mathcal B_q(z)\); hence the
outer derivative produces \(\mathcal B_q'(z)\).
Applying Lemma~\ref{la2} to the remaining \(A\)-coordinates and
including the symmetrization and index-counting factors, we obtain,
after division by \(n^{k_1+k_2}\),
\begin{align}
&\sum_{p=0}^{k_1-1}\sum_{q=0}^{k_2}
\frac{k_1}{(p+1)!}
\binom{k_1-1}{p}\binom{k_2}{q}
\notag\\[-1mm]
&\qquad\times
\left.
\partial_z^p
\left[
    g_1(z)^{k_1}F(z)^{k_1-1-p}\mathcal B_q'(z)
\right]
\right|_{z=1}.
\label{eq:one-intersection-finite-sum}
\end{align}
Notice that $\mathcal B_{k_2}(z)=1$, and hence the endpoint $q=k_2$ in
\eqref{eq:one-intersection-finite-sum} is zero.  

Applying \eqref{eq:one-variable-residue-identity} with
\[
    \Psi(z)
    =
    \sum_{q=0}^{k_2}
    \binom{k_2}{q}\mathcal B_q'(z)
\]
shows that \eqref{eq:one-intersection-finite-sum} equals
\begin{equation}
    \frac1{2\pi\mathbf i}
    \oint_{|z-1|=\epsilon}
    \Phi_1(z)^{k_1}
    \frac{d}{dz}
    \left[
        \sum_{q=0}^{k_2}
        \binom{k_2}{q}\mathcal B_q(z)
    \right]dz.
    \label{eq:one-intersection-positive-part}
\end{equation}
Indeed,
\[
    \Phi_2(z)^{k_2}
    =
    \sum_{q=0}^{k_2}
    \binom{k_2}{q}
    \frac{
        g_2(z)^{k_2}F(z)^{k_2-q}
    }{
        (z-1)^q
    },
\]
so the expression in brackets in
\eqref{eq:one-intersection-positive-part} is obtained by retaining the
nonnegative powers in the Laurent expansion of
\(\Phi_2(z)^{k_2}\) at \(z=1\). Consequently, for
\(0<\delta<\epsilon\),
\[
    \sum_{q=0}^{k_2}
    \binom{k_2}{q}\mathcal B_q(z)
    =
    \Phi_2(z)^{k_2}
    -
    \frac1{2\pi\mathbf i}
    \oint_{|w-1|=\delta}
    \frac{\Phi_2(w)^{k_2}}{z-w}\,dw.
\]
Substitute this identity into
\eqref{eq:one-intersection-positive-part} and move the \(w\)-contour
from \(|w-1|=\delta\) to \(|w-1|=2\epsilon\). The only crossed pole is
\(w=z\), whose residue is precisely the term obtained by
differentiating \(\Phi_2(z)^{k_2}\). Hence
\begin{equation}
\eqref{eq:one-intersection-finite-sum}
=
\frac1{(2\pi\mathbf i)^2}
\oint_{|z-1|=\epsilon}
\oint_{|w-1|=2\epsilon}
\frac{\Phi_1(z)^{k_1}\Phi_2(w)^{k_2}}
     {(z-w)^2}\,dw\,dz.
\label{eq:one-intersection-contour}
\end{equation}
The corresponding residue identity for the standard factor
\(1+\zeta\) appears in \cite[Lemma~6.2]{bg16}; the calculation above
establishes the level-dependent form required here.

Finally, if $|A\cap B|\ge2$, the summation has at least one fewer free
index while the derivative factors obey the same degree bounds.
Hence its total contribution is $o(n^{k_1+k_2})$.  Combining
\eqref{eq:disjoint-contour} and
\eqref{eq:one-intersection-contour} yields
\[
\frac1{(2\pi\mathbf i)^2}
\oint_{|z-1|=\epsilon}\oint_{|w-1|=2\epsilon}
\Phi_1(z)^{k_1}\Phi_2(w)^{k_2}
\left(G(z,w)+\frac1{(z-w)^2}\right)
\,dz\,dw,
\]
which is exactly \eqref{covf}.  
\endgroup
\end{proof}

\begingroup

For \(a<b\) in the chosen operator order, write
\begin{equation}
    \mathcal G_{a,b,n}(\mathbf z)
    :=
    \mathscr J_{a,n}[\mathcal F_{b,n}](\mathbf z).
    \label{eq:affine-pair-function}
\end{equation}

\begin{lemma}[Two-operator bounds]\label{lem:affine-pair-bounds}
Suppose \(s_{a,n}\to s_a\in(0,1)\) and
\(s_{b,n}\to s_b\in(0,1)\). In the sense of
Definition~\ref{def:derivative-stable-degree},
\[
    \deg_n\mathcal G_{a,b,n}\leq k_a+k_b,
    \qquad
    \deg_n(\partial_i\mathcal G_{a,b,n})<k_a+k_b.
\]
\end{lemma}

\begin{proof}
The first estimate follows from Lemmas~\ref{adfd} and
\ref{lem:affine-interaction}, applied with
\(H_n=\mathcal F_{b,n}\) and \(d=k_b\).

After applying a fixed derivative \(\partial_i\), terms whose index set
contains \(i\) have one fewer free summation index. In all other terms,
\(\partial_i\) either creates a mixed derivative of \(L_n\) or acts on
\(\partial_{a_1}\mathcal F_{b,n}\); the corresponding strict bounds in
Lemma~\ref{adfd} lower the degree by one. Hence
\[
    \deg_n(\partial_i\mathcal G_{a,b,n})<k_a+k_b.
\]
The same counting remains valid after applying any fixed number of
additional coordinate derivatives, which proves the stated degree
bounds.
\end{proof}

For a finite set $I$, let $\operatorname{Pair}(I)$ denote the set of
complete pairings of $I$.  We set $\operatorname{Pair}(I)=\varnothing$
when $|I|$ is odd, while $\operatorname{Pair}(\varnothing)$ consists of
the unique empty pairing, whose associated product is one.

\begin{lemma}[Function-valued multilevel Wick expansion]\label{le47}
Let \(m\geq1\), let \(\kappa_1\geq\cdots\geq\kappa_m\), and let
\(k_1,\ldots,k_m\) be positive integers.  Assume
\[
    s_{a,n}
    =
    \frac{|\mathbf x^{(n)}_{\kappa_a}|}{|\mathbf x^{(n)}|}
    \longrightarrow s_a\in(0,1),
    \qquad a\in[m].
\]
Let
\begin{equation}
    P_{a,n}
    :=
    \sum_{i=1}^{n}
    \left(\lambda_i^{(\kappa_a)}+n-i\right)^{k_a},
    \qquad a\in[m].
    \label{eq:abstract-level-power-sum}
\end{equation}
For each \(a\in[m]\), let \(\mathscr D_{a,n}\) be the operator
associated with \((\kappa_a,k_a)\) in
\eqref{dsag}.
For an ordered subset \(I=\{i_1<\cdots<i_\ell\}\subseteq[m]\), define, in a neighborhood of \(\mathbf z=\mathbf1^n\),
\begin{equation}
    \mathcal M_{I,n}(\mathbf z)
    :=
    \frac{
        \mathscr D_{i_1,n}\cdots\mathscr D_{i_\ell,n}
        \bigl(V_n(\mathbf z)\widehat S_n(\mathbf z)\bigr)
    }{
        V_n(\mathbf z)\widehat S_n(\mathbf z)
    },
    \qquad
    \mathcal M_{\varnothing,n}:=1.
    \label{eq:function-valued-moment}
\end{equation}
Define \(\mathcal A_{I,n}\) in the same neighborhood by
\begin{equation}
    \mathcal A_{I,n}
    :=
    \sum_{J\subseteq I}
    (-1)^{|J|}
    \left(\prod_{j\in J}\mathcal F_{j,n}\right)
    \mathcal M_{I\setminus J,n}.
    \label{eq:connected-germ-definition}
\end{equation}
Equivalently, before any specialization of the variables,
\begin{equation}
    \mathcal M_{I,n}
    =
    \sum_{J\subseteq I}
    \left(\prod_{j\in J}\mathcal F_{j,n}\right)
    \mathcal A_{I\setminus J,n}.
    \label{eq:function-valued-singleton-expansion}
\end{equation}
Set \(\mathcal R_{\varnothing,n}:=0\), and, for nonempty \(I\), set
\begin{equation}
    \mathcal R_{I,n}
    :=
    \mathcal A_{I,n}
    -
    \sum_{P\in\operatorname{Pair}(I)}
    \prod_{(a,b)\in P}\mathcal G_{a,b,n},
    \label{eq:function-valued-wick}
\end{equation}
where each pair is oriented by \(a<b\).  Then
\begin{equation}
    \deg_n\mathcal R_{I,n}
    <
    \sum_{a\in I}k_a
    \qquad(I\ne\varnothing).
    \label{eq:wick-remainder-all-jets}
\end{equation}
Consequently, the exact function-valued expansion is
\begin{align}
    \mathcal M_{I,n}
    ={}&
    \sum_{J\subseteq I}
    \left(\prod_{j\in J}\mathcal F_{j,n}\right)
    \left(
        \sum_{P\in\operatorname{Pair}(I\setminus J)}
        \prod_{(a,b)\in P}\mathcal G_{a,b,n}
        +\mathcal R_{I\setminus J,n}
    \right).
    \label{eq:raw-function-valued-wick}
\end{align}
Finally,
\begin{align}
&\mathbb E
\prod_{a=1}^{m}
\left(
    P_{a,n}-\mathbb E P_{a,n}
\right)
\notag\\
&\quad=
\left.
\left(
    \sum_{P\in\operatorname{Pair}([m])}
    \prod_{(a,b)\in P}\mathcal G_{a,b,n}
    +\mathcal R_{[m],n}
\right)
\right|_{\mathbf z=1^n}.
\label{eq:centered-wick-expansion}
\end{align}
\end{lemma}

\begin{proof}
Equation~\eqref{eq:function-valued-singleton-expansion} follows from
\eqref{eq:connected-germ-definition} by inclusion--exclusion over the
subsets of \(I\).
We prove \eqref{eq:wick-remainder-all-jets} simultaneously for all
ordered sets \(I\), by induction on \(|I|\). Throughout the induction,
all identities hold as analytic function identities in a neighborhood
of \(\mathbf z=\mathbf1^n\); we set \(\mathbf z=\mathbf1^n\) only
after the induction is complete. Building on the induction scheme of
\cite[Proposition~5.10]{bg16}, we prove below the function-valued
expansion for the ordered, level-dependent operators
\(\mathscr D_{a,n}\), with remainder bounds preserved under every
fixed number of coordinate derivatives.

For \(|I|=1\), the assertion follows from the definition of
\(\mathcal F_{a,n}\), and \(\mathcal R_{I,n}=0\).  

For \(I=\{a,b\}\) with \(a<b\), by (\ref{dti}),
\[
    \mathcal M_{I,n}
    =\mathscr T_{a,n}[\mathcal F_{b,n}],
    \qquad
    \mathcal A_{I,n}
    =\mathscr I_{a,n}[\mathcal F_{b,n}].
\]
Since
\(\mathcal G_{a,b,n}
=\mathscr J_{a,n}[\mathcal F_{b,n}]\), we have
\[
    \mathcal R_{I,n}
    =
    \mathscr I_{a,n}[\mathcal F_{b,n}]
    -
    \mathscr J_{a,n}[\mathcal F_{b,n}].
\]
Lemma~\ref{lem:affine-interaction}, applied with
\(H_n=\mathcal F_{b,n}\) and \(d=k_b\), therefore gives
\[
    \deg_n\mathcal R_{I,n}<k_a+k_b.
\]

For the induction step, let \(a\) be the first occurrence in \(I\),
put \(I'=I\setminus\{a\}\), and use the exact identity
\[
    \mathcal M_{I,n}
    =
    \frac{1}{V_n\widehat S_n}
    \mathscr D_{a,n}
    \left(V_n\widehat S_n\,\mathcal M_{I',n}\right).
\]
Substitute the function-valued expansion
\eqref{eq:raw-function-valued-wick} for \(\mathcal M_{I',n}\) into
this identity.  For a term indexed by \(J\subseteq I'\), write
\[
    H_{K,n}
    :=
    \sum_{P\in\operatorname{Pair}(K)}
    \prod_{(c,d)\in P}\mathcal G_{c,d,n}
    +\mathcal R_{K,n},
    \qquad
    K:=I'\setminus J,
\]
and put \(D_K:=\sum_{i\in K}k_i\).  The induction hypothesis and
Lemma~\ref{lem:affine-pair-bounds} give
\[
    \deg_n H_{K,n}\leq D_K,
    \qquad
    \deg_n(\partial_{z_i}H_{K,n})<D_K.
\]
Indeed, each term obtained by differentiating a pairing product
contains one differentiated \(\mathcal G\)-factor, whose degree is
strictly smaller by Lemma~\ref{lem:affine-pair-bounds}; the induction
hypothesis gives the same strict degree bound for every fixed-order
derivative of \(\mathcal R_{K,n}\).

If no derivative from \(\mathscr D_{a,n}\) falls on the inserted factor
\(\prod_{j\in J}\mathcal F_{j,n}H_{K,n}\), the sum of all such terms is
exactly
\[
    \mathcal F_{a,n}
    \left(\prod_{j\in J}\mathcal F_{j,n}\right)H_{K,n};
\]
these are precisely the summands in which \(a\) is a singleton.  It
remains to estimate the terms in which at least one inserted factor is
differentiated.

Expand the outer operator by the product rule. In a fixed term, let
\(q_0,r,f_j\), and \(h\) count the derivatives falling on
\(g_{a,n},V_n,\mathcal F_{j,n}\), and \(H_{K,n}\), respectively, and
write the contribution from \(\widehat S_n=e^{L_n}\) as
\[
    \prod_{\ell\geq1}
    \bigl(\partial^\ell L_n\bigr)^{d_\ell}.
\]
Since the operator has order \(k_a\),
\begin{equation}
    q_0+r+\sum_{j\in J}f_j+h
    +\sum_{\ell\geq1}\ell d_\ell
    =k_a.
    \label{eq:function-wick-derivative-budget}
\end{equation}

Set
\[
    A:=\{j\in J:f_j>0\},
    \qquad
    B:=J\setminus A,
\]
and write the term as
\[
    \left(\prod_{j\in B}\mathcal F_{j,n}\right)C_n.
\]
With
\[
    L:=I\setminus B=\{a\}\cup K\cup A,
    \qquad
    D_L:=k_a+D_K+\sum_{j\in A}k_j,
\]
Lemma~\ref{adfd}, the bounds for \(H_{K,n}\), and the \(r+1\) free
indices give
\begin{align}
    \deg_n C_n
    &\leq
    D_K+\sum_{j\in A}k_j+\sum_{\ell\geq1}d_\ell+r+1
    \notag\\
    &=
    D_L+1-q_0
    -\sum_{\ell\geq1}(\ell-1)d_\ell
    -\sum_{j\in A}f_j-h.
    \label{eq:function-wick-degree-bound}
\end{align}
The first inequality is strict if \(h>0\). Since at least one inserted
factor is differentiated, equality with \(D_L\) is possible only when
\[
    A=\{j\},\qquad f_j=1,\qquad h=q_0=0,\qquad
    d_\ell=0\quad(\ell\geq2).
\]
These are precisely the leading interaction terms isolated in
Lemma~\ref{lem:affine-interaction}, with
\(H_n=\mathcal F_{j,n}\). Since \(H_{K,n}\) is not differentiated, it
factors out; hence, by \eqref{eq:affine-pair-function}, their sum is
\[
    \mathcal G_{a,j,n}H_{K,n}.
\] Its pairing part creates the pair
\((a,j)\), whereas
\(\mathcal G_{a,j,n}\mathcal R_{K,n}\) and all remaining terms have
degree less than \(D_L\).

Grouping the terms by the set \(B\) of undifferentiated singleton
factors, inclusion--exclusion identifies the contribution with
\(B=\varnothing\) as \(\mathcal A_{I,n}\). The preceding degree count
shows that its only terms of degree \(D_I\) are the pairing products.
Hence, by \eqref{eq:function-valued-wick},
\[
    \deg_n\mathcal R_{I,n}<D_I,
\]
which completes the induction.

Finally, the multilevel moment identity \eqref{mmt} gives
\[
    \mathcal M_{I,n}(\mathbf1^n)
    =
    \mathbb E\prod_{a\in I}P_{a,n},
    \qquad
    \mathcal F_{a,n}(\mathbf1^n)
    =
    \mathbb E P_{a,n}.
\]
Hence, evaluating \eqref{eq:connected-germ-definition} at
\(\mathbf z=\mathbf1^n\) with \(I=[m]\) gives
\[
    \mathcal A_{[m],n}(\mathbf1^n)
    =
    \mathbb E\prod_{a=1}^{m}
    \left(P_{a,n}-\mathbb E P_{a,n}\right).
\]
Evaluating \eqref{eq:function-valued-wick} at
\(\mathbf z=\mathbf1^n\) now yields
\eqref{eq:centered-wick-expansion}.
\end{proof}
\endgroup

\begin{proof}[Proof of Theorem~\ref{t31}]
The covariance formula follows from Lemma~\ref{le43}. For $r\in[g]$, put
\[
    \kappa_r:=t_n-\lfloor\alpha_r n\rfloor.
\]
We now verify the higher moments. Fix a finite list of pairs
\[
    (r_1,k_1),\ldots,(r_m,k_m),
    \qquad r_a\in[g],\quad k_a\geq1,
\]
and write
\[
    X_a^{(n)}
    :=
    n^{-k_a}\left(
        p_{k_a}^{(\lfloor\alpha_{r_a}n\rfloor)}
        -\mathbb E p_{k_a}^{(\lfloor\alpha_{r_a}n\rfloor)}
    \right).
\]
\begingroup

Relabel the observables, if necessary, so that their levels satisfy
$\kappa_{r_1}\geq\cdots\geq\kappa_{r_m}$; equal levels may be ordered
arbitrarily.  Apply Lemma~\ref{le47} with its occurrence-level
$\kappa_a$ equal to $\kappa_{r_a}$.  Then
$P_{a,n}=p_{k_a}^{(\lfloor\alpha_{r_a}n\rfloor)}$, and the
corresponding ratios satisfy $s_{a,n}\to s_{r_a}\in(0,1)$.  By
\eqref{eq:centered-wick-expansion} and
\eqref{eq:wick-remainder-all-jets},
\begin{align*}
&n^{-\sum_{a=1}^{m}k_a}
\mathbb E\prod_{a=1}^{m}
\left(
    p_{k_a}^{(\lfloor\alpha_{r_a}n\rfloor)}
    -\mathbb E p_{k_a}^{(\lfloor\alpha_{r_a}n\rfloor)}
\right)\notag\\
&\quad=
n^{-\sum_{a=1}^{m}k_a}
\left.
\sum_{P\in\operatorname{Pair}([m])}
\prod_{(a,b)\in P}\mathcal G_{a,b,n}
\right|_{\mathbf z=1^n}
+o(1).
\end{align*}
For every pair $(a,b)$, Lemmas~\ref{adfd} and
\ref{lem:affine-interaction}, specifically
\eqref{eq:affine-interaction-replacement}, give
\[
    \mathcal G_{a,b,n}(\mathbf1^n)
    -
    \mathscr J_{a,n}
    [\mathcal F^{\mathrm{lead}}_{b,n}](\mathbf1^n)
    =o(n^{k_a+k_b}).
\]
Therefore the calculation in Lemma~\ref{le43} gives
\[
    n^{-k_a-k_b}
    \mathcal G_{a,b,n}(1^n)
    \longrightarrow
    \Sigma_{(r_a,k_a),(r_b,k_b)}.
\]
It follows that
\[
    \lim_{n\to\infty}
    \mathbb E\left[\prod_{a=1}^{m}X_a^{(n)}\right]
    =
    \begin{cases}
        0, & m\text{ odd},\\[1mm]
        \displaystyle
        \sum_{P\in\operatorname{Pair}([m])}
        \prod_{(a,b)\in P}
        \Sigma_{(r_a,k_a),(r_b,k_b)},
        & m\text{ even}.
    \end{cases}
\]
These are exactly the Wick moments of the centered Gaussian vector with
covariance \eqref{covf}. Hence every finite subcollection of
\eqref{crv} converges in moments to that Gaussian vector.
\endgroup
\end{proof}

\begin{assumption}\label{ap32}Let $l\geq2$ be a fixed integer. Assume there exists
\begin{align*}
    0=a_1<b_1<a_2<b_2<\ldots<a_{l}<b_l
\end{align*}
such that the probability measure $\mathbf m_0$ associated with the
regular bottom-boundary sequence satisfies
\[
    \frac{d\mathbf m_0}{dx}
    =
    \mathbf 1_{\bigcup_{i=1}^l(a_i,b_i)}(x).
\]
\end{assumption}

\begin{lemma}
\label{le49}
Suppose Assumption~\ref{ap32} holds. For every $\chi\in\mathbb R$, the
equation $U_y(z)=\chi$, with $U_y$ defined by \eqref{duyz}, has at most
one nonreal complex-conjugate pair of roots.
\end{lemma}

\begin{proof}
Under Assumption~\ref{ap32}, equation \eqref{hmd} gives
\[
    H_{\mathbf m_0}'(z)
    =
    -\frac{1}{z-1}+\frac{\zeta}{z},
\]
where the spectral parameter $\zeta$ is determined by
\begin{equation}
    z
    =
    \prod_{i=1}^{l}
    \frac{\zeta-a_i}{\zeta-b_i}.
    \label{ee1}
\end{equation}
Consequently,
\begin{equation}
    U_y(z)
    =
    \frac{(z-1+s)\zeta}{z}.
    \label{ee2}
\end{equation}
Every nonreal root $z$ of $U_y(z)=\chi$ determines a nonreal value of
$\zeta$: if $\zeta$ were real, then \eqref{ee1} would make $z$ real.
Eliminating $z$ from \eqref{ee1}--\eqref{ee2}, and using $a_1=0$, shows
that $\zeta$ is a zero of
\begin{equation}
    P_{\chi,s}(\zeta)
    :=
    (\zeta-\chi)
    \prod_{i=2}^{l}(\zeta-a_i)
    -(1-s)
    \prod_{i=1}^{l}(\zeta-b_i).
    \label{eq:spectral-polynomial}
\end{equation}
The polynomial $P_{\chi,s}$ has degree $l$ and leading coefficient
$s>0$. For $j=2,\ldots,l$,
\[
    P_{\chi,s}(a_j)
    =
    -(1-s)\prod_{i=1}^{l}(a_j-b_i),
\]
which is nonzero and has sign $(-1)^{l-j}$. Thus consecutive values
$P_{\chi,s}(a_j)$ and $P_{\chi,s}(a_{j+1})$ have opposite signs.
When $l\geq3$, the intermediate value theorem gives one real zero in
each interval
\[
    (a_j,a_{j+1}),
    \qquad j=2,\ldots,l-1.
\]
Hence $P_{\chi,s}$ has at least $l-2$ distinct real zeros. For $l=2$
the same conclusion is vacuous, while $P_{\chi,s}$ has degree two.
In every case, $P_{\chi,s}$ has at most two nonreal zeros, which must
form a conjugate pair because its coefficients are real. Since a fixed spectral parameter determines $z$ uniquely through
\eqref{ee1}, distinct nonreal roots $z$ give distinct spectral parameters, the equation $U_y(z)=\chi$ also has at most one nonreal
complex-conjugate pair.
\end{proof}

\begin{lemma}
\label{lem:spectral-liquid-equivalence}
Suppose that Assumption~\ref{ap32} and
Assumption~\ref{ass:tail-profile}\ref{ass:tail-profile-invertibility}
hold. Let $y\in(0,\alpha)$ and set $s=s(y)$. Define
\begin{equation}
    V_s(\zeta)
    :=
    \zeta\left[
        1-(1-s)e^{-\mathrm{St}_{\bm_0}(\zeta)}
    \right].
    \label{dvy}
\end{equation}
Then, for every $\chi\in\mathbb R$, the equation
$V_s(\zeta)=\chi$ has either zero or one root in $\mathbb H$.
Whenever such a root exists, it is simple. Moreover,
$(\chi,s)\in\widetilde{\mathcal L}$ if and only if
$V_s(\zeta)=\chi$ has a nonreal complex-conjugate pair of roots.
\end{lemma}

\begin{proof}
Set
\[
    E:=\bigcup_{i=1}^{l}(a_i,b_i).
\]
Under Assumption~\ref{ap32}, for $\zeta\in\mathbb H$ we have
\begin{equation}
    \mathrm{St}_{\bm_0}(\zeta)
    =
    \sum_{i=1}^{l}
    \log\frac{\zeta-a_i}{\zeta-b_i},
    \qquad
    z(\zeta)
    :=
    e^{\mathrm{St}_{\bm_0}(\zeta)}
    =
    \prod_{i=1}^{l}
    \frac{\zeta-a_i}{\zeta-b_i},
    \label{eq:spectral-z-physical-proof}
\end{equation}
where the logarithmic branches are fixed by continuation from infinity.
Equations~\eqref{ee1}--\eqref{ee2} therefore give
\begin{equation}
    V_s(\zeta)=U_y\bigl(z(\zeta)\bigr).
    \label{eq:Vs-Uy-physical-proof}
\end{equation}

We first identify the solution selected by the Stieltjes transform.
Write
\[
    \mathscr R(\zeta)
    :=
    \zeta e^{-\mathrm{St}_{\bm_0}(\zeta)}
    =
    \frac{\prod_{i=1}^{l}(\zeta-b_i)}
         {\prod_{i=2}^{l}(\zeta-a_i)}.
\]
Its partial-fraction expansion has the form
\begin{equation}
    \mathscr R(\zeta)
    =
    \zeta+c+
    \sum_{j=2}^{l}\frac{A_j}{\zeta-a_j},
    \qquad
    A_j
    =
    \frac{\prod_{i=1}^{l}(a_j-b_i)}
         {\prod_{\substack{2\leq i\leq l\\i\neq j}}(a_j-a_i)}
    <0.
    \label{eq:R-partial-fraction-physical-proof}
\end{equation}
Indeed, the numerator in the expression for $A_j$ has
$l-j+1$ negative factors, whereas the denominator has $l-j$
negative factors. Hence
\begin{equation}
    V_s(\zeta)
    =
    s\zeta-(1-s)c+
    \sum_{j=2}^{l}\frac{B_j}{\zeta-a_j},
    \qquad
    B_j:=-(1-s)A_j>0.
    \label{eq:Vs-partial-fraction-physical-proof}
\end{equation}

Let $w\in\mathbb H$. Clearing the poles in $V_s(\zeta)=w$ gives
\begin{equation}
    P_{w,s}(\zeta)
    :=
    (\zeta-w)\prod_{i=2}^{l}(\zeta-a_i)
    -(1-s)\prod_{i=1}^{l}(\zeta-b_i)
    =0.
    \label{eq:spectral-polynomial-complex-w}
\end{equation}
This polynomial has degree $l$ and leading coefficient $s>0$. It has
no real zero: away from the poles, $V_s$ is real on the real axis,
whereas $w\notin\mathbb R$, and at every pole $a_j$ the polynomial in
\eqref{eq:spectral-polynomial-complex-w} is nonzero.

Let $N(w)$ be the number of zeros of
\eqref{eq:spectral-polynomial-complex-w} in $\mathbb H$, counted with
multiplicity. Since no zero can cross the real axis and the leading
coefficient is independent of $w$, $N(w)$ is locally constant on the
connected set $\mathbb H$. To determine it, take $w=\mathbf iT$ with
$T\to\infty$. For the root tending to infinity, put $\zeta=w\xi$;
then $V_s(w\xi)/w=s\xi+O(T^{-1})$ near $\xi=s^{-1}$. For a root near
$a_j$, put $\zeta=a_j+\eta/w$; then
\[
    \frac{V_s(a_j+\eta/w)}{w}
    =
    \frac{B_j}{\eta}+O(T^{-1})
\]
uniformly near $\eta=B_j$. Rouch\'e's theorem therefore gives one root
\[
    \zeta_\infty(w)=\frac{w}{s}+O(1)
\]
and, for each $j=2,\ldots,l$, one root
\[
    \zeta_j(w)
    =
    a_j+\frac{B_j}{w}+O(|w|^{-2}).
\]
Thus, for $w=\mathbf iT$ and $T$ sufficiently large,
$\zeta_\infty(w)\in\mathbb H$, while
$\zeta_j(w)\in\{\Im\zeta<0\}$ for $j=2,\ldots,l$. Consequently,
\[
    N(w)=1,
    \qquad w\in\mathbb H.
\]
Denote the unique upper-half-plane root by $\zeta_s(w)$. Since it is
counted with multiplicity one, it is simple, and the implicit-function
theorem shows that $w\mapsto\zeta_s(w)$ is holomorphic on $\mathbb H$.

For $\zeta=x+\mathbf i\eta\in\mathbb H$,
\begin{equation}
    0
    <
    -\Im\mathrm{St}_{\bm_0}(\zeta)
    =
    \int_E
    \frac{\eta}{(x-u)^2+\eta^2}\,du
    <
    \int_{\mathbb R}
    \frac{\eta}{(x-u)^2+\eta^2}\,du
    =
    \pi.
    \label{eq:St-m0-imaginary-range}
\end{equation}
Hence
\[
    z_s(w)
    :=
    e^{\mathrm{St}_{\bm_0}(\zeta_s(w))}
\]
lies in the lower half-plane. By
\eqref{eq:Vs-Uy-physical-proof},
\[
    U_y\bigl(z_s(w)\bigr)=w.
\]
Moreover, as $w\to\infty$ in $\mathbb H$,
\[
    \zeta_s(w)=\frac{w}{s}+O(1),
    \qquad
    z_s(w)=1+\frac{s}{w}+O(w^{-2}).
\]
Thus $z_s(w)$ agrees near infinity with the solution branch
$\mathbf z_1(w,y)$ from
Proposition~\ref{prop:stieltjes-root}, and hence gives its holomorphic
continuation throughout $\mathbb H$.

Since $(z_s(w)-1+s)/s$ also lies in the lower half-plane, the principal
logarithm defines a holomorphic function
\[
    G_s(w)
    :=
    \operatorname{Log}\left(
        \frac{z_s(w)-1+s}{s}
    \right),
    \qquad w\in\mathbb H,
\]
with
\[
    -\pi<\Im G_s(w)<0.
\]
Near infinity, Proposition~\ref{prop:stieltjes-root} gives
$G_s(w)=\mathrm{St}_{\bm_y}(w)$. The identity theorem therefore yields
\begin{equation}
    \mathrm{St}_{\bm_y}(w)=G_s(w),
    \qquad w\in\mathbb H.
    \label{eq:physical-spectral-branch-identity}
\end{equation}

We now return to real $\chi$. The preceding upper-half-plane root count
implies that $V_s(\zeta)=\chi$ has at most one root in $\mathbb H$.
Indeed, two distinct upper-half-plane roots, or one root of multiplicity
at least two, would persist under a sufficiently small perturbation of
$\chi$ into $\mathbb H$ by Rouch\'e's theorem, contradicting $N(w)=1$.
Thus every upper-half-plane root for real $\chi$ is simple. Since
$V_s$ has real coefficients, such a root exists if and only if there is
a nonreal complex-conjugate pair.

Suppose first that $(\chi,s)\in\widetilde{\mathcal L}$, and let
$y=y(s)$. By Lemma~\ref{lem:liquid-nonreal-root}, the physical branch
$z_0:=\mathbf z_1(\chi,y)$ lies in the lower half-plane. All particle
coordinates are nonnegative, so the strict positivity of the density at
a regular liquid point implies $\chi>0$. The spectral parameter
corresponding to $z_0$ in \eqref{ee1}--\eqref{ee2} is
\[
    \zeta_0
    =
    \frac{\chi z_0}{z_0-1+s}.
\]
Since $\Im z_0<0$ and $1-s>0$,
\[
    \Im\zeta_0
    =
    -\frac{\chi(1-s)\Im z_0}
          {|z_0-1+s|^2}
    >0.
\]
Therefore $V_s(\zeta_0)=\chi$ has an upper-half-plane root and hence a
nonreal conjugate pair.

Conversely, suppose that $V_s(\zeta_0)=\chi$ for some
$\zeta_0\in\mathbb H$. The root is simple, so it extends to a
holomorphic solution $\zeta(w)$ in a complex neighborhood $\Omega$ of
$\chi$. On $\Omega\cap\mathbb H$, uniqueness gives
$\zeta(w)=\zeta_s(w)$. Set
\[
    z(w)
    :=
    e^{\mathrm{St}_{\bm_0}(\zeta(w))}.
\]
After shrinking $\Omega$, the function $(z(w)-1+s)/s$ is nonzero, and
there is a holomorphic logarithm $G$ on $\Omega$ that agrees with
$G_s$ on $\Omega\cap\mathbb H$. By
\eqref{eq:physical-spectral-branch-identity},
\[
    G(w)=\mathrm{St}_{\bm_y}(w),
    \qquad w\in\Omega\cap\mathbb H.
\]
Equation~\eqref{eq:St-m0-imaginary-range} gives
$z(\chi)\in\{\Im z<0\}$, and hence
\[
    -\pi<\Im G(\chi)<0.
\]
After shrinking to a real interval $I\ni\chi$,
Proposition~\ref{prop:density-existence} gives a continuous density
\[
    f_y(x)
    =
    -\frac{1}{\pi}\Im G(x),
    \qquad x\in I,
\]
satisfying $0<f_y(x)<1$ on $I$.

Lemma~\ref{lem:limit-measure-domain-support} gives
\[
    \operatorname{supp}(\bm_y)\subset\mathcal R_y.
\]
Since $f_y$ is strictly positive on $I$, we have
$I\subset\operatorname{supp}(\bm_y)$, and therefore
$I\subset\mathcal R_y$. Definition~\ref{df28} now shows that
$(\chi,y)\in\mathcal L$, or equivalently
$(\chi,s)\in\widetilde{\mathcal L}$.
\end{proof}

\begin{lemma}\label{lrd}
Assume that conditions
\ref{ass:tail-profile-convergence} and
\ref{ass:tail-profile-invertibility} of
Assumption~\ref{ass:tail-profile}, as well as
Assumption~\ref{ap32}, hold. Let
\[
    \mathcal T_{\mathcal L}:\widetilde{\mathcal L}\longrightarrow\mathbb H
\]
map $(\chi,s)$ to the unique upper-half-plane root of
$V_s(\zeta)=\chi$. Then $\mathcal T_{\mathcal L}$ is a
diffeomorphism, with inverse given by
\begin{align}
    &\chi_{\mathcal{L}}(\zeta):=\zeta\left[1+\frac{e^{-\mathrm{St}_{\bm_0}(\zeta)}(\zeta-\ol{\zeta})}{\ol{\zeta}e^{-\mathrm{St}_{\bm_0}(\ol{\zeta})}-\zeta e^{-\mathrm{St}_{\bm_0}(\zeta)}}\right]\label{dc}\\
    &s_{\mathcal{L}}(\zeta):=1+\frac{\zeta-\ol{\zeta}}{\ol{\zeta}e^{-\mathrm{St}_{\bm_0}(\ol{\zeta})}-\zeta e^{-\mathrm{St}_{\bm_0}(\zeta)}},\label{ds}
\end{align}
 where $\mathrm{St}_{\bm_0}$ represents the Stieltjes transform of the measure $\bm_0$.
\end{lemma}

\begin{proof}The proof is an adaptation of the proof of Theorem 2.1 in \cite{dm14}. We shall show that
\begin{enumerate}
    \item $\widetilde{\mathcal L}$ is nonempty,
    \item $\widetilde{\mathcal L}$ is open,
    \item $\mathcal T_{\mathcal L}$ is continuous,
    \item $\mathcal T_{\mathcal L}$ is injective,
    \item the formulas below give the inverse on $\mathcal T_{\mathcal L}(\widetilde{\mathcal L})$,
    \item $\mathcal T_{\mathcal L}(\widetilde{\mathcal L})=\mathbb H$.
\end{enumerate}
\bigskip
\textbf{Proof of (1).} Set
\begin{align*}
    \mathscr R(\zeta)
    :=
    \zeta\exp\bigl(-\mathrm{St}_{\bm_0}(\zeta)\bigr).
\end{align*}
Since $\bm_0$ is compactly supported, as $|\zeta|\to\infty$,
\begin{align}
    \mathrm{St}_{\bm_0}(\zeta)
    =
    \frac{1}{\zeta}
    +\frac{M_1}{\zeta^2}
    +\frac{M_2}{\zeta^3}
    +O\bigl(|\zeta|^{-4}\bigr),
    \label{est}
\end{align}
where
\begin{align*}
    M_1=\int_{\RR}x\bm_0[dx],
    \qquad
    M_2=\int_{\RR}x^2\bm_0[dx].
\end{align*}
Consequently,
\begin{align*}
    \mathscr R(\zeta)
    =
    \zeta-1
    +\frac{\frac12-M_1}{\zeta}
    +\frac{M_1-M_2-\frac16}{\zeta^2}
    +O\bigl(|\zeta|^{-3}\bigr).
\end{align*}
To prove nonemptiness, let $\zeta=\mathbf iR$ and let $R\to\infty$.
Then
\begin{align*}
    \Im\mathscr R(\mathbf iR)
    =
    R+\frac{M_1-\frac12}{R}+O(R^{-3}).
\end{align*}
By Schwarz symmetry, the denominator in \eqref{dc}--\eqref{ds} equals
\begin{align*}
    \overline{\mathscr R(\mathbf iR)}-\mathscr R(\mathbf iR)
    =
    -2\mathbf i\,\Im\mathscr R(\mathbf iR),
\end{align*}
and is therefore nonzero for all sufficiently large $R$.  Equations
\eqref{dc}--\eqref{ds} now give
\begin{align*}
    s_{\mathcal L}(\mathbf iR)
    &=
    1-
    \frac{R}{\Im\mathscr R(\mathbf iR)}
    =
    \frac{M_1-\frac12}{R^2}+O(R^{-4}),\\
    \chi_{\mathcal L}(\mathbf iR)
    &=
    -\frac{R\,\Re\mathscr R(\mathbf iR)}
          {\Im\mathscr R(\mathbf iR)}
    =
    1+O(R^{-2}).
\end{align*}
Because $\bm_0$ is a probability measure,
\begin{align*}
    \sum_{i=1}^{l}(b_i-a_i)=1.
\end{align*}
Writing $\ell_i=b_i-a_i$ and
$L_{i-1}=\ell_1+\cdots+\ell_{i-1}$, we have
$a_i\geq L_{i-1}$ for every $i$, with strict inequality for some
$i\geq2$. Hence
\begin{align*}
    M_1
    =
    \sum_{i=1}^{l}
    \left(a_i\ell_i+\frac{\ell_i^2}{2}\right)
    >
    \sum_{i=1}^{l}
    \left(L_{i-1}\ell_i+\frac{\ell_i^2}{2}\right)
    =
    \int_0^1x\,dx
    =
    \frac12.
\end{align*}
The same strict gap also gives $b_l>1$. Thus, for all sufficiently
large $R$,
\begin{align*}
    \bigl(
        \chi_{\mathcal L}(\mathbf iR),
        s_{\mathcal L}(\mathbf iR)
    \bigr)
    \in
    (0,b_l)\times(0,1).
\end{align*}
Finally, if
\begin{align*}
    q_R
    :=
    \frac{\mathbf iR-\overline{\mathbf iR}}
         {\overline{\mathscr R(\mathbf iR)}-\mathscr R(\mathbf iR)},
\end{align*}
then $s_{\mathcal L}(\mathbf iR)=1+q_R$ and
\begin{align*}
    \chi_{\mathcal L}(\mathbf iR)
    =
    \mathbf iR+q_R\mathscr R(\mathbf iR)
    =
    V_{s_{\mathcal L}(\mathbf iR)}(\mathbf iR).
\end{align*}
Since $\chi_{\mathcal L}(\mathbf iR)$ and
$s_{\mathcal L}(\mathbf iR)$ are real, Schwarz symmetry supplies the
conjugate lower-half-plane root.  The preceding lemma therefore gives
\begin{align*}
    \bigl(
        \chi_{\mathcal L}(\mathbf iR),
        s_{\mathcal L}(\mathbf iR)
    \bigr)
    \in\widetilde{\mathcal L}.
\end{align*}
Hence $\widetilde{\mathcal L}$ is nonempty.

\bigskip
\noindent\textbf{Proofs of (2) and (3).}
Let $(\chi_1,s_1)\in\widetilde{\mathcal L}$ and
$\zeta_1=\mathcal T_{\mathcal L}(\chi_1,s_1)$. Choose $\epsilon>0$ so small
that $\overline{B(\zeta_1,\epsilon)}\subset\mathbb H$ and
$V_{s_1}(\zeta)-\chi_1$ has exactly one zero in this disk and no zero
on its boundary. Set
\[
    m_\epsilon
    :=
    \min_{\zeta\in\partial B(\zeta_1,\epsilon)}
    |V_{s_1}(\zeta)-\chi_1|
    >0.
\]
By uniform continuity on the boundary circle, whenever
$|\chi_1-\chi_2|+|s_1-s_2|$ is sufficiently small,
\[
    |(V_{s_2}(\zeta)-\chi_2)-(V_{s_1}(\zeta)-\chi_1)|
    <m_\epsilon,
    \qquad
    \zeta\in\partial B(\zeta_1,\epsilon).
\]
Rouch\'e's theorem implies that $V_{s_2}(\zeta)-\chi_2$ has exactly
one zero in $B(\zeta_1,\epsilon)$. Hence
$(\chi_2,s_2)\in\widetilde{\mathcal L}$. Since $\epsilon$ may be
chosen arbitrarily small, the corresponding upper-half-plane root
converges to $\zeta_1$ as $(\chi_2,s_2)\to(\chi_1,s_1)$. This proves
both openness of $\widetilde{\mathcal L}$ and continuity of
$\mathcal T_{\mathcal L}$.

\bigskip
\noindent\textbf{Proofs of (4) and (5).}
Recall that
\begin{align*}
    \mathscr R(\zeta)
    =
    \zeta e^{-\mathrm{St}_{\bm_0}(\zeta)}.
\end{align*}
Under Assumption \ref{ap32}, equation \eqref{ee1} and $a_1=0$ give
\begin{align*}
    \mathscr R(\zeta)
    =
    \frac{\prod_{i=1}^{l}(\zeta-b_i)}
         {\prod_{i=2}^{l}(\zeta-a_i)}
    =
    \zeta+c+\sum_{j=2}^{l}\frac{A_j}{\zeta-a_j},
\end{align*}
where $c\in\RR$ and
\begin{align*}
    A_j
    =
    \frac{\prod_{i=1}^{l}(a_j-b_i)}
         {\prod_{\substack{2\leq i\leq l\\i\neq j}}(a_j-a_i)}
    <0,
    \qquad j=2,\ldots,l.
\end{align*}
Indeed, the numerator in the expression for $A_j$ has
$l-j+1$ negative factors, whereas the denominator has $l-j$ negative
factors.  Hence, for $\zeta=x+\mathbf i\eta\in\HH$,
\begin{align}
    \Im\mathscr R(\zeta)
    =
    \eta+
    \sum_{j=2}^{l}
    \frac{-A_j\eta}{(x-a_j)^2+\eta^2}
    >
    \eta.
    \label{eq:Im-R-positive}
\end{align}
Since $\mathrm{St}_{\bm_0}$ satisfies Schwarz symmetry, the common
denominator in \eqref{dc}--\eqref{ds} is
\begin{align*}
    D(\zeta)
    &:={}
    \overline{\zeta}
    e^{-\mathrm{St}_{\bm_0}(\overline{\zeta})}
    -
    \zeta e^{-\mathrm{St}_{\bm_0}(\zeta)}\\
    &=
    \overline{\mathscr R(\zeta)}-\mathscr R(\zeta)
    =
    -2\mathbf i\,\Im\mathscr R(\zeta),
\end{align*}
which is nonzero for every $\zeta\in\HH$.  Moreover,
\begin{align*}
    q(\zeta)
    :=
    \frac{\zeta-\overline\zeta}{D(\zeta)}
    =
    -\frac{\eta}{\Im\mathscr R(\zeta)}
    \in(-1,0)
\end{align*}
by \eqref{eq:Im-R-positive}. Thus
\begin{align*}
    s_{\mathcal L}(\zeta)=1+q(\zeta)\in(0,1).
\end{align*}
Also,
\begin{align*}
    \chi_{\mathcal L}(\zeta)
    =
    \zeta+q(\zeta)\mathscr R(\zeta),
\end{align*}
and therefore
\begin{align*}
    \Im\chi_{\mathcal L}(\zeta)
    =
    \eta+q(\zeta)\Im\mathscr R(\zeta)
    =0.
\end{align*}
Since $q(\zeta)=-(1-s_{\mathcal L}(\zeta))$, we further have
\begin{align*}
    V_{s_{\mathcal L}(\zeta)}(\zeta)
    &=
    \zeta\left[
        1-
        \bigl(1-s_{\mathcal L}(\zeta)\bigr)
        e^{-\mathrm{St}_{\bm_0}(\zeta)}
    \right]\\
    &=
    \zeta+q(\zeta)\mathscr R(\zeta)
    =
    \chi_{\mathcal L}(\zeta).
\end{align*}
Thus $\chi_{\mathcal L}(\zeta)\in\RR$, and Schwarz symmetry shows
that $\zeta$ and $\overline\zeta$ form a conjugate pair of roots of
\begin{align*}
    V_{s_{\mathcal L}(\zeta)}(w)
    =
    \chi_{\mathcal L}(\zeta).
\end{align*}
By the preceding lemma,
\begin{align*}
    \bigl(
        \chi_{\mathcal L}(\zeta),
        s_{\mathcal L}(\zeta)
    \bigr)
    \in\widetilde{\mathcal L}.
\end{align*}
Finally, if $(\chi,s)\in\widetilde{\mathcal L}$ and
$\zeta=\mathcal T_{\mathcal L}(\chi,s)$, then the equations
\begin{align*}
    V_s(\zeta)=\chi,
    \qquad
    V_s(\overline\zeta)=\chi
\end{align*}
form a nonsingular linear system for the real variables $\chi$ and
$s$, because its determinant is $D(\zeta)\neq0$. Solving this system
gives precisely \eqref{dc} and \eqref{ds}. Hence
\begin{align*}
    \chi_{\mathcal L}(\zeta)=\chi,
    \qquad
    s_{\mathcal L}(\zeta)=s,
\end{align*}
which proves (4) and (5). The simple-root implicit-function theorem
shows that $\mathcal T_{\mathcal L}$ is real analytic, while the explicit inverse
formulas \eqref{dc}--\eqref{ds} are real analytic on $\mathbb H$.
Thus the resulting homeomorphism is a diffeomorphism.

\bigskip
\noindent\textbf{Proof of (6).} From (1)--(5) we see that $\mathcal T_{\mathcal L}(\widetilde{\mathcal{L}})$ is open and homeomorphic to $\widetilde{\mathcal{L}}$. Assume there exists
\[
\zeta\in\partial \mathcal T_{\mathcal L}(\widetilde{\mathcal{L}})
\cap\bigl(\HH\setminus \mathcal T_{\mathcal L}(\widetilde{\mathcal{L}})\bigr).
\]
Choose $\zeta_n\in \mathcal T_{\mathcal L}(\widetilde{\mathcal{L}})$ with $\zeta_n\to\zeta$. Since the denominator in \eqref{dc}--\eqref{ds} is nonzero at $\zeta$, there exists $\delta>0$ such that
\[
K:=\overline{B(\zeta,\delta)}\subset\HH
\]
and this denominator is nonzero on $K$. Hence the map
\[
\eta\longmapsto
\bigl(\chi_{\mathcal{L}}(\eta),s_{\mathcal{L}}(\eta)\bigr)
\]
is continuous on the compact set $K$, so its image is compact. For all sufficiently large $n$, $\zeta_n\in K$; therefore a subsequence of
$\bigl(\chi_{\mathcal{L}}(\zeta_n),s_{\mathcal{L}}(\zeta_n)\bigr)$ converges to some $(\chi,s)$. By continuity of \eqref{dc}--\eqref{ds},
\[
(\chi,s)=\bigl(\chi_{\mathcal{L}}(\zeta),s_{\mathcal{L}}(\zeta)\bigr),
\]
and hence $\zeta=\mathcal T_{\mathcal L}(\chi,s)$. Since $\zeta\in\HH$, the preceding characterization gives $(\chi,s)\in\widetilde{\mathcal{L}}$, contradicting $\zeta\notin \mathcal T_{\mathcal L}(\widetilde{\mathcal{L}})$.
\end{proof}

\begingroup

\begin{lemma}[Contour-to-Green covariance identity]
\label{lem:contour-to-green}
Assume that conditions
\ref{ass:tail-profile-convergence} and
\ref{ass:tail-profile-invertibility} of
Assumption~\ref{ass:tail-profile}, as well as
Assumption~\ref{ap32}, hold. Let
\(\widetilde{\mathcal L}\) and
\(\mathcal T_{\mathcal L}\) be as in
\eqref{dll} and Lemma~\ref{lrd}, respectively. For \(s\in(0,1)\), set
\[
    I_s
    :=
    \{\chi\in\mathbb R:(\chi,s)\in\widetilde{\mathcal L}\},
    \qquad
    z_s(\chi)
    :=
    \mathcal T_{\mathcal L}(\chi,s),
    \quad \chi\in I_s.
\]

For \(0<s_1\leq s_2<1\) and \(k_1,k_2\geq1\), let
\(\mathfrak C_{s_1,s_2}^{k_1,k_2}\) denote the double contour integral
in \eqref{covf} after the substitutions
\(s_{r_1}=s_1\) and \(s_{r_2}=s_2\). Then
\begin{equation}
    \mathfrak C_{s_1,s_2}^{k_1,k_2}
    =
    \frac{k_1k_2}{\pi}
    \int_{I_{s_1}}\int_{I_{s_2}}
    \chi^{k_1-1}\eta^{k_2-1}
    G_{\mathbb H}
    \bigl(z_{s_1}(\chi),z_{s_2}(\eta)\bigr)
    \,d\eta\,d\chi.
    \label{eq:contour-to-green}
\end{equation}
\end{lemma}

\begin{proof}
We first identify the physical inverse branch of the weighted spectral
map and establish the nesting of its level-dependent spectral contours.
The final passage from the resulting contour integral to the Green
kernel uses the contour-to-Green argument of
\cite[Lemma~9.2 and equations~(9.4)--(9.6)]{bg16}.

Define the local meromorphic change of variables
\begin{equation}
    \mathscr A(x)
    :=
    xH_{\mathbf m_0}'(x)+\frac{x}{x-1}.
    \label{eq:A-transform-gff}
\end{equation}
Choose \(\operatorname{Log}x\) near \(x=1\) by
\(\operatorname{Log}1=0\).  Equation~\eqref{hmd} and the identity
\(\mathrm{St}_{\mathbf m_0}(u)=S_{\mathbf m_0}(u^{-1})\) give
\begin{equation}
    \mathscr A(x)
    =
    \frac{1}{S_{\mathbf m_0}^{(-1)}(\operatorname{Log}x)}
    =
    \mathrm{St}_{\mathbf m_0}^{(-1)}(\operatorname{Log}x).
    \label{eq:A-is-St-inverse}
\end{equation}
Here \(\mathrm{St}_{\mathbf m_0}^{(-1)}\) is the inverse branch at
infinity characterized by
\begin{equation}
    \mathrm{St}_{\mathbf m_0}(\mathscr A(x))
    =\operatorname{Log}x,
    \qquad
    \mathscr A(x)=\frac{1}{x-1}+O(1)
    \quad(x\to1).
    \label{eq:A-physical-branch}
\end{equation}

By the definition of \(Q\) following \eqref{covf},
\begin{align}
    Q(x,y)
    &=
    \frac{\partial^2}{\partial x\,\partial y}
    \log\bigl(\mathscr A(y)-\mathscr A(x)\bigr)
    \notag\\
    &=
    \frac{\mathscr A'(x)\mathscr A'(y)}
         {(\mathscr A(x)-\mathscr A(y))^2}.
    \label{eq:Q-under-A}
\end{align}
Indeed,
\begin{align*}
&(x-y)
\left(
1-
\frac{(x-1)(y-1)}{x-y}
\bigl[xH_{\mathbf m_0}'(x)-yH_{\mathbf m_0}'(y)\bigr]
\right)
\\
&\hspace{35mm}
=(x-1)(y-1)\bigl(\mathscr A(y)-\mathscr A(x)\bigr),
\end{align*}
and the factors depending on only one variable disappear after the
mixed derivative.  Moreover, if \(\zeta=\mathscr A(x)\), then
\(x=\exp(\mathrm{St}_{\mathbf m_0}(\zeta))\), and therefore
\begin{equation}
    (x-1+s)
    \left(
        H_{\mathbf m_0}'(x)+\frac{1}{x-1}
    \right)
    =V_s(\zeta).
    \label{eq:moment-factor-under-A}
\end{equation}

By symmetry of the covariance, assume first that \(s_1<s_2\) and
order the two pairs in \eqref{covf} accordingly.  By
\eqref{eq:A-physical-branch}, a positively oriented small circle about
\(1\) is mapped by \(\mathscr A\) to a clockwise large contour.  The
smaller circle \(|x-1|=\epsilon\) is mapped to the outer contour, while
\(|y-1|=2\epsilon\) is mapped to the inner one.

Set
\[
    \zeta=\mathscr A(x),
    \qquad
    \omega=\mathscr A(y).
\]
Equations~\eqref{eq:Q-under-A} and
\eqref{eq:moment-factor-under-A} transform the integrand in
\eqref{covf} into
\[
    V_{s_1}(\zeta)^{k_1}
    V_{s_2}(\omega)^{k_2}
    \frac{d\zeta\,d\omega}{(\zeta-\omega)^2}.
\]
Both image contours are clockwise; reversing both orientations
produces no net sign.
  For sufficiently small
\(\epsilon\), they are disjoint Jordan curves near infinity and may be
deformed, without crossing a pole or each other, to positively oriented
round circles while preserving their nesting.  Thus
\begin{equation}
    \mathfrak C_{s_1,s_2}^{k_1,k_2}
    =
    \frac{1}{(2\pi\mathbf i)^2}
    \oint_{|\zeta|=2C}
    \oint_{|\omega|=C}
    \frac{V_{s_1}(\zeta)^{k_1}V_{s_2}(\omega)^{k_2}}
         {(\zeta-\omega)^2}
    \,d\zeta\,d\omega,
    \label{eq:large-spectral-covariance}
\end{equation}
where both circles are positively oriented and \(C\) is sufficiently
large.  

We now deform the two large circular contours in
\eqref{eq:large-spectral-covariance} to the level contours
\(s_{\mathcal L}=s_1\) and \(s_{\mathcal L}=s_2\). We first show that
these contours enclose all finite poles and are nested in the required
order.  In the notation of
\eqref{eq:R-partial-fraction-physical-proof}, put \(c_j:=-A_j>0\) and
\begin{equation}
    \mathfrak q(\zeta)
    :=
    \sum_{j=2}^{l}\frac{c_j}{|\zeta-a_j|^2}.
    \label{eq:spectral-nesting-function}
\end{equation}
with the convention \(\mathfrak q(a_j)=+\infty\).
The partial-fraction expansion there reads
\[
    \mathscr R(\zeta)
    :=
    \zeta e^{-\mathrm{St}_{\mathbf m_0}(\zeta)}
    =
    \zeta+c-\sum_{j=2}^{l}\frac{c_j}{\zeta-a_j}.
\]
For \(\zeta=x+\mathbf i v\in\mathbb H\), it follows that
\[
    \Im\mathscr R(\zeta)
    =v\bigl(1+\mathfrak q(\zeta)\bigr).
\]
Formula~\eqref{ds} consequently gives
\begin{equation}
    s_{\mathcal L}(\zeta)
    =
    1-\frac{v}{\Im\mathscr R(\zeta)}
    =
    \frac{\mathfrak q(\zeta)}{1+\mathfrak q(\zeta)}.
    \label{eq:spectral-level-function}
\end{equation}
For \(s\in(0,1)\), define
\[
    \Omega_s
    :=
    \left\{
        \zeta\in\mathbb C:
        \mathfrak q(\zeta)>\frac{s}{1-s}
    \right\},
    \qquad
    \mathscr C_s:=\partial\Omega_s,
\]
where each component of \(\mathscr C_s=\partial\Omega_s\) is oriented
so that \(\Omega_s\) lies to its left.  Since
\[
    \mathfrak q(\zeta)\longrightarrow0
    \quad\text{as }|\zeta|\to\infty,
    \qquad
    \mathfrak q(\zeta)\longrightarrow\infty
    \quad\text{as }\zeta\to a_j,
\]
the set \(\Omega_s\) is bounded and contains a neighborhood of every
finite pole \(a_j\) of \(V_s\).  Its upper boundary is
\[
    \mathscr C_s^+
    :=
    \mathscr C_s\cap\mathbb H
    =
    \{\zeta\in\mathbb H:s_{\mathcal L}(\zeta)=s\}.
\]
The function \(\mathfrak q\) is real analytic on \(\mathbb H\), and,
for \(\zeta=x+\mathbf i v\) with \(v>0\),
\[
    \frac{\partial\mathfrak q}{\partial v}(x+\mathbf i v)
    =
    -2v\sum_{j=2}^{l}
    \frac{c_j}{\bigl((x-a_j)^2+v^2\bigr)^2}
    <0.
\]
Thus the gradient of \(\mathfrak q\) does not vanish on
\[
    \mathscr C_s^+
    =
    \left\{
        \zeta\in\mathbb H:
        \mathfrak q(\zeta)=\frac{s}{1-s}
    \right\}.
\]
By the real-analytic implicit function theorem, every component of
\(\mathscr C_s^+\) is locally the graph of a real-analytic function;
in particular, \(\mathscr C_s^+\) is a real-analytic curve.
Since \(s\mapsto s/(1-s)\) is strictly increasing,
\begin{equation}
    0<s_1<s_2<1
    \quad\Longrightarrow\quad
    \overline{\Omega_{s_2}}\subset\Omega_{s_1}.
    \label{eq:nested-spectral-domains}
\end{equation}
Thus \(\mathscr C_{s_2}\) is strictly inside
\(\mathscr C_{s_1}\), as an oriented boundary chain.
Since \(\Omega_s\) contains all finite poles of \(V_s\), its oriented
boundary chain is homologous to a sufficiently large positively
oriented circle in the complement of those poles.  If \(\partial\Omega_s\) is not smooth at a point of the real axis, we
interpret \(\mathscr C_s\) as the oriented boundary of \(\Omega_s\),
obtained as the limit of nearby smooth level curves.  By Cauchy's
theorem, replace the inner circle in
\eqref{eq:large-spectral-covariance} by the homologous boundary chain
\(\mathscr C_{s_2}\), and then replace the outer circle by
\(\mathscr C_{s_1}\).  The strict nesting in
\eqref{eq:nested-spectral-domains} ensures that this replacement
crosses neither a pole nor the diagonal \(\zeta=\omega\).  Hence
\begin{equation}
    \mathfrak C_{s_1,s_2}^{k_1,k_2}
    =
    \frac{1}{(2\pi\mathbf i)^2}
    \oint_{\mathscr C_{s_1}}
    \oint_{\mathscr C_{s_2}}
    \frac{V_{s_1}(z)^{k_1}V_{s_2}(w)^{k_2}}
         {(z-w)^2}
    \,dz\,dw.
    \label{eq:covariance-on-level-contours}
\end{equation}

On \(\mathscr C_s\), the inverse formulas
\eqref{dc}--\eqref{ds} give
\begin{equation}
    V_s(z)=\chi_{\mathcal L}(z)\in\mathbb R.
    \label{eq:V-on-level-curve}
\end{equation}
On each positively oriented upper component of \(\mathscr C_s\),
\[
    \Im V_s(z)
    =
    \Im z\left(s-(1-s)\mathfrak q(z)\right)
\]
is negative on the interior side and positive on the exterior side.
The Cauchy--Riemann equations therefore show that
\(\chi_{\mathcal L}=\Re V_s\) decreases along the boundary. We
parametrize both upper contours by increasing \(\chi_{\mathcal L}\);
the two orientation reversals cancel in the double integral.

For \(z,w\in\mathbb H\), set
\begin{equation}
    \mathscr K(z,w)
    :=
    \log
    \frac{(z-w)(\overline z-\overline w)}
         {(z-\overline w)(\overline z-w)}
    =
    2\log\left|\frac{z-w}{z-\overline w}\right|,
    \label{eq:green-cross-ratio}
\end{equation}
where \(\log\) denotes the real natural logarithm.  Writing
\[
    \mathscr C_s
    =
    \mathscr C_s^+\cup\overline{\mathscr C_s^+},
    \qquad
    \mathscr C_s^+=z_s(I_s),
\]
the double contour integral splits into the upper--upper,
upper--lower, lower--upper, and lower--lower contributions.
Parametrizing the upper curves by \(z_s(\chi)\), the lower curves by
\(\overline{z_s(\chi)}\), and using
\eqref{eq:V-on-level-curve}, their sum is
\begin{align}
    \mathfrak C_{s_1,s_2}^{k_1,k_2}
    ={}&
    \frac{1}{(2\pi\mathbf i)^2}
    \int_{I_{s_1}}\int_{I_{s_2}}
    \chi^{k_1}\eta^{k_2}
    \frac{\partial^2}{\partial\chi\,\partial\eta}
    \mathscr K
    \bigl(z_{s_1}(\chi),z_{s_2}(\eta)\bigr)
    \,d\eta\,d\chi.
    \label{eq:closed-contours-to-cross-ratio}
\end{align}

We integrate by parts once in \(\chi\) and once in \(\eta\).  All
boundary terms vanish.  Indeed, every endpoint of a component of
\(I_s\) is mapped by \(z_s\) to the real boundary, and
\[
    \mathscr K(x,w)=0,
    \qquad
    \mathscr K(z,y)=0,
    \qquad x,y\in\mathbb R,\quad z,w\in\mathbb H.
\]
In the first integration, for example,
\(\partial_\eta\mathscr K(x,z_{s_2}(\eta))=0\) at a real endpoint.
Therefore
\begin{align*}
    \mathfrak C_{s_1,s_2}^{k_1,k_2}
    ={}&
    \frac{k_1k_2}{(2\pi\mathbf i)^2}
    \int_{I_{s_1}}\int_{I_{s_2}}
    \chi^{k_1-1}\eta^{k_2-1}
    \mathscr K\bigl(z_{s_1}(\chi),z_{s_2}(\eta)\bigr)
    \,d\chi\,d\eta
    \\
    ={}&
    \frac{k_1k_2}{\pi}
    \int_{I_{s_1}}\int_{I_{s_2}}
    \chi^{k_1-1}\eta^{k_2-1}
    G_{\mathbb H}\bigl(z_{s_1}(\chi),z_{s_2}(\eta)\bigr)
    \,d\chi\,d\eta.
\end{align*}
The last equality uses
\[
    \mathscr K(z,w)=-4\pi G_{\mathbb H}(z,w),
    \qquad
    \frac{1}{(2\pi\mathbf i)^2}=-\frac{1}{4\pi^2}.
\]
This proves \eqref{eq:contour-to-green} when \(s_1<s_2\); the case
\(s_2<s_1\) follows by symmetry.

It remains to consider \(s_1=s_2=s\). Apply the identity already
proved with \(s<\sigma\) and let \(\sigma\downarrow s\). The contour
side is continuous in \(\sigma\).

Choose a compact interval \(J\) containing \(I_s\) and all
\(I_\sigma\) for \(\sigma\) sufficiently close to \(s\). Regard each
integrand as a function on \(J^2\) by setting it equal to zero outside
\(I_s\times I_\sigma\). The relevant root branches form a compact semialgebraic family by
\eqref{eq:spectral-polynomial}. Since
\(\mathcal T_{\mathcal L}\) is injective,
\[
    z_s(\chi)=z_\sigma(\eta)
    \quad\Longrightarrow\quad
    s=\sigma,\qquad \chi=\eta.
\]
The Łojasiewicz inequality on this compact family therefore gives
constants \(C_0>0\) and \(N\geq1\) such that
\[
    |z_s(\chi)-z_\sigma(\eta)|
    \geq C_0|\chi-\eta|^N
\]
uniformly for \(\sigma\) near \(s\). Hence
\[
    G_{\mathbb H}
    \bigl(z_s(\chi),z_\sigma(\eta)\bigr)
    \leq
    C\bigl(1+|\log|\chi-\eta||\bigr),
\]
and the right-hand side is integrable on \(J^2\). The extended
integrands converge pointwise almost everywhere as
\(\sigma\downarrow s\). Dominated convergence therefore proves
\eqref{eq:contour-to-green} for \(s_1=s_2\).
\end{proof}
\endgroup

\begin{theorem}
\label{thm:gff}
Suppose that the data $t_n$, $\mathbf x^{(n)}$, and
$\lambda^{(0)}(n)$ satisfy the standing hypotheses of
Theorem~\ref{t15}. In particular, the bottom-boundary sequence is
regular in the sense of Definition~\ref{def:regular-bottom}. Assume in
addition that
Assumption~\ref{ass:tail-profile}\ref{ass:tail-profile-invertibility}
and Assumption~\ref{ap32} hold. For $s\in(0,1)$ and all sufficiently large
$n$, set
\begin{equation}
    \kappa_n(s)
    :=
    t_n-\lfloor\vartheta(s)n\rfloor,
    \label{eq:kappa-general-tail}
\end{equation}
where $\vartheta$ is defined in
\eqref{eq:theta-tail-profile}. Define
\begin{equation}
    \Delta_{n,s}(\chi)
    :=
    \sqrt{\pi}
    \left|
        \left\{
            g\in[n]:
            \lambda_g^{(\kappa_n(s))}+n-g\geq n\chi
        \right\}
    \right|,
    \qquad \chi\in\mathbb R.
    \label{defDelta}
\end{equation}
The particles on the slice $\kappa_n(s)$ are located at
\[
    x_g^{(n,s)}
    :=
    \lambda_g^{(\kappa_n(s))}+n-g,
    \qquad g\in[n].
\]
Whenever $n\chi$ is not a particle position,
\[
    h^{\mathrm{disc}}(n\chi,\kappa_n(s))
    =
    \left|
        \left\{g\in[n]:x_g^{(n,s)}<n\chi\right\}
    \right|,
\]
and hence
\begin{equation}
    \Delta_{n,s}(\chi)
    =
    \sqrt\pi\left(
        n-h^{\mathrm{disc}}(n\chi,\kappa_n(s))
    \right)
    \label{eq:Delta-height}
\end{equation}
away from the finitely many jump points.

For $s\in(0,1)$, let
\[
    I_s
    :=
    \{\chi\in\mathbb R:(\chi,s)\in\widetilde{\mathcal L}\},
    \qquad
    z_s(\chi)
    :=
    \mathcal T_{\mathcal L}(\chi,s),
    \quad \chi\in I_s.
\]
For each integer $j\geq0$, define the finite signed measure on
$\mathbb H$
\begin{equation}
    \nu_{s,j}
    :=
    (z_s)_*\bigl(\chi^j\mathbf 1_{I_s}(\chi)\,d\chi\bigr).
    \label{eq:gff-curve-measure}
\end{equation}
Under Assumption~\ref{ap32}, the horizontal liquid slices are bounded.
Moreover, the logarithmic diagonal singularity of $G_{\mathbb H}$ is
integrable along the smooth curves $z_s(I_s)$, while the Green kernel
vanishes at the real boundary. Hence
\[
    \int_{\mathbb H}\int_{\mathbb H}
    |G_{\mathbb H}(z,w)|\,|\nu_{s,j}|(dz)\,|\nu_{s,j}|(dw)
    <\infty,
\]
so $\Xi(\nu_{s,j})$ is well-defined by
\eqref{eq:gff-measure-covariance}.

For every $m\geq1$, every $s_1,\ldots,s_m\in(0,1)$, and every
$j_1,\ldots,j_m\in\mathbb Z_{\geq0}$, the random vector
\begin{equation}
    \left(
        \int_{\mathbb R}\chi^{j_r}
        \left(
            \Delta_{n,s_r}(\chi)
            -\mathbb E\Delta_{n,s_r}(\chi)
        \right)d\chi
    \right)_{r=1}^{m}
    \label{eq:gff-observable-vector}
\end{equation}
converges, in the sense of moments, to the centered Gaussian vector
\begin{equation}
    \left(
        \Xi(\nu_{s_r,j_r})
    \right)_{r=1}^{m}.
    \label{eq:gff-limit-vector}
\end{equation}
Equivalently, these horizontal polynomial observables converge to the
corresponding observables of the pullback of the zero-boundary GFF under
$\mathcal T_{\mathcal L}$.
\end{theorem}

\begin{proof}
Fix $m$, $s_1,\ldots,s_m$, and $j_1,\ldots,j_m$ as in the theorem. For
each fixed $s\in(0,1)$, equations
\eqref{eq:kappa-general-tail} and
\eqref{eq:theta-tail-profile} give
\[
    \frac{\kappa_n(s)}{n}
    \longrightarrow
    \alpha-\vartheta(s)
    =y(s).
\]
Moreover, by
Assumption~\ref{ass:tail-profile}\ref{ass:tail-profile-convergence},
\begin{equation}
    \frac{|\mathbf x^{(n)}_{\kappa_n(s)}|}
         {|\mathbf x^{(n)}|}
    =
    s_n\left(\frac{\kappa_n(s)}{n}\right)
    \longrightarrow
    s(y(s))
    =s.
    \label{eq:tail-ratio-general-gff}
\end{equation}

Integration by parts for the step function $\Delta_{n,s}$ gives
\begin{align}
&\int_{\mathbb R}\chi^j
\left(
    \Delta_{n,s}(\chi)-\mathbb E\Delta_{n,s}(\chi)
\right)d\chi
\label{cx}\\
&\quad=
\frac{\sqrt\pi}{(j+1)n^{j+1}}
\left(
    p_{j+1}^{(\lfloor\vartheta(s)n\rfloor)}
    -\mathbb E p_{j+1}^{(\lfloor\vartheta(s)n\rfloor)}
\right).
\label{dlp}
\end{align}
Indeed, $\Delta_{n,s}$ decreases by $\sqrt\pi$ at each point
$\chi=n^{-1}(\lambda_g^{(\kappa_n(s))}+n-g)$, and the boundary term
vanishes because the centered step function is zero outside a bounded
interval.

Set $k_r:=j_r+1$ and
$\alpha_r:=\vartheta(s_r)\in(0,\alpha)$. By
\eqref{eq:tail-ratio-general-gff}, the hypotheses of
Theorem~\ref{t31} hold for these levels. Equations
\eqref{dlp} and \eqref{crv} therefore imply that the vector in
\eqref{eq:gff-observable-vector} has a joint Gaussian moment limit. It
remains to identify its covariance.

\begingroup

For two indices \(a,b\in[m]\), apply
Lemma~\ref{lem:contour-to-green} to the covariance in
Theorem~\ref{t31}, with \((s_1,k_1)=(s_a,k_a)\) and
\((s_2,k_2)=(s_b,k_b)\).  We obtain
\begin{align}
&\lim_{n\to\infty}
\frac{
\operatorname{Cov}\left(
    p_{k_a}^{(\lfloor\vartheta(s_a)n\rfloor)},
    p_{k_b}^{(\lfloor\vartheta(s_b)n\rfloor)}
\right)
}{n^{k_a+k_b}}
\label{ccov}\\
&\quad=
\frac{k_ak_b}{\pi}
\int_{z\in\mathbb H:\,s_{\mathcal L}(z)=s_a}
\int_{w\in\mathbb H:\,s_{\mathcal L}(w)=s_b}
\chi_{\mathcal L}(z)^{k_a-1}
\chi_{\mathcal L}(w)^{k_b-1}
G_{\mathbb H}(z,w)\,
 d\chi_{\mathcal L}(z)\,d\chi_{\mathcal L}(w).
\notag
\end{align}
The level curves are oriented so that \(\chi_{\mathcal L}\) increases.
The proof of Lemma~\ref{lem:contour-to-green} also shows explicitly
that the two height prefactors in \eqref{dlp} cancel the factor
\(k_ak_b/\pi\): indeed,
\[
    \frac{\sqrt\pi}{k_a}\,
    \frac{\sqrt\pi}{k_b}\,
    \frac{k_ak_b}{\pi}=1.
\]
\endgroup
Multiplying \eqref{ccov} by the two prefactors in \eqref{dlp} and using
$k_r=j_r+1$, we obtain
\begin{align*}
&\lim_{n\to\infty}
\operatorname{Cov}\left(
    \int_{\mathbb R}\chi^{j_a}
    (\Delta_{n,s_a}-\mathbb E\Delta_{n,s_a})\,d\chi,
    \int_{\mathbb R}\eta^{j_b}
    (\Delta_{n,s_b}-\mathbb E\Delta_{n,s_b})\,d\eta
\right)\\
&\quad=
\int_{I_{s_a}}\int_{I_{s_b}}
\chi^{j_a}\eta^{j_b}
G_{\mathbb H}\bigl(z_{s_a}(\chi),z_{s_b}(\eta)\bigr)
\,d\chi\,d\eta\\
&\quad=
\mathbb E\left[
    \Xi(\nu_{s_a,j_a})\Xi(\nu_{s_b,j_b})
\right].
\end{align*}
Theorem~\ref{t31} supplies Wick's rule for every joint moment. Hence the
vector in \eqref{eq:gff-observable-vector} converges in moments to the
Gaussian vector in \eqref{eq:gff-limit-vector}.
\end{proof}

\begin{corollary}[Uniform weights]
\label{cor:gff-uniform}
Suppose that the data $t_n$, $\mathbf x^{(n)}$, and
$\lambda^{(0)}(n)$ satisfy the standing hypotheses of
Theorem~\ref{t15} and Assumption~\ref{ap32}, and assume the uniform
specialization \eqref{uw}. Then
condition~\ref{ass:tail-profile-invertibility} of
Assumption~\ref{ass:tail-profile} holds, with
\[
    s(y)=1-\frac{y}{\alpha},
    \qquad
    y(s)=\alpha(1-s),
    \qquad
    \vartheta(s)=\alpha s.
\]
Consequently, the conclusion of Theorem~\ref{thm:gff} holds with
\[
    \kappa_n(s)=t_n-\lfloor\alpha sn\rfloor.
\]
\end{corollary}

\begin{proof}
Under \eqref{uw}, the discrete tail profile in
\eqref{eq:discrete-tail-profile} is
\[
    s_n(y)
    =
    \frac{(t_n-\lfloor ny\rfloor)_+}{t_n},
\]
and it converges uniformly on $[0,\alpha]$ to
$1-y/\alpha$. The result follows from
Theorem~\ref{thm:gff}.
\end{proof}

\appendix
\section{Technical Results}\label{sa}

We use $\YY$ to denote the set of all   partitions and $\YY_N$ to denote the set of all   partitions of length $N$.  The symbol $N$ is used in this appendix for the general length parameter.

\begin{lemma}\label{la1}
Let $\{\lambda(N)\}_{N\geq1}$, with $\lambda(N)\in\YY_N$, be regular in the sense of
Definition~\ref{def:regular-bottom}, and let $\bm$ be its associated
compactly supported probability measure.
Then there exists an explicit function $H_{\bm}$, analytic in a
neighborhood of $1$ and depending only on $\bm$, such that, for every
fixed $k\geq1$,
\begin{equation}
  \lim_{N\rightarrow\infty}
  \frac{1}{N} \log\left(%
  \frac{s_{\lambda(N)}(u_1,\ldots,u_k,1,\ldots,1)}{s_{\lambda(N)}(1,\ldots,1)}
  \right) = H_{\bm}(u_1)+\cdots+H_{\bm}(u_k),
  \label{nlc}
\end{equation}
and the convergence is locally uniform in a complex neighborhood of
$(1,\dotsc,1)$.
\end{lemma}

\begin{proof}
Lemma~\ref{lem:regular-implies-moment-convergence} gives
$m(\lambda(N))\Longrightarrow\bm$. The assertion is then
Theorem~8.1 of \cite{bg16}, the unitary specialization of
Theorem~4.2 of \cite{bg}.
\end{proof}

Precisely, $H_{\bm}$ is constructed as follows.  Recall from Definition
\ref{def:transforms} the moment generating function
$S_{\bm}(z)=z+\sum_{k=1}^\infty M_k(\bm)z^{k+1}$ and its local compositional
inverse $S_{\bm}^{(-1)}$.  Let $R_{\bm}(z)$ be the
\emph{Voiculescu R-transform} of $\bm$ defined as
\begin{equation*}
  R_{\bm}(z) = \frac{1}{S_\bm^{(-1)}(z)} - \frac{1}{z}.
\end{equation*}
Then
\begin{equation}
  \label{hmz}
  H_{\bm}(u) = \int_{0}^{\ln u} R_\bm(t)dt+ \ln\left( \frac{\ln u}{u-1} \right).
\end{equation}
In particular, $H_{\bm}(1)=0$, and
\begin{equation}
  H'_\bm(u) = \frac{1}{u S_\bm^{(-1)}(\ln u)} - \frac{1}{u-1}.\label{hmd}
\end{equation}
Here all logarithms are taken in their holomorphic branches near
\(u=1\), normalized by continuity at \(u=1\); these branches agree
locally with the principal logarithm.

\begin{lemma}\label{la2}Let $n$ be a positive integer and let $g(z)$ be an analytic function defined in a neighborhood of 1. Then
\begin{eqnarray*}
\lim_{(z_1,\ldots,z_n)\rightarrow (1,\ldots,1)}\left(\sum_{i=1}^n\frac{g(z_i)}{\prod_{j\in[n],j\neq i}(z_i-z_j)}\right)=\left.\frac{\partial^{n-1}}{\partial z^{n-1}}\left(\frac{g(z)}{(n-1)!}\right)\right|_{z=1}.
\end{eqnarray*}
\end{lemma}

\begin{proof}See Lemma 5.5 of \cite{bg}.
\end{proof}

\begin{lemma}[Partial confluence of divided differences]
\label{lem:partial-confluent-divided-difference}
Let \(r\geq0\), and let \(h\) be analytic near \(1\), with
\[
    h(z)=\sum_{m\geq0}c_m(z-1)^m.
\]
Keeping \(z_1\) free, one has
\begin{align}
&\lim_{z_2,\ldots,z_{r+1}\to1}
\sum_{i=1}^{r+1}
\frac{h(z_i)}
{\prod_{\substack{1\leq j\leq r+1\\j\neq i}}(z_i-z_j)}
=
\frac{
    h(z_1)-\displaystyle\sum_{m=0}^{r-1}c_m(z_1-1)^m
}{
    (z_1-1)^r
}.
\label{eq:partial-confluent-divided-difference}
\end{align}
For \(r=0\), the sum in the numerator is empty. The right-hand side
extends analytically to \(z_1=1\).
\end{lemma}

\begin{proof}
For \(r=0\), the identity is immediate. For \(r\geq1\), this is
\cite[Lemma~5.4]{bg16}.
\end{proof}

\begin{lemma}\label{la3}
Let $\{\lambda(N)\}_{N\geq1}$, with $\lambda(N)\in\YY_N$, be regular in the sense of
Definition~\ref{def:regular-bottom}, and let $\bm$ be its associated
compactly supported probability measure.
Then, for every fixed $k\geq3$,
\begin{align*}
    &\lim_{N\rightarrow\infty}\frac{\partial^2}{\partial x_1\partial x_2}\log\frac{s_{\lambda(N)}(x_1,\ldots,x_{k},1^{N-k})}{s_{\lambda(N)}(1^N)}
    \\
    &=\frac{\partial^2}{\partial x_1\partial x_2}\log\left(1-(x_1-1)(x_2-1)\frac{x_1 H_{\bm}'(x_1)-x_2H_{\bm}'(x_2)}{x_1-x_2}\right)
\end{align*}
and
\begin{align*}
    \lim_{N\rightarrow\infty}\frac{\partial^3}{\partial x_1\partial x_2\partial x_3}\log\frac{s_{\lambda(N)}(x_1,\ldots,x_{k},1^{N-k})}{s_{\lambda(N)}(1^N)}=0.
\end{align*}
Both convergences are locally uniform in a complex neighborhood of
$(1,\ldots,1)$.
\end{lemma}

\begin{proof}
Lemma~\ref{lem:regular-implies-moment-convergence} gives
$m(\lambda(N))\Longrightarrow\bm$, so the assertion follows from
Theorem~8.2 of \cite{bg16}.
\end{proof}

\section{A Detailed Frozen-Boundary Example}
\label{app:linear-bottom-boundary}

\begin{example}\label{ex:linear-bottom-boundary}
Assume that the bottom boundary partition is
\[
    \lambda^{(0)}(n)
    =
    \bigl((p-1)n,(p-1)(n-1),\ldots,p-1\bigr)
    \in\YY_n,
\]
where $p$ is a fixed positive integer and $n\geq1$. This sequence is regular in
the sense of Definition~\ref{def:regular-bottom}, with limiting
bottom-boundary profile
\[
    f_0(u)=(p-1)(1-u),
    \qquad 0\leq u\leq1.
\]

The limiting bottom-boundary
measure is the uniform probability measure on $(0,p)$, namely

\[
     \varrho_0(x)
    :=
    \frac{d\mathbf m_0}{dx}(x)
    =
    \begin{cases}
        \dfrac1p, & x\in(0,p),\\[1mm]
        0, & \text{otherwise}.
    \end{cases}
\]

Consequently,
\[
    M_k(\mathbf m_0)
    =
    \int_{\mathbb R}x^k\,\mathbf m_0(dx)
    =
    \frac{p^k}{k+1},
    \qquad k\geq0.
\]
Therefore, for the moment series introduced in
Definition~\ref{def:transforms},
\[
    S_{\mathbf m_0}(z)
    =
    \sum_{k=0}^{\infty}M_k(\mathbf m_0)z^{k+1}
    =
    -\frac1p\log(1-pz),
\]
and hence
\[
    S_{\mathbf m_0}^{(-1)}(u)
    =
    \frac{1-e^{-pu}}{p}.
\]
Applying \eqref{hmd} with $\mathbf m=\mathbf m_0$ gives
\[
    H'_{\mathbf m_0}(u)
    =
    \frac{pu^{p-1}}{u^p-1}
    -
    \frac1{u-1}.
\]
Thus, for fixed $y$ and the corresponding value $s=s(y)$,
\[
    U_y(z)
    =
    \frac{pz^{p-1}(z-1+s)}{z^p-1}.
\]

We now specialize to $p=3$ and write
\[
    U_s(z)
    :=
    \frac{3z^2(z-1+s)}{z^3-1},
    \qquad 0<s<1.
\]
The equation $U_s(z)=\chi$ is equivalent to
\begin{equation}
    P_{\chi,s}(z)
    :=
    (3-\chi)z^3
    -3(1-s)z^2
    +\chi
    =0.
\end{equation}
Its discriminant is
\begin{equation}
    \Delta_s(\chi)
    =
    27\chi
    \left[
        4(1-s)^3
        -\chi(3-\chi)^2
    \right].
\end{equation}
Since the function
\[
    \chi\longmapsto \chi(3-\chi)^2
\]
is strictly increasing on $(0,1)$, strictly decreasing on $(1,3)$,
and has maximum value $4$ at $\chi=1$, there are unique numbers
\[
    \chi_-(s)\in(0,1),
    \qquad
    \chi_+(s)\in(1,3),
\]
satisfying
\begin{equation}
    \chi_\pm(s)\bigl(3-\chi_\pm(s)\bigr)^2
    =
    4(1-s)^3.
\end{equation}
It follows from the discriminant formula that, for $0<\chi<3$, the
equation $U_s(z)=\chi$ has exactly one real root and one nonreal
complex-conjugate pair precisely when
\[
    \chi_-(s)<\chi<\chi_+(s).
\]

\medskip
\noindent\textit{Auxiliary lemma (identification of the Stieltjes branch).}
For every
\[
    \chi\in\bigl(\chi_-(s),\chi_+(s)\bigr),
\]
the solution branch $\mathbf z_1(w,y)$ introduced in
Proposition~\ref{prop:stieltjes-root} and characterized by
\[
    \mathbf z_1(w,y)\longrightarrow1
    \qquad\text{as }w\longrightarrow\infty
\]
continues through the upper half of the $w$-plane to the unique root of
$U_s(z)=\chi$ in the lower half of the $z$-plane. In other words, the
boundary value
\[
    \mathbf z_1(\chi+\mathbf i0,y)
    :=
    \lim_{\delta\downarrow0}
    \mathbf z_1(\chi+\mathbf i\delta,y)
\]
exists and satisfies
\[
    \Im\mathbf z_1(\chi+\mathbf i0,y)<0.
\]

\begin{proof}[Proof of the auxiliary lemma]
Set $a:=1-s\in(0,1)$. A direct differentiation gives
\begin{equation}
    U_s'(z)
    =
    \frac{3z\bigl(az^3+2a-3z\bigr)}{(z^3-1)^2}.
\end{equation}
Let
\[
    g(z):=az^3+2a-3z.
\]
Then
\[
    g(1)=-3s<0,
    \qquad
    g'(z)=3(az^2-1).
\]
Since $a^{-1/2}>1$, the function $g$ is strictly decreasing on
$(1,a^{-1/2})$ and strictly increasing on
$(a^{-1/2},\infty)$. Moreover,
\[
    g(a^{-1/2})=2a-2a^{-1/2}<0,
    \qquad
    \lim_{z\to\infty}g(z)=+\infty.
\]
Hence $g$ has exactly one zero in $(1,\infty)$. Denote this zero by
$r_+(s)$. In fact,
\[
    r_+(s)>a^{-1/2}>1.
\]
Consequently, $U_s$ is strictly decreasing on $(1,r_+(s))$, strictly
increasing on $(r_+(s),\infty)$, and
\[
    \chi_+(s)=U_s(r_+(s)).
\]
At the critical point,
\begin{equation}
    U_s''(r_+(s))
    =
    \frac{18r_+(s)}{
        (r_+(s)-1)(r_+(s)^3+2)(r_+(s)^2+r_+(s)+1)
    }
    >0.
\end{equation}

Since
\[
    U_s(z)=\frac{s}{z-1}+O(1)
    \qquad\text{as }z\to1,
\]
the solution branch selected at infinity satisfies
\[
    \mathbf z_1(w,y)
    =
    1+\frac{s}{w}+O(w^{-2}).
\]
We next justify that this branch is real for every real
$w>\chi_+(s)$. 

Since \(U_s\) has real coefficients, the conjugate branch
\[
    \widehat{\mathbf z}_1(w,y)
    :=
    \overline{\mathbf z_1(\overline w,y)}
\]
also satisfies
\[
    U_s\bigl(\widehat{\mathbf z}_1(w,y)\bigr)=w.
\]
Moreover, both \(\widehat{\mathbf z}_1\) and \(\mathbf z_1\) have the
same expansion
\[
    1+\frac{s}{w}+O(w^{-2})
\]
as \(w\to\infty\). By uniqueness of the local inverse branch near
\(z=1\), they agree for sufficiently large \(\lvert w\rvert\).
Therefore, \(\mathbf z_1(w,y)\) is real for sufficiently large
positive real \(w\), and the expansion shows that
\(\mathbf z_1(w,y)>1\).

Moreover,
\[
    U_s:(1,r_+(s))\longrightarrow(\chi_+(s),\infty)
\]
is a strictly decreasing bijection, and $U_s'$ does not vanish on
$(1,r_+(s))$. Its real inverse therefore agrees with
$\mathbf z_1(w,y)$ for all sufficiently large positive $w$ and, by
local uniqueness and continuation along the interval, for every
$w>\chi_+(s)$. Thus
\[
    \mathbf z_1(w,y)\in(1,r_+(s)),
    \qquad w>\chi_+(s).
\]

Near the critical value $\chi_+(s)$, the local inverse has the
square-root expansion
\begin{equation}
    \mathbf z_1(w,y)
    =
    r_+(s)
    -
    \sqrt{
        \frac{2\bigl(w-\chi_+(s)\bigr)}
             {U_s''(r_+(s))}
    }
    +
    O\bigl(w-\chi_+(s)\bigr),
\end{equation}
where the square root is positive for real $w>\chi_+(s)$. Let
\[
    a:=\chi_+(s)-\chi>0,
    \qquad
    c_s:=
    \sqrt{\frac{2}{U_s''(r_+(s))}}>0.
\]
The square root in the preceding expansion is the branch that is
positive for real \(w>\chi_+(s)\). Hence, as \(w=\chi+\mathbf i\delta\)
approaches the real axis from the upper half-plane,
\[
    \sqrt{w-\chi_+(s)}
    \longrightarrow
    \sqrt{-a+\mathbf i0}
    =
    \mathbf i\sqrt a.
\]
Consequently,
\[
    \mathbf z_1(\chi+\mathbf i0,y)
    =
    r_+(s)-\mathbf i c_s\sqrt a+O(a),
\]
and therefore
\[
    \Im\mathbf z_1(\chi+\mathbf i0,y)<0
\]
for \(\chi<\chi_+(s)\) sufficiently close to \(\chi_+(s)\).

It remains to propagate this conclusion across the whole interval
$(\chi_-(s),\chi_+(s))$. The discriminant is nonzero throughout this
interval, so every root is simple. Since the
continued branch lies in the lower half-plane near \(\chi_+(s)\),
suppose that it meets the real axis at an interior point
\(\chi_0\).  As \(\chi\to\chi_0\) from the side on which the branch is
nonreal,
\[
    \mathbf z_1(\chi+\mathbf i0,y)
    \quad\text{and}\quad
    \overline{\mathbf z_1(\chi+\mathbf i0,y)}
\]
are two distinct conjugate roots that converge to the same real root.
Hence \(P_{\chi_0,s}\) has a multiple root, and therefore
\[
    \Delta_s(\chi_0)=0,
\]
a contradiction.  Thus the branch cannot meet the real axis in the
interior of the interval and remains in the lower half-plane
throughout.
\end{proof}

Define the horizontal slices
\[
    \mathcal R_y
    :=
    \{x\in\mathbb R:(x,y)\in\mathcal R\},
    \qquad
    \mathcal L_y
    :=
    \{x\in\mathbb R:(x,y)\in\mathcal L\},
\]
and let
\[
    J_s:=\bigl(\chi_-(s),\chi_+(s)\bigr).
\]
For every $\chi\in J_s$, the lower-half-plane root identified in the
auxiliary lemma is simple. Hence it admits a holomorphic local
continuation in $w$, and
\[
    \frac{\mathbf z_1(w,y)-1+s}{s}
\]
is nonzero in a sufficiently small complex neighborhood. The
logarithmic branch selected in
Proposition~\ref{prop:stieltjes-root} therefore also admits a local
holomorphic continuation.

For $w>\chi_+(s)$, the ratio
\[
    \frac{\mathbf z_1(w,y)-1+s}{s}
\]
is positive real and the selected logarithm is real. After crossing
$\chi_+(s)$ through the upper half of the $w$-plane, the root lies in
the lower half-plane, so the imaginary part of the continued logarithm
belongs to $(-\pi,0)$ near the right endpoint. Since the root remains
in the lower half-plane on the connected interval $J_s$, the ratio
above never meets the real axis there; consequently the same
logarithmic branch has imaginary part in $(-\pi,0)$ throughout $J_s$.
Proposition~\ref{prop:density-existence} then gives
\[
    0<f_y(\chi)<1,
    \qquad \chi\in J_s.
\]
Thus every point of
\[
    J_s\cap\operatorname{int}_{\mathbb R}(\mathcal R_y)
\]
is a regular liquid point.

Conversely, for
\[
    \chi\in(0,\chi_-(s))\cup(\chi_+(s),3),
\]
all roots of $U_s(z)=\chi$ are real. By
Lemma~\ref{lem:liquid-nonreal-root}, no such point can be a regular
liquid point. Consequently, whenever
\[
    [\chi_-(s),\chi_+(s)]
    \subset
    \operatorname{int}_{\mathbb R}(\mathcal R_y),
\]
we have
\[
    \mathcal L_y\cap(0,3)
    =
    (\chi_-(s),\chi_+(s)),
\]
and the two endpoints are the corresponding internal horizontal
frozen-boundary points.

\begin{figure}[t]
    \centering
    \includegraphics[width=.8\textwidth]{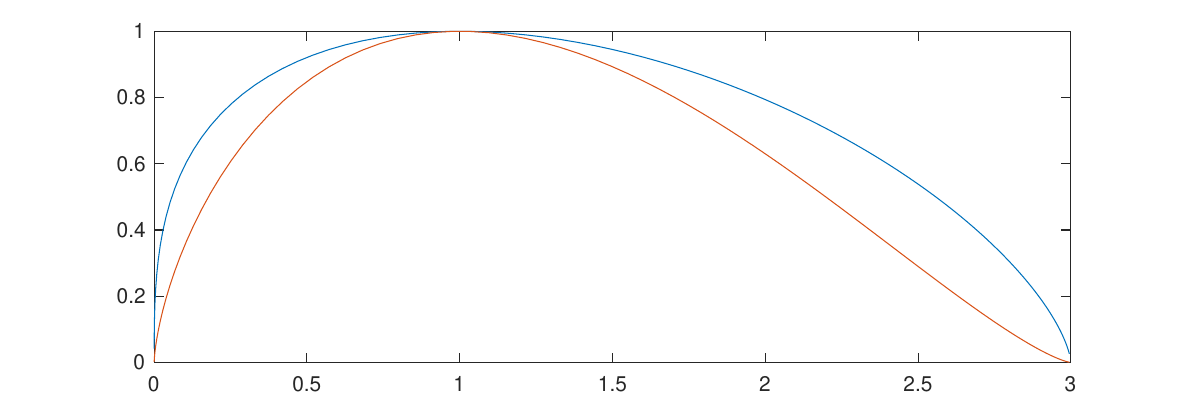}
    \caption{Internal frozen-boundary curves for the scaling limit of
    weighted non-intersecting paths in
    Example~\ref{ex:linear-bottom-boundary}. The blue curve corresponds
    to uniform weights, for which $s=1-y$; the red curve corresponds to
    the relation $y=(1-s)^2$. In both cases the two branches are
    parametrized by $z\in(0,1)\cup(1,\infty)$.}
    \label{fig:fb1}
\end{figure}

We finally compute these endpoints explicitly. The nonzero critical
points satisfy
\[
    (1-s)(z^3+2)=3z.
\]
Combining this relation with $U_s(z)=\chi$ gives
\begin{equation}
    \chi
    =
    \frac{3z^3}{z^3+2},
    \qquad
    s
    =
    \frac{z^3-3z+2}{z^3+2}.
\end{equation}
The two internal frozen-boundary branches displayed in
Figure~\ref{fig:fb1} correspond precisely to
\[
    z\in(0,1)
    \qquad\text{and}\qquad
    z\in(1,\infty).
\]
The first interval parametrizes the branch with $0<\chi<1$, and the
second parametrizes the branch with $1<\chi<3$. The critical point
$z=0$ gives the outer locus $\chi=0$, while
\[
    z\in(-\infty,-2)
\]
gives an additional algebraic double-root branch with $\chi>3$;
neither belongs to the internal frozen boundary shown in
Figure~\ref{fig:fb1}.

In the two weight regimes used in the figure, the parametrizations in
$(\chi,y)$-coordinates are therefore as follows.
\begin{enumerate}
    \item If the weights $x_0,x_1,\ldots,x_{t-1}$ are all equal and
    $\alpha=1$, then $s=1-y$, and
    \begin{equation}
        \chi(z)
        =
        \frac{3z^3}{z^3+2},
        \qquad
        y(z)
        =
        \frac{3z}{z^3+2},
        \qquad
        z\in(0,1)\cup(1,\infty).
    \end{equation}
    This is the blue curve in Figure~\ref{fig:fb1}.

    \item If $\alpha=1$ and $y=(1-s)^2$, then
    \begin{equation}
        \chi(z)
        =
        \frac{3z^3}{z^3+2},
        \qquad
        y(z)
        =
        \left(\frac{3z}{z^3+2}\right)^2,
        \qquad
        z\in(0,1)\cup(1,\infty).
    \end{equation}
    This is the red curve in Figure~\ref{fig:fb1}.
\end{enumerate}
\end{example}

\bigskip

\bigskip
\noindent\textbf{Acknowledgements.} We thank Sylvie Corteel for asking the questions solved in the paper and for helpful discussions. We thank an anonymous referee for a careful reading and detailed comments on an earlier version of the manuscript. ZL acknowledges support from the National Science Foundation under grant 1608896 and from Simons Foundation under grant 638143. DK and IP are grateful to the Workshop on ‘Randomness, Integrability, and Universality’, held in Spring 2022 at the Galileo Galilei Institute for Theoretical Physics, for hospitality and support at some stage of this work. DK acknowledges support from NSF RTG grant DMS-1937241. IP acknowledges support from the Research Council of Finland under grant 365297.

\medskip
Declaration of competing interest.
The authors declare that they have no known competing
financial interests or personal relationships that could have appeared to
influence the work reported in this paper.

\medskip
Data availability.
No data were used for the research described in this article.

\medskip
\noindent\textbf{Declaration of generative AI and AI-assisted technologies in the manuscript preparation process.}
During the preparation of this work, OpenAI ChatGPT was used to assist with revising and polishing portions of the manuscript, including the abstract and introduction, and with systematically checking the consistency of notation, cross-references, and definitions throughout the manuscript. After using these tools, the authors reviewed and edited the content as needed. The authors take full responsibility for the content of the manuscript.

 
\bibliography{lht}
\bibliographystyle{plain}

\end{document}